\documentclass[letterpaper, paper,11pt]{AAS}
\usepackage{subcaption}
\usepackage{bm}
\usepackage{amsmath}
\usepackage{amsfonts}
\usepackage[colorlinks=true, pdfstartview=FitV, linkcolor=black, citecolor= black, urlcolor= black]{hyperref}
\usepackage{overcite}
\usepackage{footnpag}
\usepackage{makecell}
\usepackage{multirow}
\usepackage{fancyhdr}

\fancypagestyle{firstpage}{
    \fancyhf{}
    \fancyhead[L]{%
        \footnotesize
        Presented at the 2026 AAS/AIAA Astrodynamics Specialist Conference,\\
        Whistler, BC, July 26--30, 2026.
    }
    
}

\PaperNumber{26-880}

\begin{document}

\title{Many-Revolution Low-Thrust Transfers within Periodic Orbit Families}

\author{Ian M. Down\thanks{Astrodynamics/Satellite Navigation Engineer, Advanced Space, LLC, 1400 W 122nd Ave, Westminster, CO 80234}}

\maketitle{} 
\thispagestyle{firstpage}

\begin{abstract}
This paper explores the use of a reduced manifold space initialization method for the generation of both time and fuel optimal many-revolution, low-thrust transfers between three-body periodic orbits within a family. Libration point orbit families are approximated by multi-segment Chebyshev polynomials and 1-dimensional Fourier series. Control is mapped into the family space, and a 1-dimensional targeting scheme in the family space dynamics is used to generate an initial guess for the phase space, directly accounting for winding. An integral collocation-based direct method is then used to first continue a feasible trajectory into the phase space, and then optimize for a particular cost function. The methodology is first applied to the distant retrograde orbit family. Modifications for time-regularized dynamics are then discussed, with an ensuing example of many-revolution transfers in the $L_2$ halo family.
\end{abstract}

\section{Introduction}

Many-revolution, low-thrust trajectory optimization has traditionally been developed within the context of planetocentric, two-body dynamics, often leveraging conic section geometry. Early work established the feasibility of optimizing very-low-thrust trajectories over long durations, while highlighting the numerical challenges associated with multi-revolution transfers and sensitivity to initial guesses\cite{scheel1994optimization}. Since then, significant effort has focused on improving convergence and computational efficiency for this class of problems. Several methods have been introduced to better handle the oscillatory nature and long time horizons of many-revolution transfers. Aziz et al.\cite{aziz2018low} employ differential dynamic programming with a Sundman transformation to regularize trajectory parameterization, enabling more efficient optimization. Shannon et al.\cite{shannon2020q} incorporate Q-law guidance into direct optimization to generate high-quality initial guesses and improve solver robustness. More broadly, Ghosh\cite{ghosh2019survey} surveys some of these approaches, emphasizing the trade-offs between robustness and optimality across direct, indirect, and hybrid methods. Together, these works establish a mature framework for two-body, many-revolution, low-thrust trajectory design.

As mission design has expanded to multi-body environments, attention has shifted toward leveraging the dynamical structure of systems such as the circular restricted three-body problem (CR3BP). The low-thrust trajectory solution space has much richer geometry when compare to the perturbed two-body regime, which makes solution generation more difficult. In these settings, invariant structures like periodic orbit families and their manifolds provide new mechanisms for low-thrust trajectory design. Pritchett et al.\cite{pritchett2018impulsive} examine both impulsive and low-thrust transfers between stable and nearly stable periodic orbits using collocation, demonstrating how invariant manifolds can reduce control effort. Complementing this, Pino et al.\cite{pino2020energy} introduce an energy-informed adaptive framework for cislunar trajectory design, using system energy to guide optimization in highly nonlinear regimes. Successive convex optimization incorporating intermediate linking orbits has been investigated by several authors\cite{kayama2022low,hiraiwa2023halo,bernardini2024state}. A similar approach is taken by Aziz et al. using hybrid differential dynamic programming\cite{aziz2019hybrid}. Invariant manifolds are directly parametrized for trajectory optimization using orthogonal polynomials by Kelly et al.\cite{kelly2023orthogonal}. This concept is most recently applied by Fanger and Russell\cite{fanger2026optimal} to extend trajectory design from individual cislunar periodic orbits to entire orbit families, highlighting the added flexibility and complexity introduced when the design space is expanded.

Recent work by Down et al.\cite{down2025manifold,down2026manifold} has explored solving trajectory optimization problems in reduced manifold spaces within the context of minimal torus deviation during quasi-periodic orbit rephasing. The concept of leveraging a reduced manifold space for trajectory design and optimization is now applied to many-revolution, low-thrust transfers over a periodic orbit family. This particular type of transfer is motivated by ride share opportunities as cislunar transport becomes more frequent as well as constellation deployment and reconfiguration. Orbit families as 2-dimensional manifolds are parametrized by multiple segments of Chebyshev polynomials and 1-dimensional Fourier series in the circular restricted three-body problem (CR3BP). Approximation orders are selected by the evaluation of an invariance error metric. Control is mapped into the family space, and a 1-dimensional targeting scheme in the family space dynamics is used to generated an initial guess, directly estimating the number of revolutions during the transfer. An integral collocation-based direct method\cite{ahrens2025conf,ahrens2026direct} is then used to first continue a feasible trajectory into the phase space, and then optimize for a particular cost function. Both minimum time and minimum fuel time free trajectories with both free and fixed arrival phases are initialized using the family space method. Modifications to the mathematical framework are presented for time-regularized dynamics in order to handle sub-family regimes with varying fast/slow dynamics (e.g. NRHOs). The initial guess scheme is shown to be robust for the distant retrograde and halo orbit families exemplified in this paper.

\section{Periodic Orbit Families}
The CR3BP contains an abundance of orbit families because of its single energy integral \cite{howell2001families}. These families are linked together through stability bifurcations, and when examined from a global perspective, establish a useful structure to an otherwise chaotic solution space. When forcing frequencies perturb the three body problem, families of periodic orbits disappear. However, their utility in discovering and analyzing useful trajectories that can then be carried into higher fidelity, full ephemeris models makes their exploration and consideration critical in multi-body mission design. The CR3BP dynamics with acceleration level control parametrized by spherical thrust angles and throttling are given as
\begin{align}
    \ddot{x} - 2\dot{y} - x &= -\frac{(1-{\mu})(x+{\mu})}{d_1^3} - \frac{{\mu}(x-1+{\mu})}{d_2^3} +u_{\max}\,\delta\cos\phi\cos\psi\label{eq:cr3bpcnt}\\
    \ddot{y} + 2\dot{x} - y &= -\frac{(1-{\mu})y}{d_1^3} - \frac{{\mu} y}{d_2^3} +u_{\max}\,\delta\sin\phi\cos\psi\\
    \ddot{z} &= -\frac{(1-{\mu})z}{d_1^3} - \frac{{\mu} z}{d_2^3}+u_{\max}\,\delta\sin\psi
\end{align}
\noindent where $d_1 = \sqrt{(x+{\mu})^2 + y^2 + z^2}$, $d_2 = \sqrt{(x-1+{\mu})^2 + y^2 + z^2}$, and $\delta\in[0,\,1]$ is the throttle parameter.

\subsection{The Family Function and Invariance}
Let an autonomous dynamical system with continuous affine control input be defined as
\begin{equation}
    \dot{\boldsymbol{x}}=\boldsymbol{f}(\boldsymbol{x}) + \boldsymbol{g}(\boldsymbol{x})\boldsymbol{u},\quad\boldsymbol{x}\in\mathbb{R}^n
\end{equation}
A family of periodic orbits embedded within this system in the absence of control is defined here by its family function $\boldsymbol{\Gamma}$. For a 1-parameter family of 1-dimensional tori, $\boldsymbol{\Gamma}$ is a function of the family's parameter $p$, and the torus angle of a particular periodic orbit $\theta$, such that the resulting topological space is the cylinder.
\begin{equation}
    \boldsymbol{\Gamma}(\boldsymbol{y}):= \left(\mathbb{R}^1\times\mathbb{T}^1\right)\rightarrow\mathbb{R}^n,\quad
    \boldsymbol{y}=[p,\:\theta]^\text{T}\quad p\in[p_{\text{min}},p_{\text{max}}],\quad \theta\in[0,2\pi)
\end{equation}
The parameter $p$ must be strictly monotonic over the defined family set such that each orbit within the family is uniquely defined by this parameter. Examples of possible parameters for families in the CR3BP include perilune distance, Jacobi's Constant, orbital period/frequency, etc.

A family function is invariant with respect to the system's natural dynamics if it obeys the following field invariance condition in general form.
\begin{equation}
    \left[\frac{\partial \boldsymbol{\Gamma}}{\partial \boldsymbol{y}}\right]\dot{\boldsymbol{y}} =\boldsymbol{f}\left(\boldsymbol{\Gamma}(\boldsymbol{y})\right)\label{eq:faminv}
\end{equation}
The left hand side of Eq. \ref{eq:faminv} can be viewed as the family manifold's tangent bundle as a subset of the system's dynamics. The family parameter describing each orbit does not naturally change with respect to time, i.e. $\dot{p}=0$. If the torus function for each orbit is parametrized by its constant, fundamental frequency, it is said to be the fundamental torus function such that $\dot{\theta}(p,\theta)=\omega(p)$ depends only on which orbit is selected. This gives the following time explicit relationship in angle evolution.
\begin{equation}
    \dot{\boldsymbol{y}}=\boldsymbol{l}(\boldsymbol{y}) =
    \begin{bmatrix}
        0 \\ \omega(p)
    \end{bmatrix},\quad \theta(p,t) = \omega(p)(t-t_0) + \theta(t_0)\label{eq:angintime}
\end{equation}
The reduced family invariance condition then becomes
\begin{equation}
    \left[\frac{\partial \boldsymbol{\Gamma}}{\partial \theta}\right]\omega(p)=\boldsymbol{f}\left(\boldsymbol{\Gamma}(\boldsymbol{y})\right)
\end{equation}

\subsection{Approximating the Family Function}
Closed form solutions to family functions in nonlinear systems do not exist, generally speaking. In order to leverage a reduced family space manifold to initialize optimal transport solutions, an orbit family's function must be approximated in closed form. This was done by Hirani and Russell\cite{hirani2006approximations} by using Fourier series and Chebyshev polynomials for distant retrograde orbits, and more recently expanded on by Bourjeili and Russell\cite{bourjeili2025approximations} to include regularization. A similar approach is taken in this work. The family parameter space is represented using Chebyshev polynomials of the first kind at Chebyshev-Gauss-Lobatto (CGL) nodes, which include domain end points. The family torus angle is represented by a discrete Fourier series, as the domain is naturally $2\pi$ periodic. Multiple segments of Chebyshev polynomials are used to mitigate the nonlinear behavior of the family parameter.

\subsubsection{Chebyshev Polynomials}
A function $f$ can be approximated by a set of orthogonal basis functions and coefficients as
\begin{gather}
    f(\tau)\approx\sum_{i=0}^m\varphi_i(\tau)c_i = \boldsymbol{\varphi}_m(\tau)\boldsymbol{c} \\
    \boldsymbol{\varphi}_m(\tau)=
    \begin{bmatrix}
        \varphi_0(\tau) & \cdots & \varphi_m(\tau)
    \end{bmatrix}\in\mathbb{R}^{1\times(m+1)},\quad
    \boldsymbol{c}=
    \begin{bmatrix}
        c_0 \\ \vdots \\ c_m
    \end{bmatrix}\in\mathbb{R}^{(m+1)\times1}
\end{gather}

Selecting Chebyshev polynomials as the basis presents many advantages in functional approximation\cite{trefethen2019approximation}. By selecting this basis, the maximum approximation error is minimized. Furthermore, careful selection of interpolation nodes (polynomial roots or extrema) leads to an interpolation matrix with full rank and avoids Runge's or Gibb's phenomena. An $m$th order Chebyshev polynomial of the first kind $\varphi_m(\tau) = \cos{(m\arccos\tau)}$ with $\tau\in[-1,1]$ is defined by its generating function:
\begin{align}
    \varphi_0(\tau) &= 1 \\
    \varphi_1(\tau) &= \tau \\
    \varphi_{i}(\tau) &= 2\tau\varphi_{i-1}(\tau)-\varphi_{i-2}(\tau),\quad i\geq2
\end{align}
The ascending CGL nodes are then defined to be the $m$th order polynomial's extrema as
\begin{equation}
    \tau_j = -\cos\left(\frac{j\pi}{m}\right),\quad j=[0,1,\cdots,m]
\end{equation}

Approximating the scalar function $f$ by a set of coefficients $\boldsymbol{c}$ can then be re-written in matrix form as
\begin{align}
    \begin{bmatrix}
        f(\tau_0) \\ \vdots \\ f(\tau_m)
    \end{bmatrix} &= 
    \begin{bmatrix}
        \varphi_0(\tau_0) & \cdots & \varphi_m(\tau_0) \\
        \vdots & \ddots & \vdots \\
        \varphi_m(\tau_0) & \cdots & \varphi_m(\tau_m)
    \end{bmatrix}
    \begin{bmatrix}
        c_0 \\ \vdots \\ c_m
    \end{bmatrix} \\
    \boldsymbol{f}(\boldsymbol{\tau}) &= \boldsymbol{\Phi}_m(\boldsymbol{\tau})\boldsymbol{c} \implies \boldsymbol{c} = \boldsymbol{\Phi}_m^{-1}(\boldsymbol{\tau})\boldsymbol{f}(\boldsymbol{\tau}),\quad\boldsymbol{c}\in\mathbb{R}^{m+1}
\end{align}
There is no analytical inverse of $\boldsymbol{\Phi}$, but explicit computation is only done once, and does not impede algorithmic efficiency in methods that rely on these fits.

In general, a function's domain with independent variable $p$, like a family function, does not align with the Chebyshev domain using the independent variable $\tau$. To transform between domains, the following expressions are used.
\begin{equation}
    \tau = \left(\frac{2}{p_{\text{max}}-p_{\text{min}}}\right)p - \left(\frac{p_{\text{max}}+p_{\text{min}}}{p_{\text{max}}-p_{\text{min}}}\right)\quad\Leftrightarrow\quad p=\left(\frac{p_{\text{max}}-p_{\text{min}}}{2}\right)\tau + \left(\frac{p_{\text{max}}+p_{\text{min}}}{2}\right)
\end{equation}
In the above fit of $f$, data is collected at all $p_j$ corresponding to $\tau_j$ of the CGL nodes, and placed into matrix form.

\subsubsection{Discrete Fourier Series}
The discrete Fourier transform is often used in the computation and approximation of higher dimensional tori \cite{down2026linear}. A uniform set of angles $\boldsymbol{\theta}$ is defined with a corresponding mode vector $\boldsymbol{k}$
\begin{align}
    \boldsymbol{\theta}&=[\theta_1,\cdots,\theta_m],\quad\theta_j = \frac{2\pi(j-1)}{m},\quad j=[1,\cdots,m] \\
    \boldsymbol{k} &= 
    \begin{bmatrix}
        -\frac{m-1}{2} & \cdots & -1 & 0 & 1 & \cdots & \frac{m-1}{2}
    \end{bmatrix} \quad \text{for} \,\, \text{odd} \,\, m \\
    \boldsymbol{k} &=
    \begin{bmatrix}
        -\frac{m}{2} & \cdots & -1 & 0 & 1 & \cdots & \left(\frac{m}{2}-1\right)
    \end{bmatrix} \quad \text{for} \,\, \text{even} \,\, m
\end{align}
with the discrete Fourier transform then expressed in matrix form as
\begin{equation}
    \boldsymbol{D} = \frac{1}{m}\text{e}^{-\sqrt{-1}\cdot\boldsymbol{k}^\text{T}\boldsymbol{\theta}}\in\mathbb{R}^{m\times m}
\end{equation}
The final approximation of the scalar function $f$ is written in the following matrix form.
\begin{equation}
    \boldsymbol{f}(\boldsymbol{\theta}) = \boldsymbol{D}^{-1}(\boldsymbol{\theta})\boldsymbol{c}\implies\boldsymbol{c}=\boldsymbol{D}(\boldsymbol{\theta})\boldsymbol{f}(\boldsymbol{\theta}),\quad \boldsymbol{c}\in\mathbb{R}^{m}
\end{equation}
Unlike the Chebyshev approximation, both $\boldsymbol{D}$ and its inverse are available analytically.

\subsubsection{Family Function Fit and Evaluation}
First, the domain of ${p}$ is broken up into $m_s$ segments. While an \textit{hp}-adaptive method is more robust, a simple hand-tuned power law segmentation is used in this work to handle parameter nonlinearity. The bounds for each segment are
\begin{align}
    p_{k}^- &= p_{\min} + (p_{\max} - p_{\min})\left(\frac{k-1}{m_s}\right)^\beta,\quad \text{for}\,\,\, k=1,\dots,m_s\\
    p_{k}^+ &= p_{\min} + (p_{\max} - p_{\min})\left(\frac{k}{m_s}\right)^\beta\label{eq:nonlinearfit},\quad\beta>0
\end{align}

The family parameter is then nondimensionalized as $\bar{p}=p/p_{\max}$ to remove units. This will be important when constructing initial guesses in the family space. A Chebyshev polynomial of degree $m_p$ is used to approximate each segment using $m_p+1$ periodic orbits evaluated on the mapped CGL grid. For each periodic orbit, an $m_\theta$ degree Fourier series is used such that $m_\theta$ states evenly spaced in time are propagated and collected according to Eq. \ref{eq:angintime}\footnote{To accommodate orbit families with fast and slow dynamics, regularization should be included at this step for an accurate fit requiring a low number of data points.}. The resulting data structure is
\begin{equation}
    \boldsymbol{\mathcal{X}} =
    \begin{bmatrix}
        \boldsymbol{X}_1 \\ \vdots \\ \boldsymbol{X}_{m_s}
    \end{bmatrix}\in\mathbb{R}^{m_s(m_\theta(m_p+1)\times6)},\quad
    \boldsymbol{X}_k = 
    \begin{bmatrix}
        \boldsymbol{x}_{0,1}^\text{T} \\ \vdots \\ \boldsymbol{x}_{0,m_\theta}^\text{T} \\ \boldsymbol{x}_{1,1}^\text{T}\\ \vdots\\ \boldsymbol{x}_{m_p,m_\theta}^\text{T}
    \end{bmatrix}\in\mathbb{R}^{m_\theta(m_p+1)\times6},\quad \boldsymbol{x}_{i,j}=
    \begin{bmatrix}
        x_{i,j} \\ y_{i,j} \\ z_{i,j} \\ \dot{x}_{i,j} \\ \dot{y}_{i,j} \\ \dot{z}_{i,j}
    \end{bmatrix}
\end{equation}

The Kronecker product of the Chebyshev and Fourier series approximation matrices for each segment are stacked to form the linear system
\begin{equation}
    \boldsymbol{\mathcal{X}} = 
    \left[\mathbb{I}_{m_s}\otimes\boldsymbol{\Phi}(\boldsymbol{\tau})\otimes\boldsymbol{D}^{-1}(\boldsymbol{\theta})
    \right]\cdot\boldsymbol{\mathcal{C}},\quad \boldsymbol{\mathcal{C}} =
    \begin{bmatrix}
        \boldsymbol{C}_1 \\ \vdots \\ \boldsymbol{C}_{m_s}
    \end{bmatrix}
    \label{eq:fit_system}
\end{equation}
where $\boldsymbol{\mathcal{C}}$ contains the stacked coefficients of each segment. The above system is square, and thus a unique solution exists given the full rank properties of the interpolation matrices. $C^0$ continuity is automatically enforced by construction in Eq. \ref{eq:nonlinearfit}, but $C^1$ continuity across segments is left open. Thus, the following equations are appended to the linear system that then requires a least-squares solution.
\begin{equation}
    \left.\frac{\partial \boldsymbol{\Gamma}}{\partial \bar{p}}\right|_{\bar{p}_k^{\,+}} - \left.\frac{\partial \boldsymbol{\Gamma}}{\partial \bar{p}}\right|_{\bar{p}_{k+1}^{\,-}} = \mathbf{0},\,\,\forall\,\,\theta
\end{equation}
The sensitivities are computed analytically by equations found in the Appendix. With the coefficient solution, a value $p$ can be dimensionalized and mapped into a particular segment given the segment bounds. The resulting family function can then be evaluated as
\begin{equation}
    \boldsymbol{\Gamma}(\bar{p},\theta) = \boldsymbol{C}^\text{T}\cdot\left[\boldsymbol{\varphi}(\tau(\bar{p}))\otimes\text{e}^{\sqrt{-1}\cdot\boldsymbol{k}\cdot\theta}\right]^\text{T}
\end{equation}
This representation of the family function allows for quick and analytically tractable sensitives with respect to the input arguments. These are given in the Appendix. In order to complete the dynamic representation of the family, a fit of the family's orbital frequency as a function of family parameter $\omega(\bar{p})$ must be computed. This is done using the same structures developed above for the states, with the data matrix collapsing to a 1-dimensional data vector. If the family's orbital frequency is monotonic, it can be used as the family parameter, and this step can be skipped.

\subsubsection{Approximation Error Measure}
The global fit for a particular orbit family is dependent on four parameters: $m_s$, $m_p$, $m_\theta$, and $\beta$. To understand the fit's approximating capabilities, an approximation error measure is needed. For this, the family invariance condition of Eq. \ref{eq:faminv} is used to compare the dynamics vector field to the tangent bundle of the family function at points precisely in between those used to construct the fit. The intermediate uniform angle distribution is
\begin{equation}
    \tilde{\theta}_i=\frac{\pi(2i-1)}{m_\theta},\quad i=1,\dots,m_\theta
\end{equation}
A $m_p$ degree Chebyshev polynomial grid has $m_p+1$ CGL nodes. Its intermediate $m_p$ nodes denoted $\tilde{\tau}_j$ lie on the Chebyshev-Gauss grid of the $m_p-1$ degree Chebyshev polynomial\cite{ahrens2026direct}. Thus, for all segments, the invariance metric becomes
\begin{equation}
    \epsilon = \frac{1}{ m_s m_p m_\theta} \sum_{k=1}^{m_s} \sum_{j=1}^{m_p} \sum_{i=1}^{m_\theta}\left|\left|\left[\frac{\partial \boldsymbol{\Gamma}}{\partial \tilde{\theta}_i}\right]\omega(\bar{p}_{k}(\tilde{\tau}_j)) -\boldsymbol{f}\left(\boldsymbol{\Gamma}(\bar{p}_{k}(\tilde{\tau}_j),\tilde{\theta}_i)\right)\right|\right|\label{eq:inverr}
\end{equation}

\subsection{Distant Retrograde Orbit Family}
\begin{figure}[htpb!]
    \centering
    \begin{subfigure}{0.34\textwidth}
         \centering
         \includegraphics[width=1\textwidth]{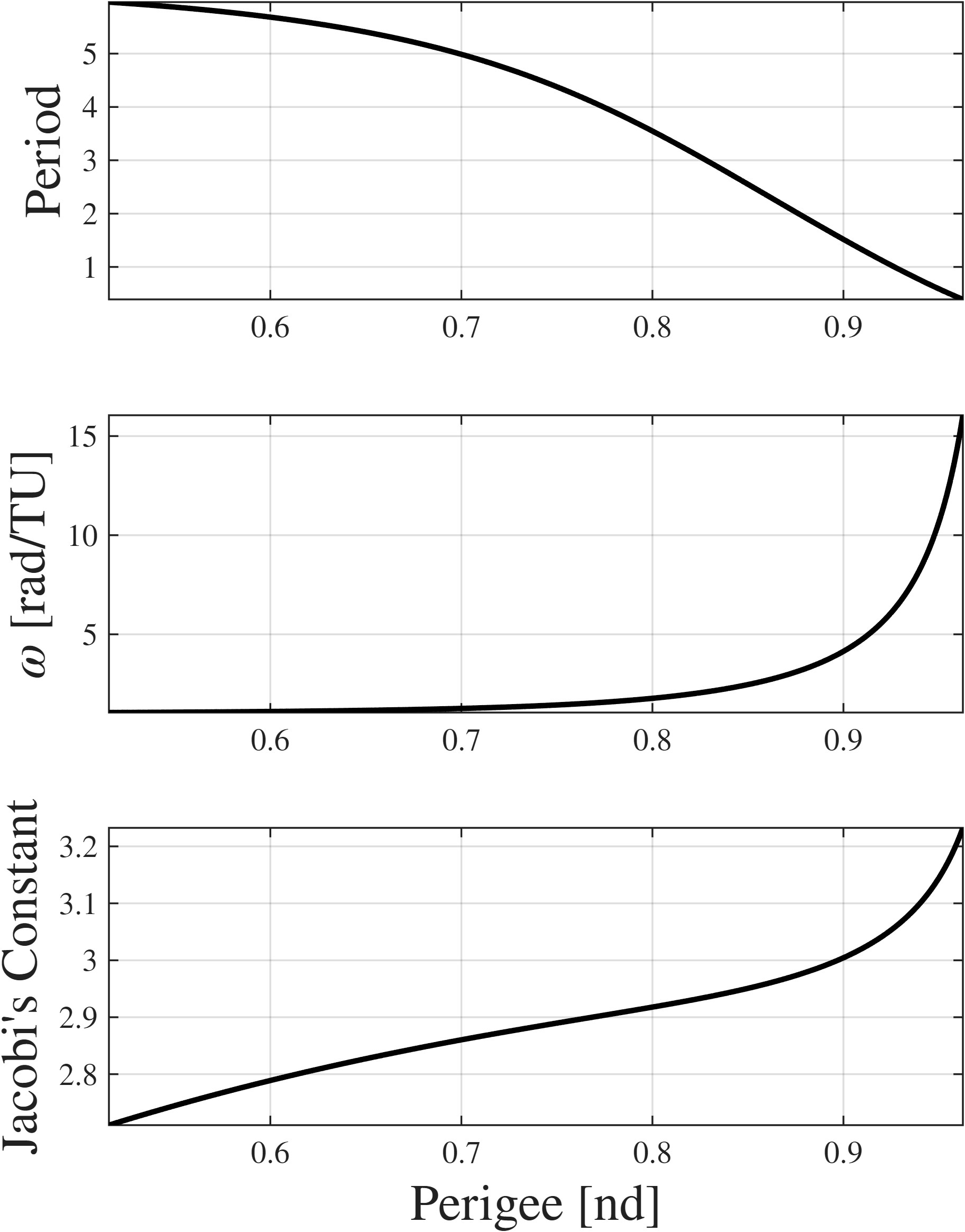}
         \caption{Family Data}
     \end{subfigure}\hfill
     \begin{subfigure}{0.275\textwidth}
         \centering
         \includegraphics[width=1\textwidth]{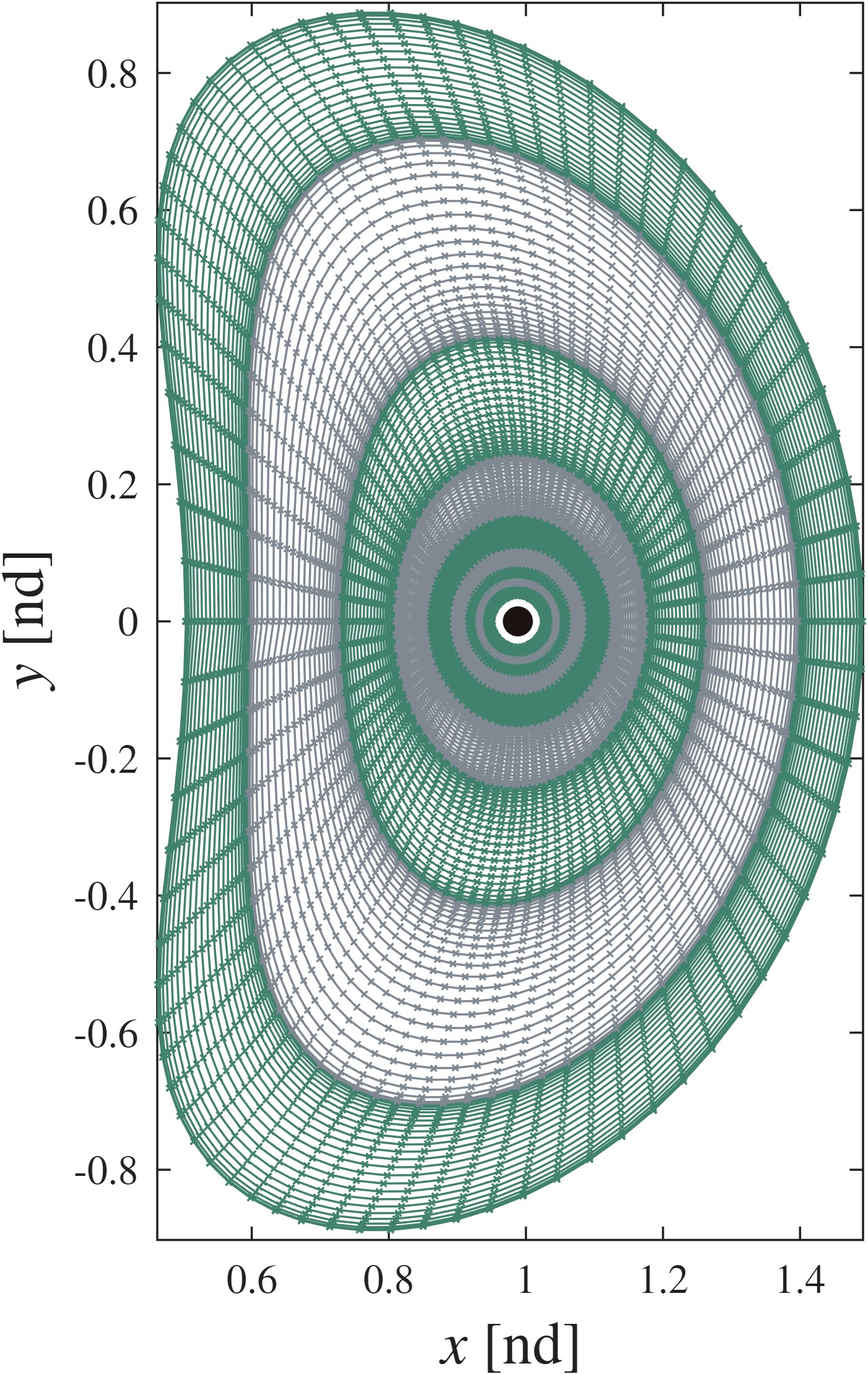}
         \caption{Family Fit}
     \end{subfigure}\hfill
     \begin{subfigure}{0.37\textwidth}
         \centering
         \includegraphics[width=1\textwidth]{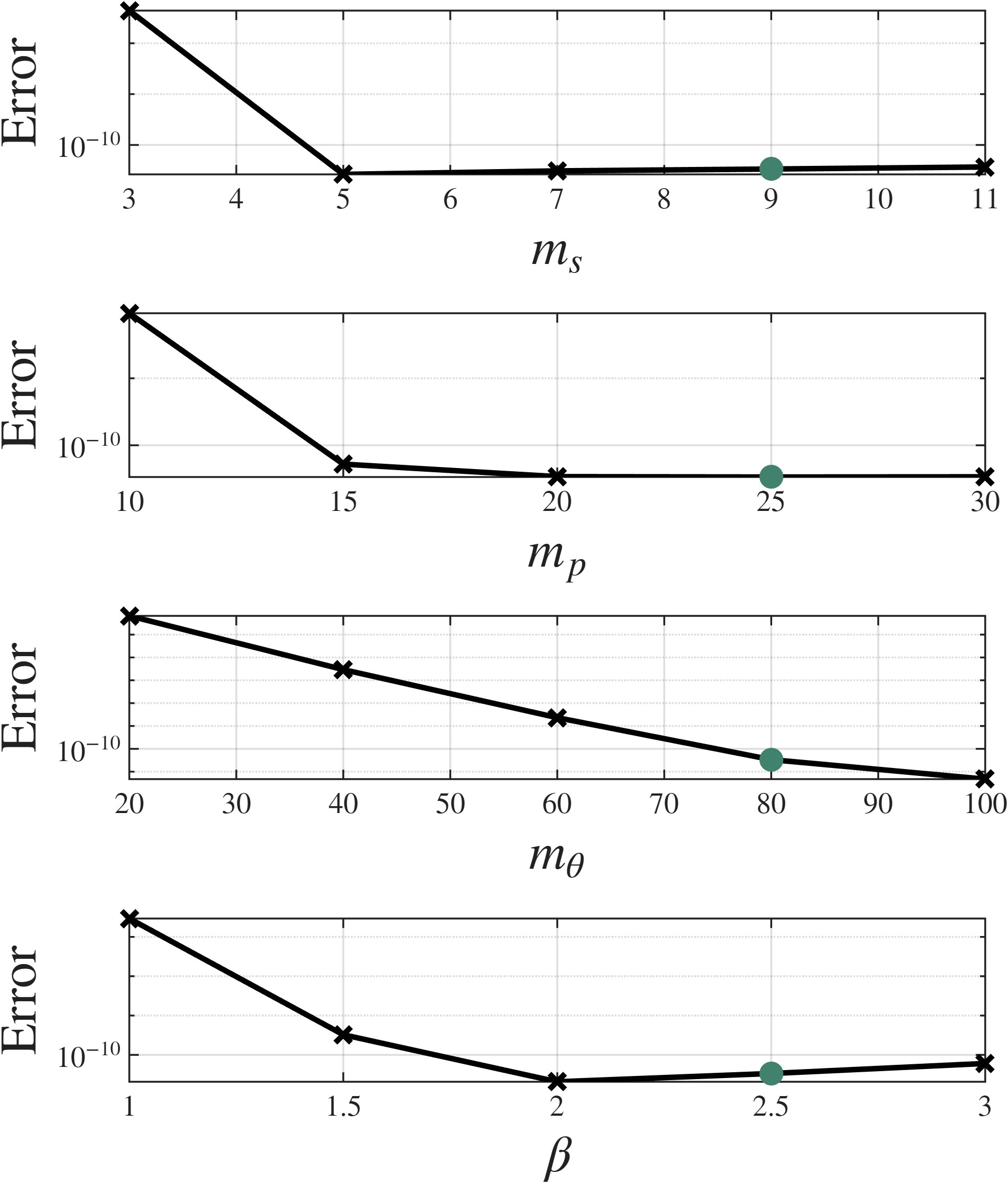}
         \caption{Error Data}
     \end{subfigure}
    \caption{Distant retrograde orbit family: characteristics and fit. Fit shown with $m_s=9$, $m_p=25$, $m_\theta=80$, and $\beta=2.5$.}
    \label{fig:dro_fam_fit}
\end{figure}

The distant retrograde orbit (DRO) family is used as an example throughout this work because time regularization is not required for an accurate fit and orbital frequency can be used as a family parameter. It is also geometrically well-behaved in the rotating frame, with clear, non-Keplerian structures in the Earth-centered and Moon-centered inertial frames. This is evident when looking at low-integer sidereal resonant orbits, shown in the Appendix. Hence, traditional planetocentric approaches to many-revolution, low-thrust transfer strategies that leverage perturbed orbital element spaces may fail. The DRO family also shows promise for many mission design applications, and transfers to this regime from Earth are well studied\cite{capdevila2014various}. A base fit of $m_s=9$, $m_p=25$, $m_\theta=80$, and $\beta=2.5$ are used, and deviations in each fit parameter are examined individually. The results are shown in Fig. \ref{fig:dro_fam_fit}.
Both Jacobi's constant and orbital frequency are monotonic as a function of perigee (also monotonic). A closer look at orbital frequency as a function of perigee reveals the need for either \textit{hp}-adaptive or power law segmentation because of its exponential like growth as family members approach the secondary body. The invariance error metric shown in Fig. \ref{fig:dro_fam_fit} is the averaged quantity from Eq. \ref{eq:inverr}. Despite $\beta=2$ resulting in a lower error measure by this quantity, $\beta=2.5$ was selected for this work because it resulted in a significantly lower absolute maximum error over all error points -- an equally important metric when considering ideal fit parameters (data not shown for brevity).

\section{Initialization in the Family Space}
A challenging aspect of many-revolution, low-thrust transfers is providing an accurate-enough initial guess for both the number of revolutions, and the total time of flight. Two-body dominant dynamics can leverage initialization schemes centered around orbital elements, many of which are described in the introduction. Three-body dynamics lack those luxuries. Even the concept of a revolution in three-body flow can be poorly defined, and with hyperbolic manifolds scattered around, optimal low-thrust solution geometry is nearly almost impossible to "guess." This work approaches the problem two-fold: first, the scope is limited to transfers within a single periodic orbit family motivated by ride share drop offs or constellation reconfiguration/rephasing. Second, the parametrized family space dynamics are directly used to generated an initial guess for minimum time-like solutions that can be mapped into the phase space via the family function. To begin, a connection between control in both spaces must be established.

In expanding the definition of manifold invariance to include control in both the manifold space $\boldsymbol{\nu}$ and the phase space $\boldsymbol{u}$\cite{down2026manifold}, it was previously shown that 
\begin{align}
    \dot{\boldsymbol{y}} &=\boldsymbol{l}(\boldsymbol{y}) + \boldsymbol{\nu}\label{eq:famdyncnt}\\
    \boldsymbol{g}(\boldsymbol{\Gamma}(\boldsymbol{y}))\boldsymbol{u}&=\left[\frac{\partial \boldsymbol{\Gamma}}{\partial \boldsymbol{y}}\right]\boldsymbol{\nu}\label{eq:famspacecntcon}
\end{align}
For any non-zero instantaneous family space control input, the invariance error becomes
\begin{align}
    \delta\dot{\boldsymbol{x}}(\boldsymbol{y},\boldsymbol{\nu}) = \left[\mathbb{I}_n - \boldsymbol{g}(\boldsymbol{\Gamma}(\boldsymbol{y}))\boldsymbol{g}(\boldsymbol{\Gamma}(\boldsymbol{y}))^\text{T}\right]\left[\frac{\partial \boldsymbol{\Gamma}}{\partial \boldsymbol{y}}\right]\boldsymbol{\nu}\label{eq:dynerr}
\end{align}
For a large class of dynamical systems, including astrodynamics models, the affine control map takes the constant form
\begin{equation}
    \boldsymbol{g}(\boldsymbol{x}) = \boldsymbol{B}= \begin{bmatrix}
        \boldsymbol{0}_{\bar{n}\times\bar{n}} & \mathbb{I}_{\bar{n}}
    \end{bmatrix}^\text{T}
\end{equation}
because control authority is not directly available at the velocity level to continuously change the position states. This means the manifold space can never be tracked perfectly under available acceleration level control input. When looking at this through the lens of family space control, this leads to the following dynamical error
\begin{equation}
    \delta\dot{\boldsymbol{x}}(\boldsymbol{y},\boldsymbol{\nu}) = 
    \begin{bmatrix}
        \left[\frac{\partial \boldsymbol{\gamma}}{\partial \boldsymbol{y}}\right]\boldsymbol{\nu}\\
        \boldsymbol{0}_{\bar{n}}
    \end{bmatrix},\quad \boldsymbol{\Gamma} =\begin{bmatrix}
        \boldsymbol{\gamma} \\ \dot{\boldsymbol{\gamma}}\label{eq:reducederror}
    \end{bmatrix}
\end{equation}
with $\bar{n}=n/2$. Previous work looked to minimize this error in the torus space on quasi-periodic invariant tori to find phase space trajectories that stuck close to the torus during low-thrust rephasing maneuvers\cite{down2026manifold}. Solving optimal control problems in the manifold space that don't penalize this mapped velocity control input error (i.e mapped minimum acceleration level fuel/energy problems) lead to erroneous trajectories that blatantly defy phase space physics because velocity discontinuity is allowed (and typically encouraged for optimality on the manifold).

Despite the facts above, the family space can still be leveraged to generate useful initial guess for number of revolutions, state histories, and time of flight to a lesser degree. This is done by examining a minimum time control structure in the family space context. The family dynamics with control of Eq. \ref{eq:famdyncnt} clearly show that the family parameter $p$ only changes when control is applied, and that change is not dependent on $\theta$. Thus, when considering a minimum time, phase-free transfer within the family, the parameter control component should be set to a maximum in the desired family parameter direction. That is
\begin{equation}
    \boldsymbol{\nu}=
    \begin{bmatrix}
        \text{sign}(p_f-p_0)\,\nu_{\max} \\ 0
    \end{bmatrix} \label{eq:famcnt}
\end{equation}
The quantity $\nu_{\max}$ is so far undefined and arbitrary. It can be given dynamical significance by connecting it to $u_{\max}$ from the phase space by using the relationship of Eq. \ref{eq:famspacecntcon}. Let
\begin{equation}
    \boldsymbol{T}=\left[\frac{\partial\boldsymbol{\gamma}}{\partial\boldsymbol{y}}\right],\quad\dot{\boldsymbol{T}}=\left[\frac{\partial\dot{\boldsymbol{\gamma}}}{\partial\boldsymbol{y}}\right]
\end{equation}
then with the control direction clearly defined by Eq. \ref{eq:famcnt}, the maximum family space control magnitude is mapped as
\begin{equation}
    \nu_{\max}(\boldsymbol{y})=\frac{\alpha(\boldsymbol{y})u_{\max}}{\sqrt{\hat{\boldsymbol{\nu}}^\text{T}\dot{\boldsymbol{T}}^\text{T}\dot{\boldsymbol{T}}\hat{\boldsymbol{\nu}}}}
\end{equation}
Because of the non-zero error in Eq. \ref{eq:reducederror}, some penalty $\alpha(\boldsymbol{y})$ related to this quantity to account for the extra control authority needed to handle the family manifold's departure when operating in the phase space. The following scaling factor applies that penalty via the average singular value of the reduced error mapping
\begin{equation}
    \alpha(\boldsymbol{y})=\frac{1}{1+\frac{1}{2}\sqrt{\text{tr}(\boldsymbol{T}^\text{T}\boldsymbol{T})}}<1
\end{equation}
which can be meaningfully applied under the assumption that the family parameter's dimensions have been removed. Using the maximum singular value of $\boldsymbol{T}$ in the penalty function proved to be too conservative, while not including this penalty function resulted in projected phase space trajectories that showed disproportionate amounts of control effort.

With the minimum time, phase-free control structure fully defined in the family space, a numerical targeting problem can be constructed and solved to generate an initial guess for seeding a phase space trajectory optimization method. The targeting problem is written as
\begin{equation}
    \text{find}\:t_f:\:\: p_0 + \int_{t_0}^{t_f}\text{sign}(p_f-p_0)\,\nu_{\max}(\boldsymbol{y})\,\text{d}t=p_f
\end{equation}
This can be solved simply with an event function in any ordinary differential equation propagator. With that solution comes a final value $\theta(t_f)$ on the real line, which is used to estimate the number of revolutions
\begin{equation}
    N_{\text{rev}}=\text{round}\left(\frac{\theta(t_f)-\theta_0}{2\pi}\right)
\end{equation}
Finally, the resulting family space time history is mapped to the phase space to arrive at an initial state history.
\begin{equation}
    \boldsymbol{x}(t) = \boldsymbol{\Gamma}(p(t),\theta(t))
\end{equation}

The quantities above can now be used to seed any direct trajectory optimization method such as multiple shooting or collocation. It's possible an indirect method could be initialized effectively as well, but more curated attention would have to paid to mapping costates in a meaningful way. A direct integral collocation framework discussed in the next section is used in this work, with the predicted number of revolutions guiding the segmentation count. In light of the maximum family space control, both minimum time and minimum fuel phase space control variables upon transition are initialized with maximum magnitude, using directions that align with the mapped phase space's instantaneous velocity vector. A phase-free feasibility problem is solved first with zero cost and the family function as the terminal manifold. Then the phase-free time or fuel optimal problem is solved as a direct continuation. Any phase-fixed or time-dependent phase problem is solved from its phase-free solution. Because the family space initialization scheme doesn't account for the variational dynamics surrounding the family, large family space initialization versus optimal phase space solution differences are expected in families with significant hyperbolicity. This makes this method suited more for stable or nearly stable orbit families. While a constant control law was used above, a Lyapunov control strategy could similarly be applied in the family space\cite{peterson2020lyapunov}.

\begin{figure}[htbp!]
     \centering
     \begin{subfigure}{0.24\textwidth}
         \centering
         \includegraphics[width=1\textwidth]{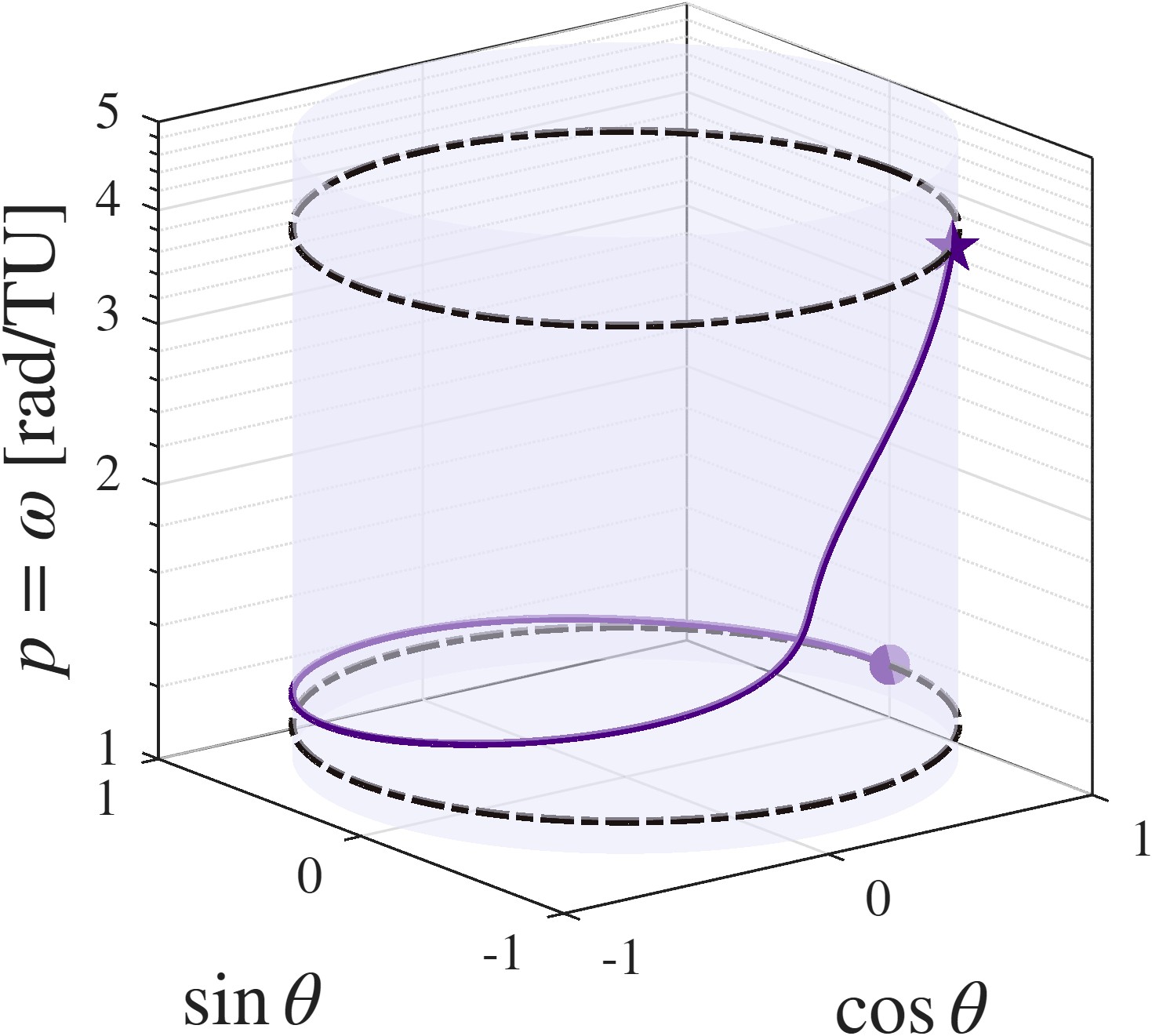}\\
         \includegraphics[width=1\textwidth]{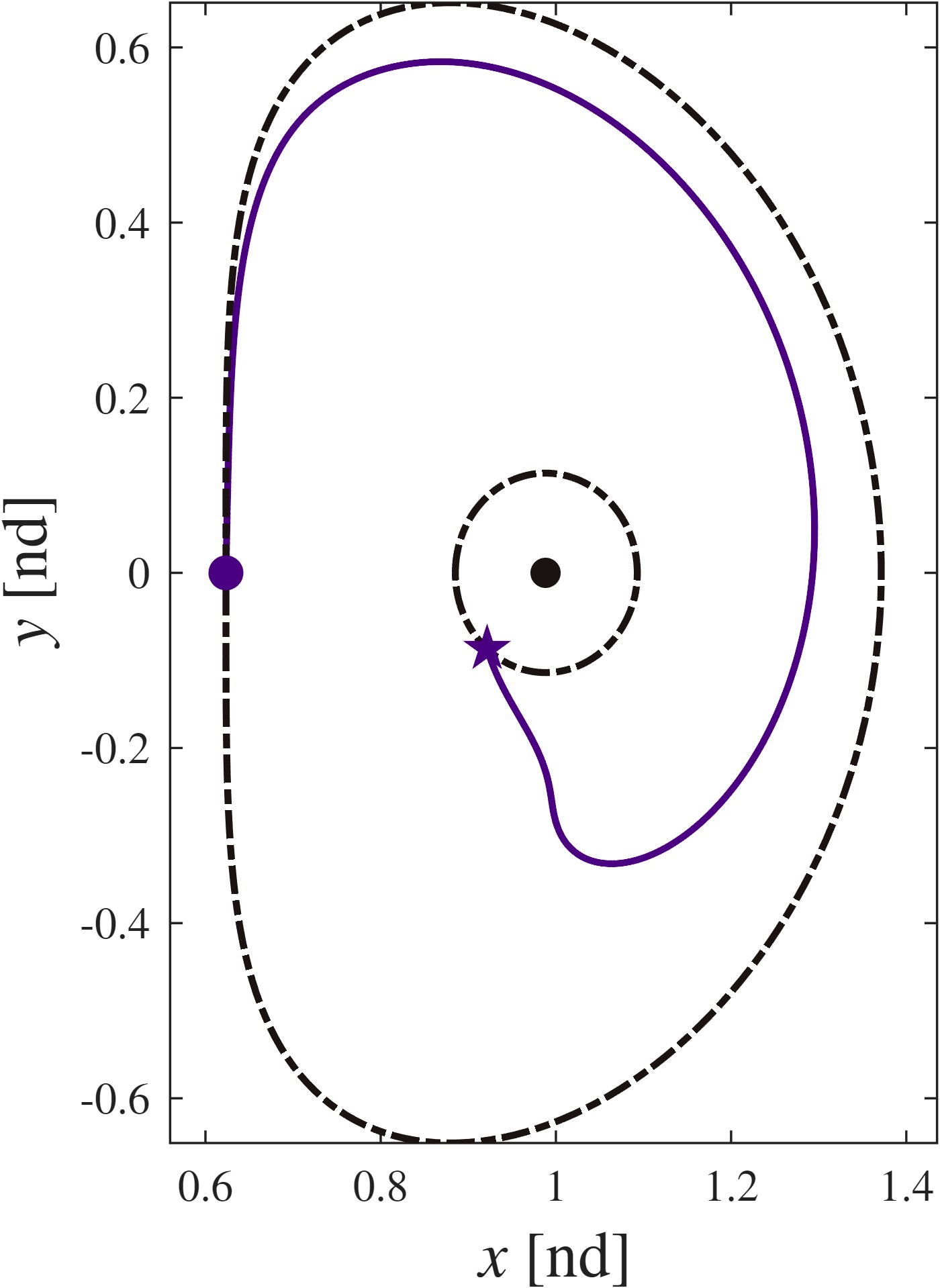}
         \caption{$u_{\max}=0.3$}
     \end{subfigure}\hfill
     \begin{subfigure}{0.24\textwidth}
         \centering
         \includegraphics[width=1\textwidth]{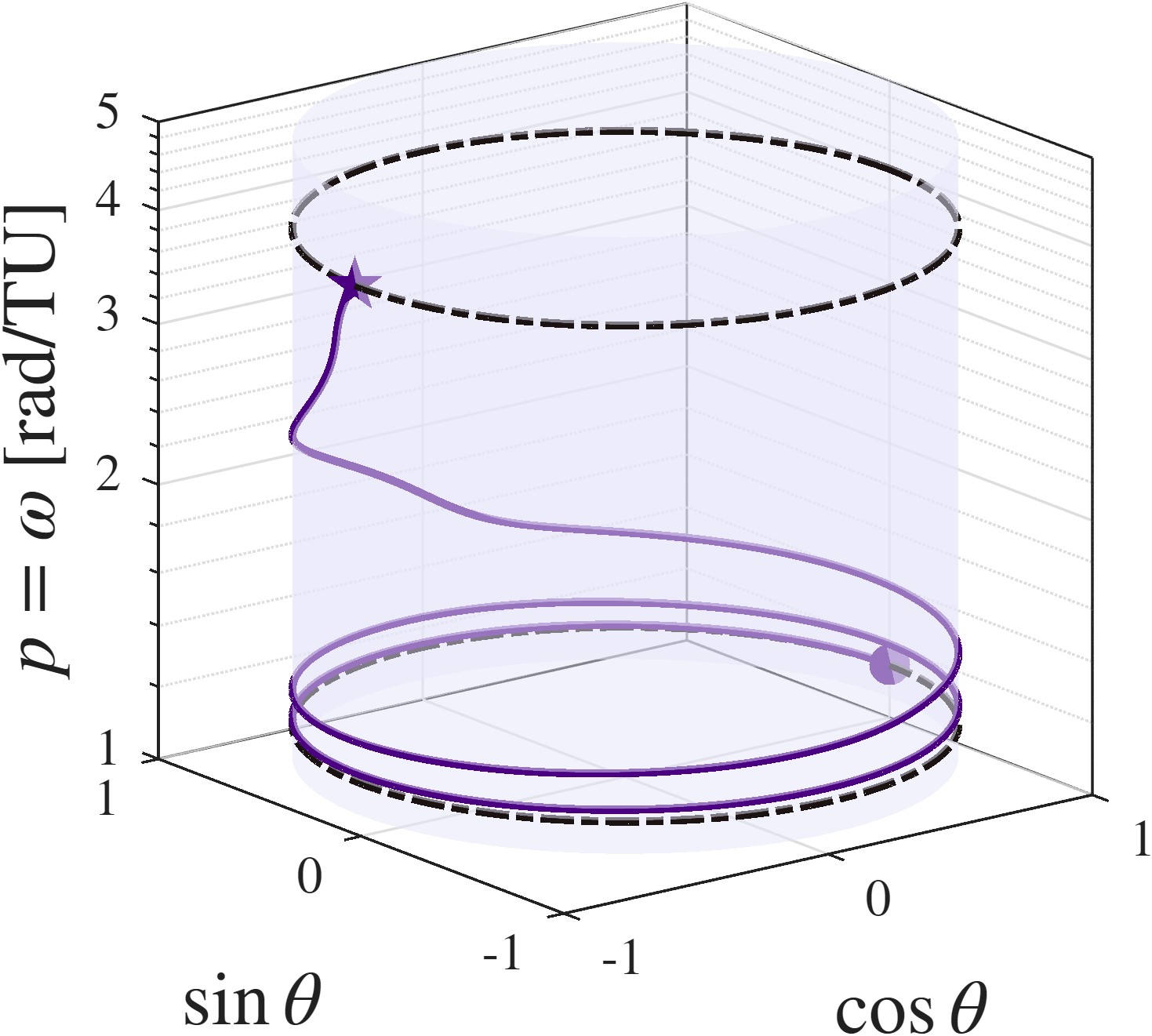}\\
         \includegraphics[width=1\textwidth]{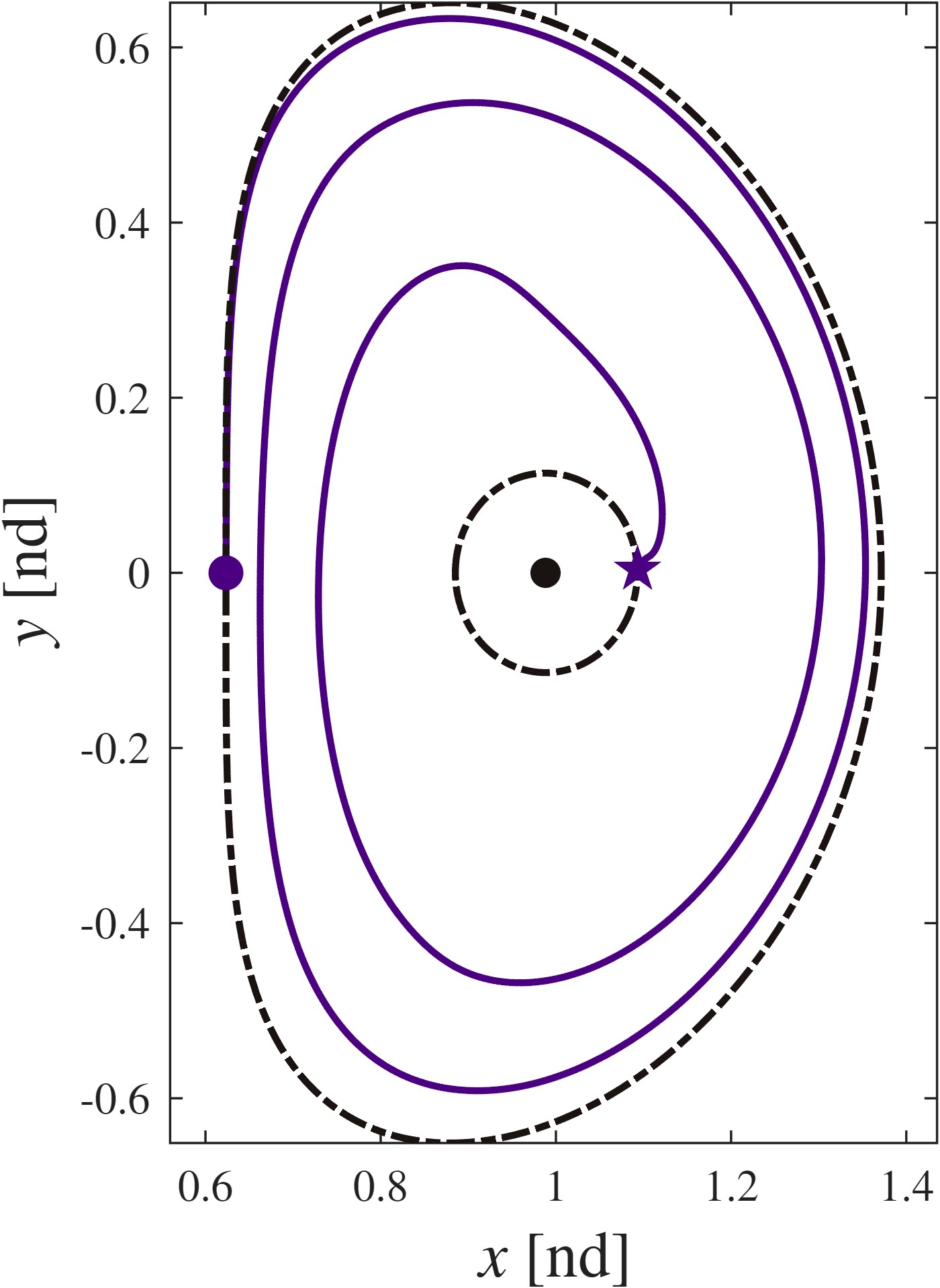}
         \caption{$u_{\max}=0.09$}
     \end{subfigure}\hfill
     \begin{subfigure}{0.24\textwidth}
         \centering
         \includegraphics[width=1\textwidth]{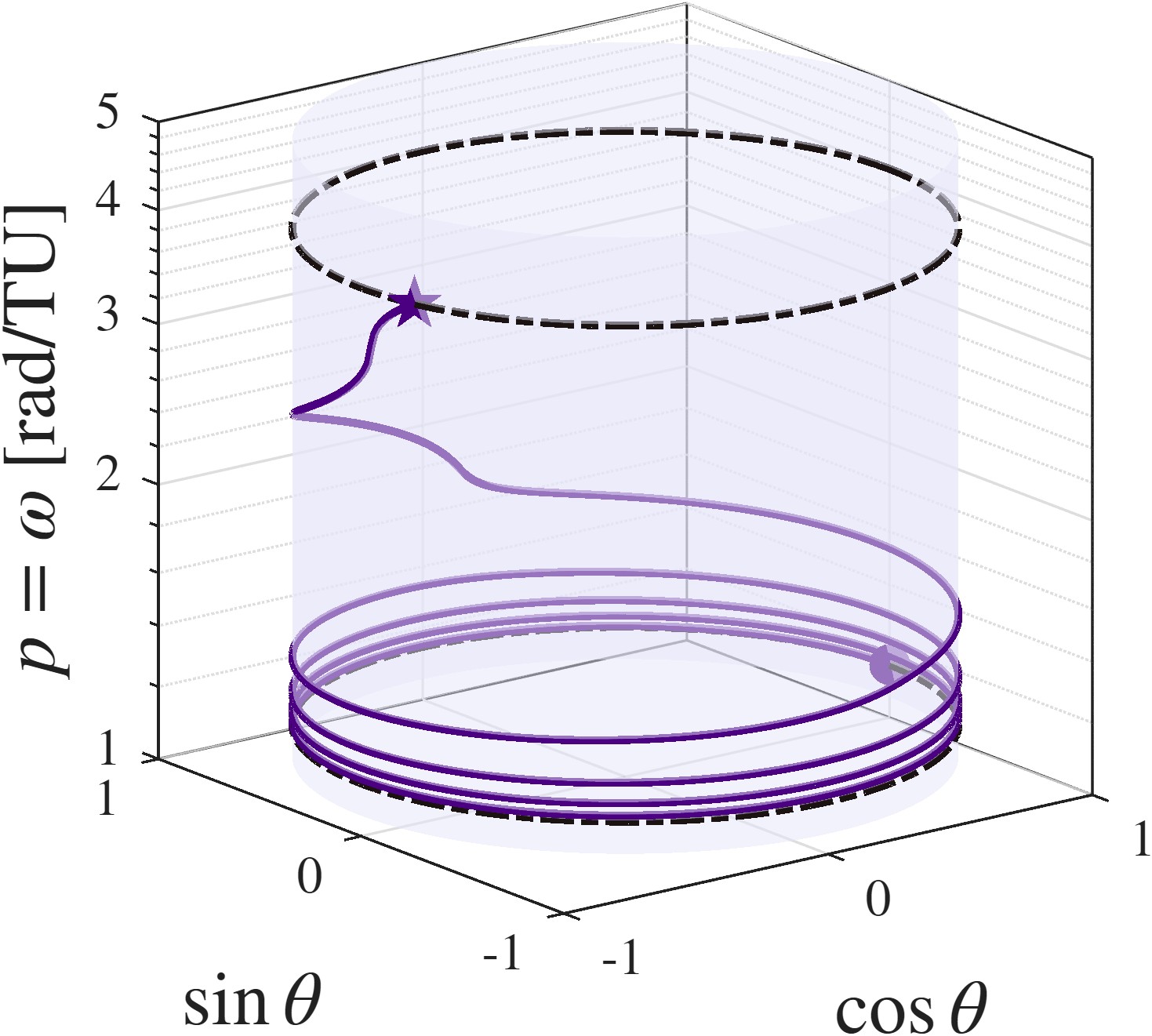}\\
         \includegraphics[width=1\textwidth]{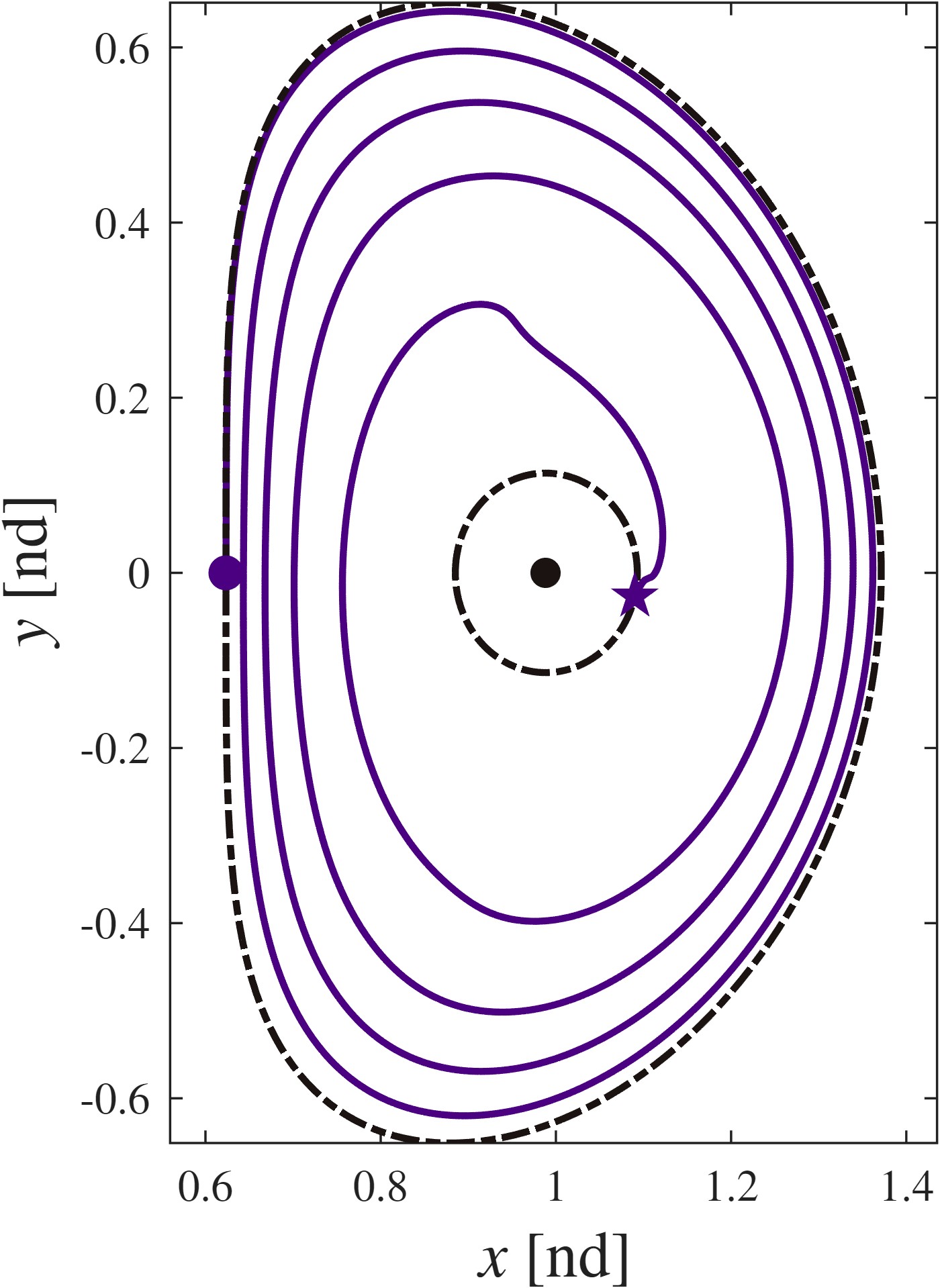}
         \caption{$u_{\max}=0.05$}
     \end{subfigure}\hfill
     \begin{subfigure}{0.24\textwidth}
         \centering
         \includegraphics[width=1\textwidth]{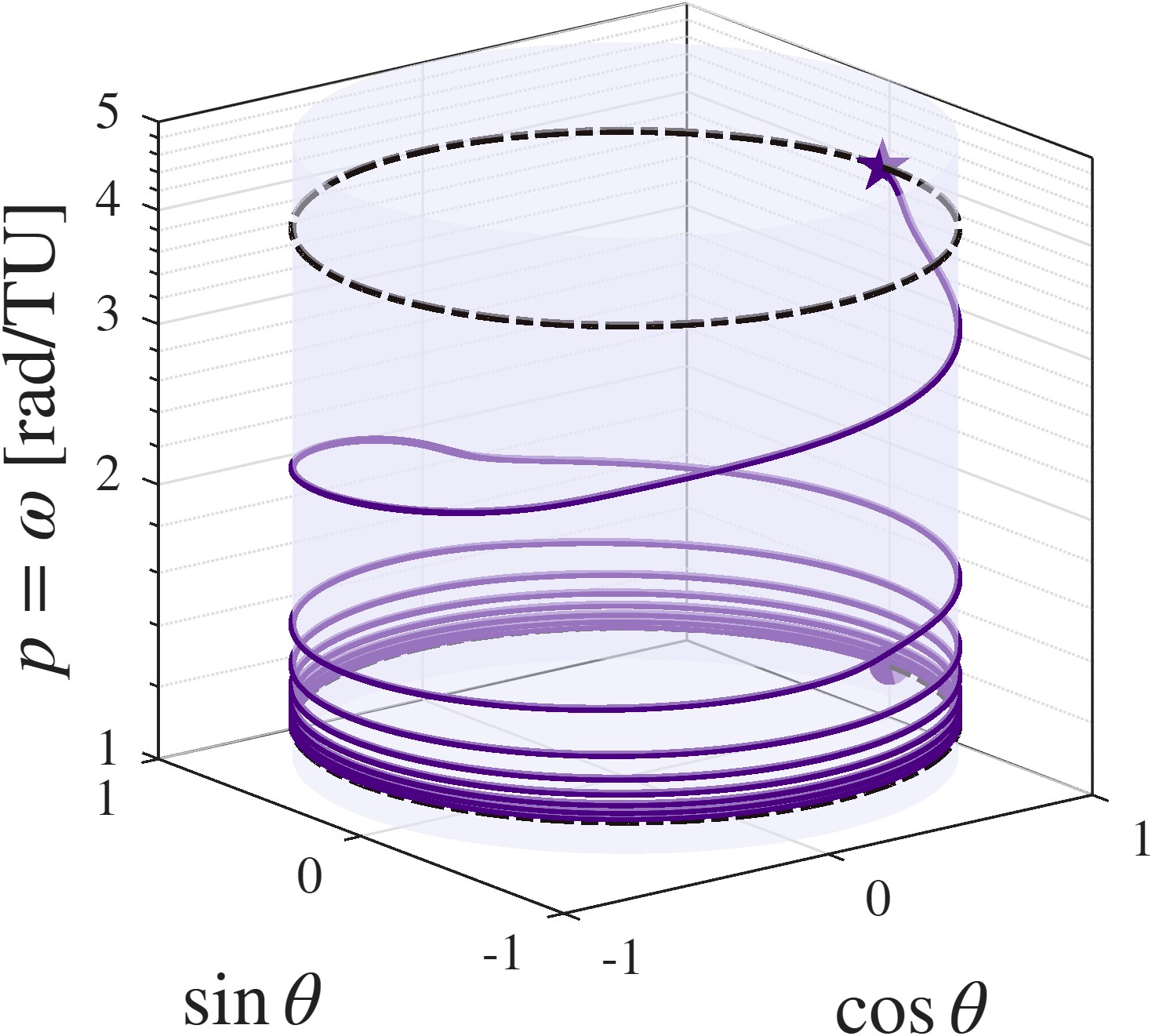}\\
         \includegraphics[width=1\textwidth]{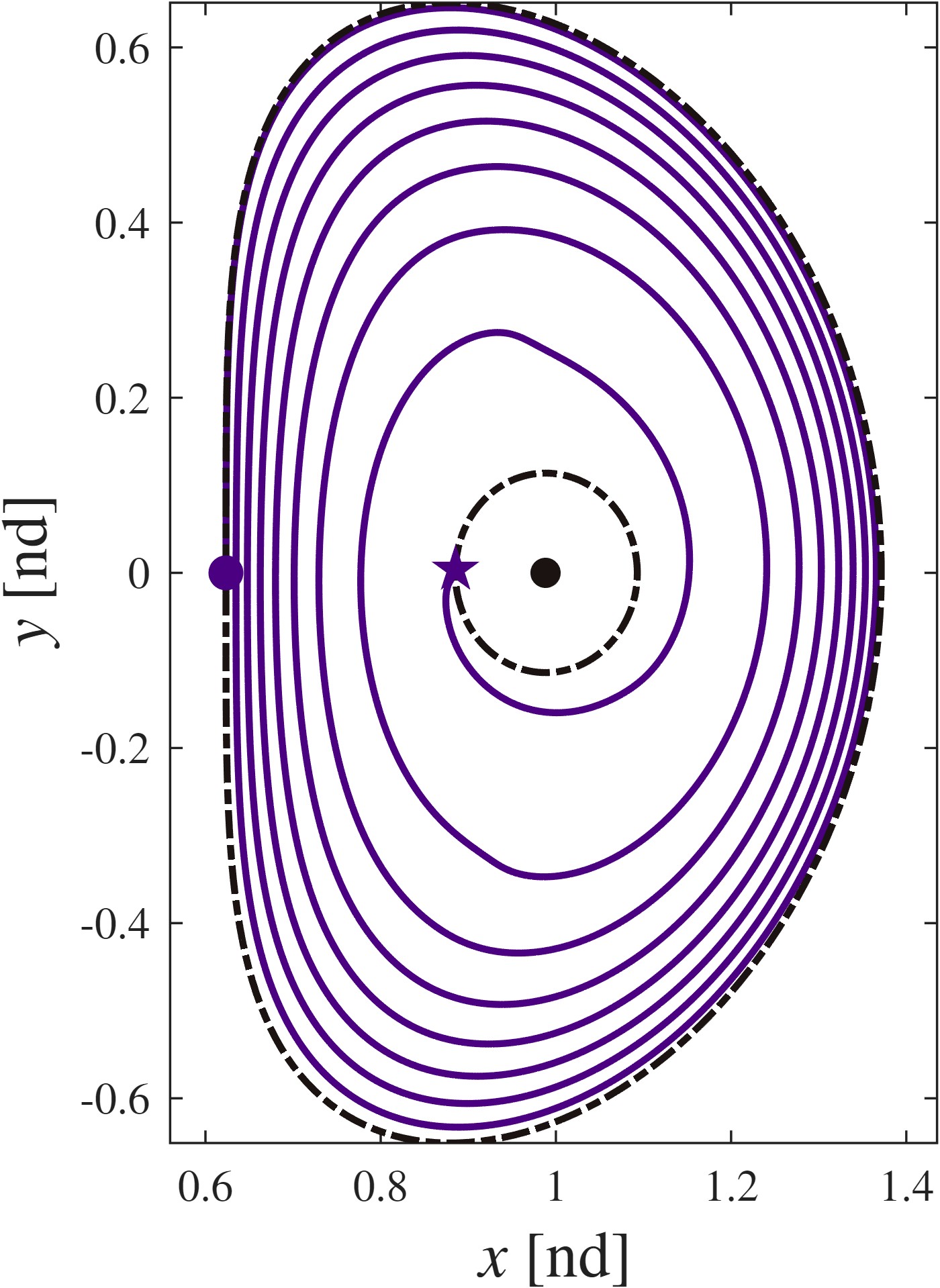}
         \caption{$u_{\max}=0.03$}
     \end{subfigure}
     \caption{Family space initialization trajectories for the 8:7 to 4:1 sidereal resonant DRO transfer.}
     \label{fig:famspacetrajs}
\end{figure}

To demonstrate the proposed initialization scheme, a 8:7 to 4:1 sidereal resonant orbit transfer is used in the DRO family with orbital frequency as the family parameter. Initialization trajectories are shown for different control magnitudes in both the family space represented by the cylinder, and the phase space in Fig. \ref{fig:famspacetrajs}.
As maximum control authority is decreased, the number of revolutions increases, with concentration decreasing as the trajectory approaches the Moon's gravity well. Abnormal phase space kinks created by Eq. \ref{eq:reducederror} can be seen in all phase space mapped trajectories.
\section{Direct Orthogonal Polynomial Integral Collocation}

Solutions are continued into the phase space via direct orthogonal polynomial integral collocation (DOPIC) on CGL nodes. An in depth explanation of the method is provided by Ahrens et al.\cite{ahrens2026direct}. A brief treatment is shown here for completeness. DOPIC takes a modal approach to the continuous optimal control problem's conversion to discrete parameter optimization. A set of orthogonal polynomial coefficients approximates only the highest-order derivative of the system dynamics and reconstructs lower-order states via successive integration matrices. This effectively halves the dimension of the state component of the transcribed nonlinear program's decision vector for second-order systems by utilizing the information-gain/smoothing property of integration while maintaining the same degree of sparsity as more common pseudospectral methods that rely on differentiation. The control is still parametrized as discrete values on the respective grid nodes.

This work uses a domain segmentation approach with multiple sets of chained polynomials guided by the number of revolutions predicted by the family space. For explanation purposes, consider a trajectory in the CR3BP approximated by a single set of polynomials. The time domain $t \in [t_0, t_f]$ is mapped to the standard domain $\tau \in [-1, 1]$ of the CGL grid. Then each component to acceleration is approximated by coefficients as
\begin{equation}
    \ddot{r}_x(\tau) \approx \sum_{j=0}^{N_{\text{poly}}} \alpha_{x,j} \varphi_j(\tau),\quad\ddot{r}_y(\tau) \approx \sum_{j=0}^{N_{\text{poly}}} \alpha_{y,j} \varphi_j(\tau),\quad\ddot{r}_z(\tau) \approx \sum_{j=0}^{N_{\text{poly}}} \alpha_{z,j} \varphi_j(\tau)
\end{equation}
where $N_{\text{poly}}$ is the polynomial degree of transcription. Including the spherical control angles and throttling of Eq. \ref{eq:cr3bpcnt} at each node, the problem's decision variable vector is written
\begin{equation}
    \boldsymbol{X} = [\boldsymbol{\alpha}_x,\, \boldsymbol{\alpha}_y,\,\boldsymbol{\alpha}_z,\,\boldsymbol{\phi},\,\boldsymbol{\psi},\,\boldsymbol{\delta},\,t_f]^\text{T}\in\mathbb{R}^{6(N_\text{poly}+1)+1}
\end{equation}

By applying standard linear integration operators, the velocity and position state profiles at any point on the internal collocation grid are structured explicitly as functions of the acceleration coefficients $\boldsymbol{\alpha}$ and the initial boundary states of the segment:
\begin{align}
    \dot{\boldsymbol{r}}(\tau)&=\frac{\Delta t}{2} \mathbf{B}_1 \boldsymbol{\alpha} + \dot{\boldsymbol{r}}_0 \\
    \boldsymbol{r}(\tau) &= \left(\frac{\Delta t}{2}\right)^2 \mathbf{B}_2 \mathbf{B}_1 \boldsymbol{\alpha} + \frac{\Delta t_s}{2}(\tau + 1)\dot{\boldsymbol{r}}_0 + \boldsymbol{r}_0
\end{align}
where $\mathbf{B}_1$ and $\mathbf{B}_2$ represent the first and second-order integral collocation matrices, respectively. After the lower order states have been computed from the coefficients, they can be plugged into Eq. \ref{eq:cr3bpcnt} with control to form dynamic defect equality constraints at every node to enforce the flow field. The minimum time and minimum fuel cost functions are expressed as
\begin{equation}
    J=t_f,\quad J = \frac{\Delta t}{2} \sum_{j=0}^{N_{\text{poly}}} w_j \delta_j
\end{equation}
where $w_j$ are the CGL quadrature weights\cite{ahrens2026direct}. For the minimum time problem, $\boldsymbol{\delta}$ components are removed from the decision vector because throttle is also set to a maximum. With the domain segmentation formulation, continuity constraints on lower order states are enforce between segments, and initial position and velocity states per segment are included as additional free variables to maintain proper equality constraint dimensionality with respect to the number of decision variables. Lastly, for phase-free problems, $\theta_f$ is included as an additional free variable with the final constraint
\begin{equation}
    \boldsymbol{x}(t_f) = \boldsymbol{\Gamma}(p_f,\theta_f)
\end{equation}

Overall, direct orthogonal polynomial integral collocation represents a powerful framework for low-thrust trajectory optimization in chaotic multi-body regimes by fundamentally reforming how the system dynamics are transcribed. By parameterizing only the highest-order acceleration profiles and reconstructing the position and velocity states analytically via successive integration matrices, the formulation avoids treating these kinematic components as independent decision variables. This second-order integration structure results in an exceptional reduction in the overall size of the nonlinear program's decision space, yielding highly sparse Jacobian matrices. Consequently, the optimization process benefits from significantly enhanced convergence stability and numerical robustness, enabling the rapid, reliable generation of complex, many-revolution trajectories over long time horizons.

\section{Numerical Results}
The family space trajectories from Fig. \ref{fig:famspacetrajs} are now transitioned to the phase space via the direct integral collocation method. First, minimum time and minimum fuel solutions with phase-free boundary conditions are presented, followed by a time-varying phase-fixed problem simulating a constellation deployment scenario.

\subsection{Time-Free, Phase-Free Minimum Time and Fuel Solutions}
\begin{figure}[htbp!]
     \centering
     \begin{subfigure}{0.24\textwidth}
         \centering
         \includegraphics[width=1\textwidth]{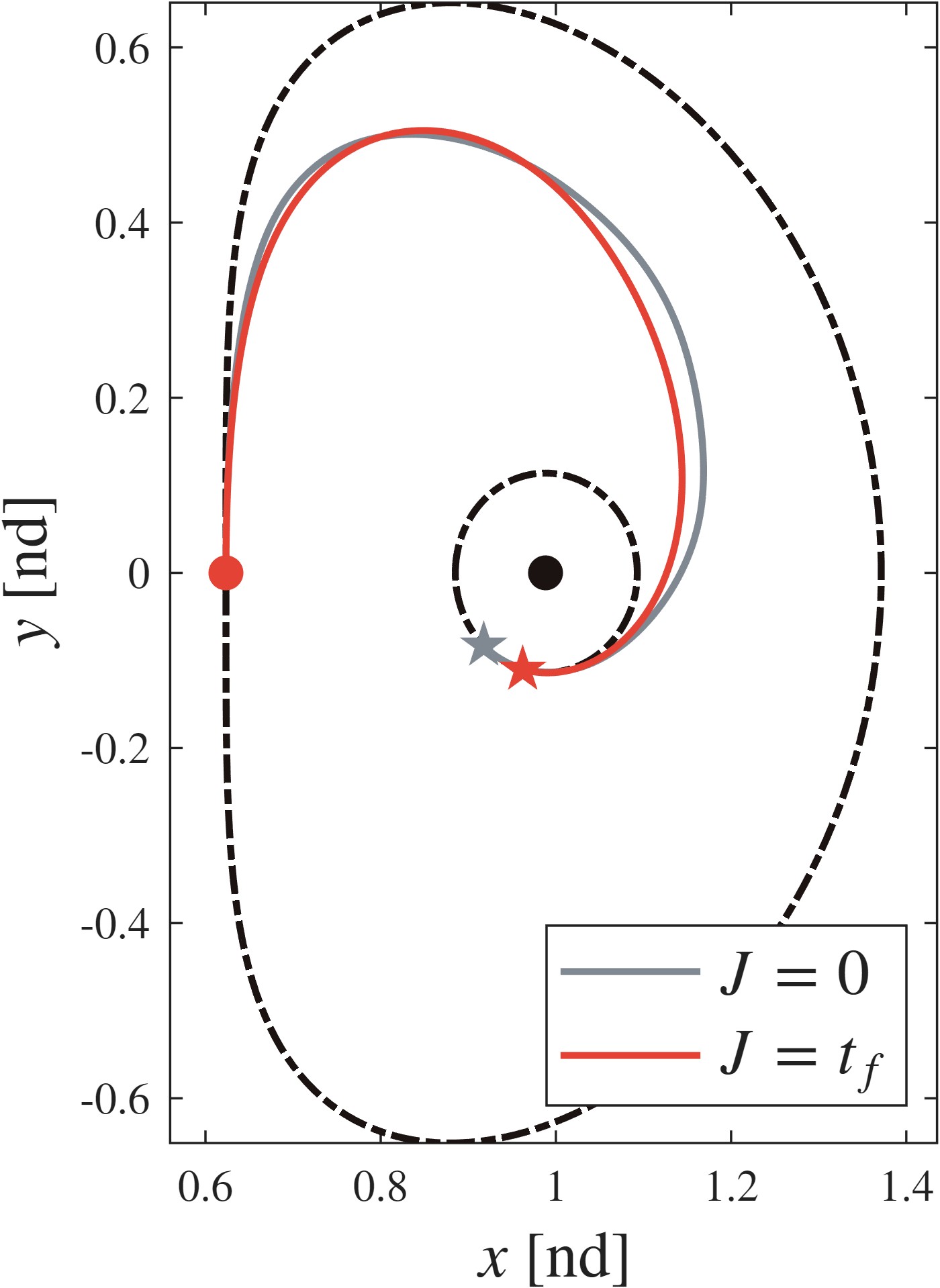}\\
         \includegraphics[width=1\textwidth]{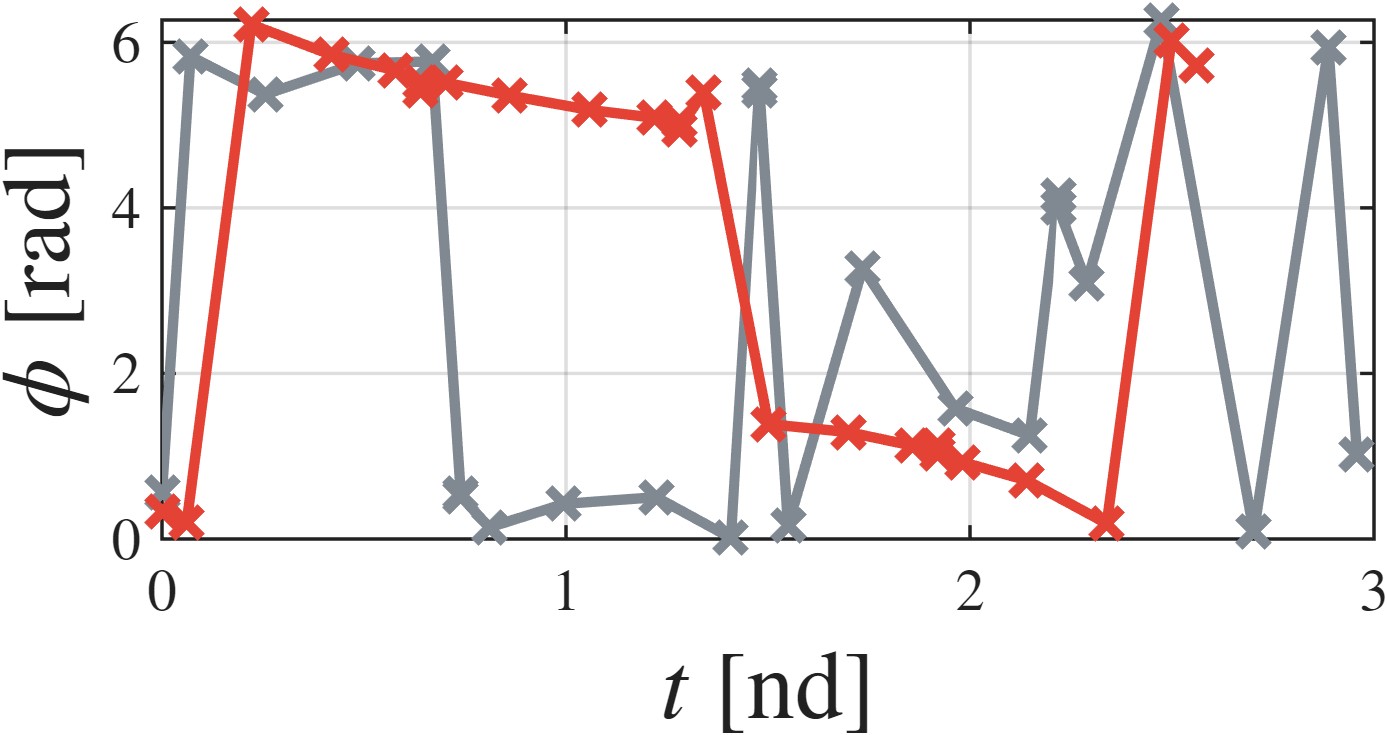}
         \caption{$u_{\max}=0.3$}
     \end{subfigure}\hfill
     \begin{subfigure}{0.24\textwidth}
         \centering
         \includegraphics[width=1\textwidth]{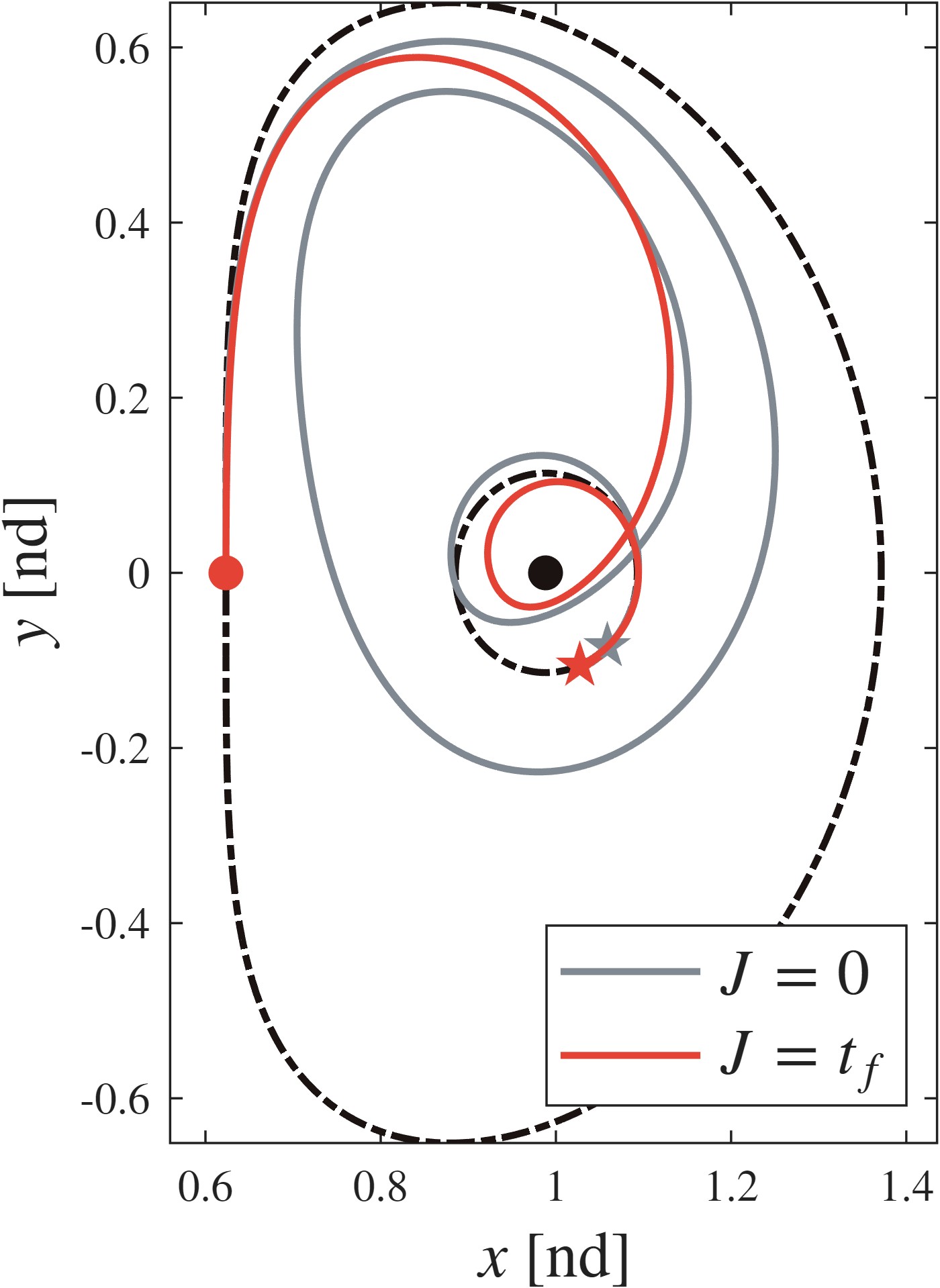}\\
         \includegraphics[width=1\textwidth]{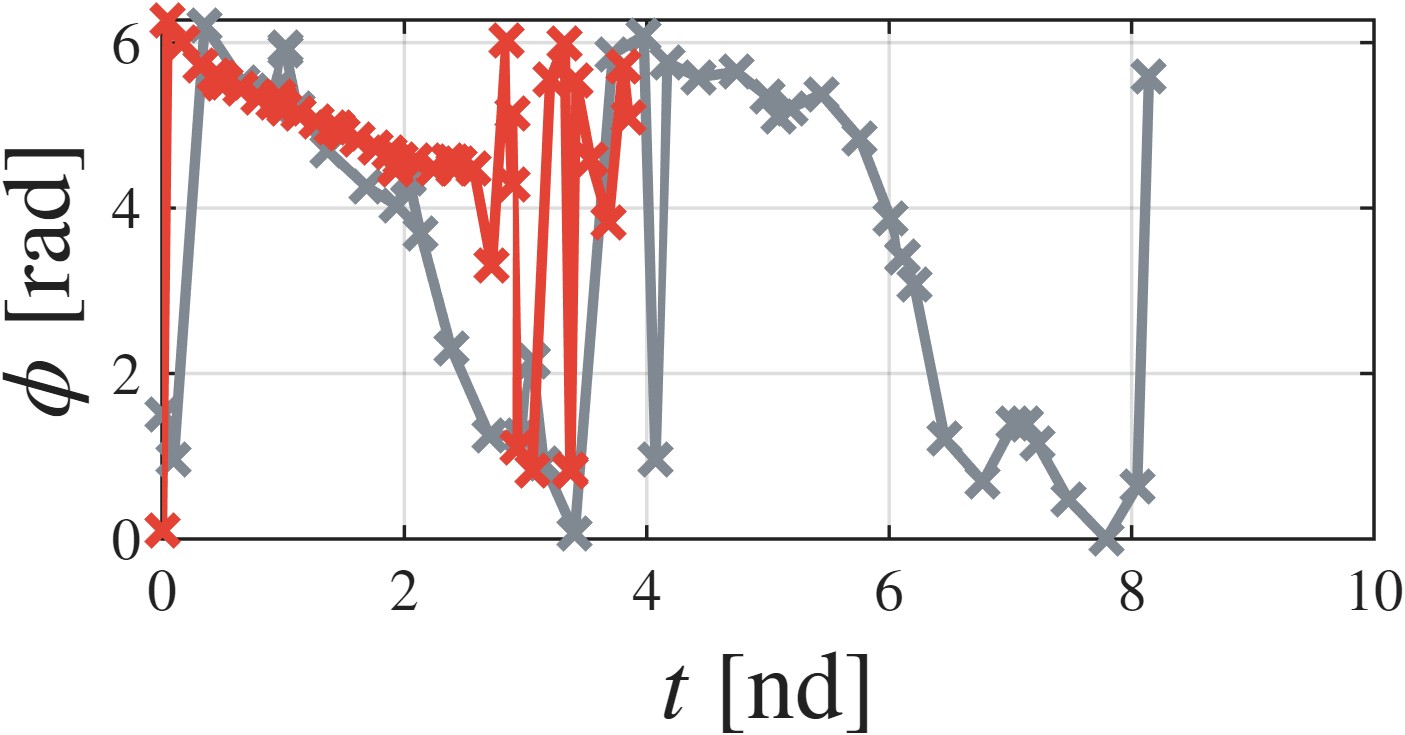}
         \caption{$u_{\max}=0.09$}
     \end{subfigure}\hfill
     \begin{subfigure}{0.24\textwidth}
         \centering
         \includegraphics[width=1\textwidth]{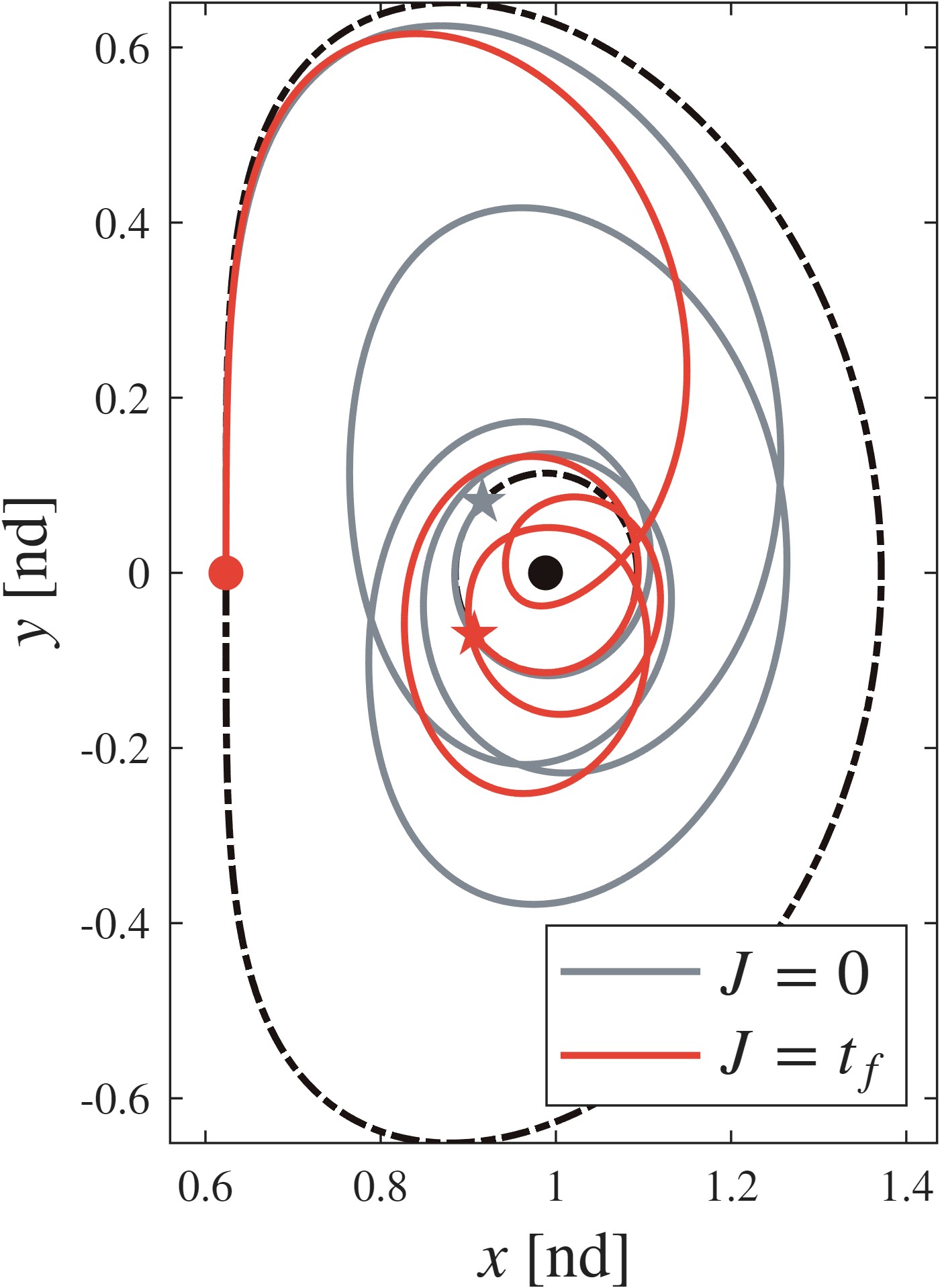}\\
         \includegraphics[width=1\textwidth]{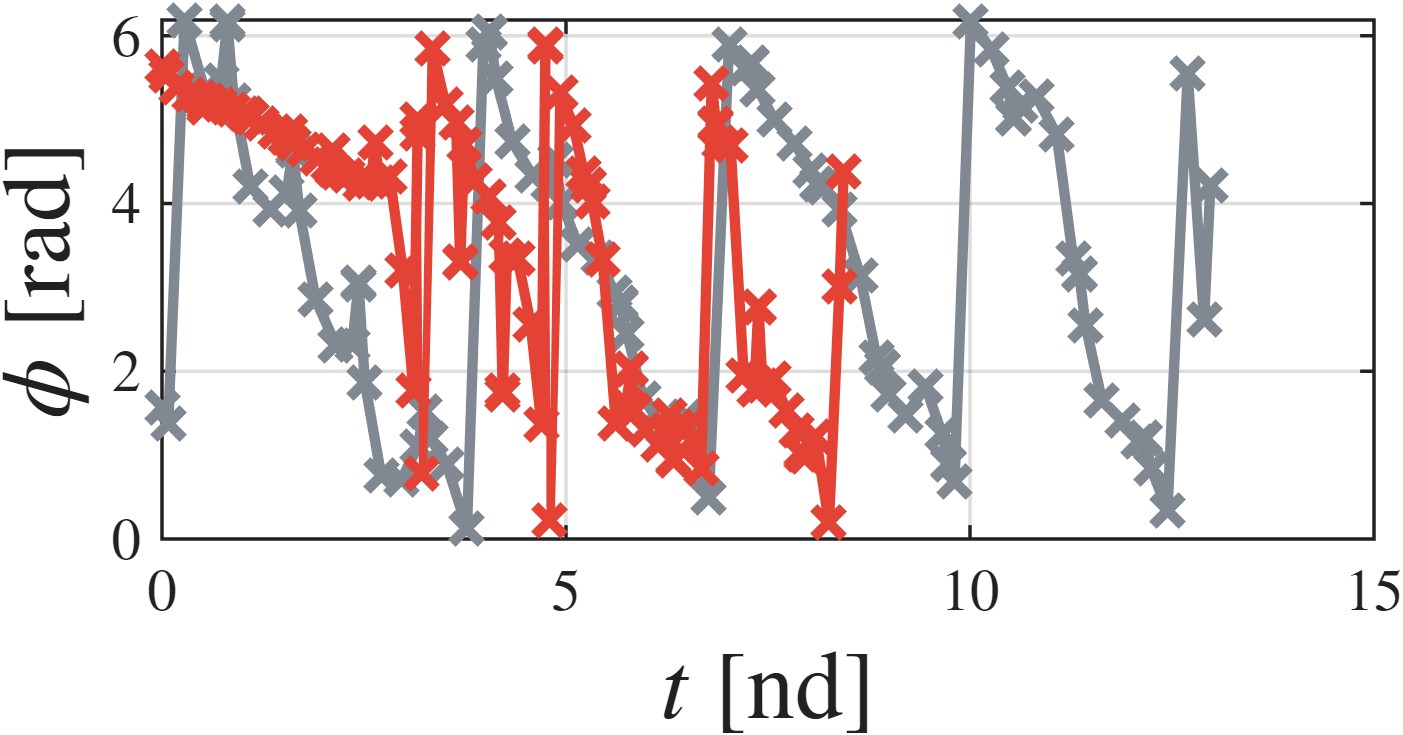}
         \caption{$u_{\max}=0.05$}
     \end{subfigure}\hfill
     \begin{subfigure}{0.24\textwidth}
         \centering
         \includegraphics[width=1\textwidth]{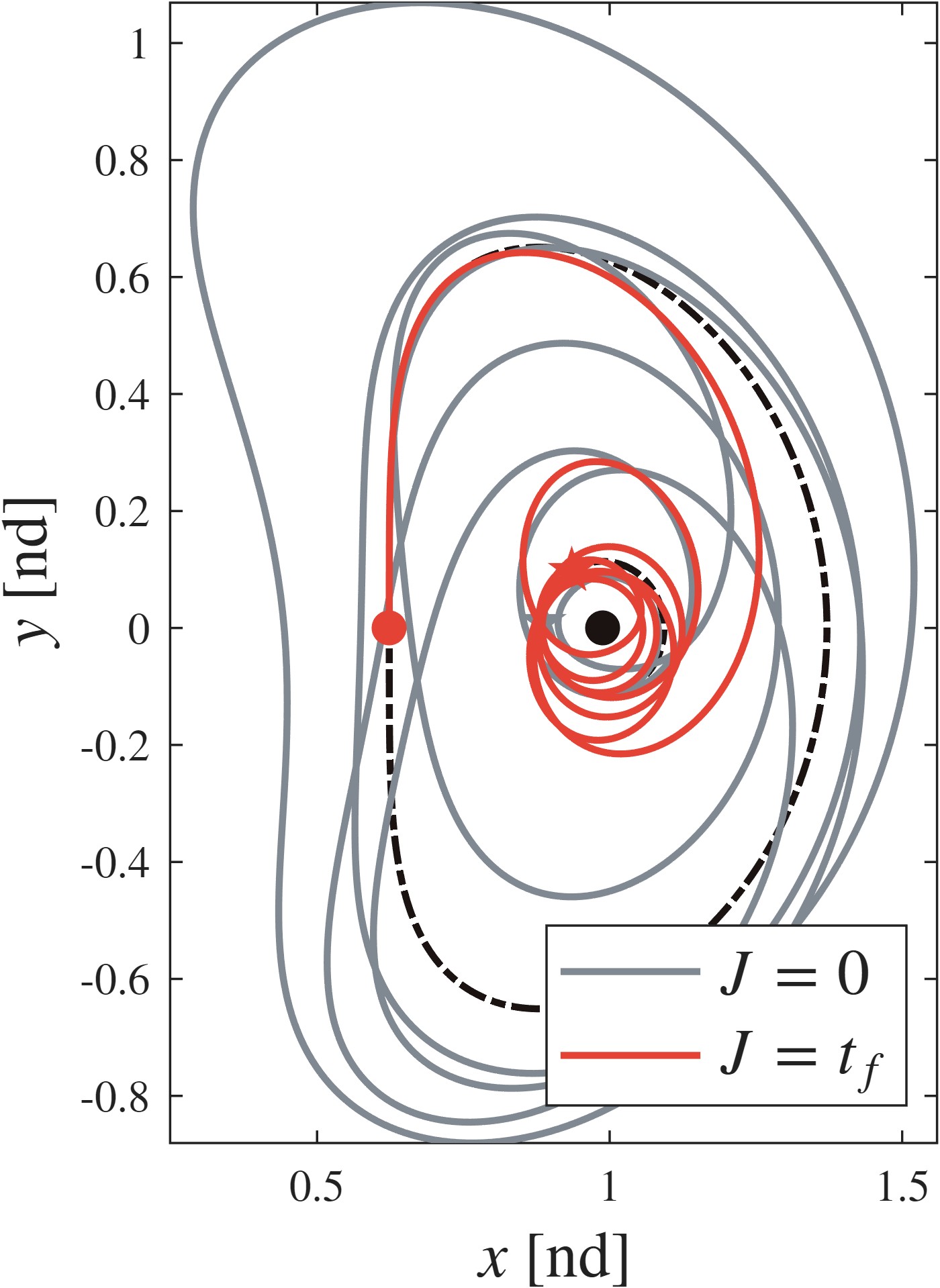}\\
         \includegraphics[width=1\textwidth]{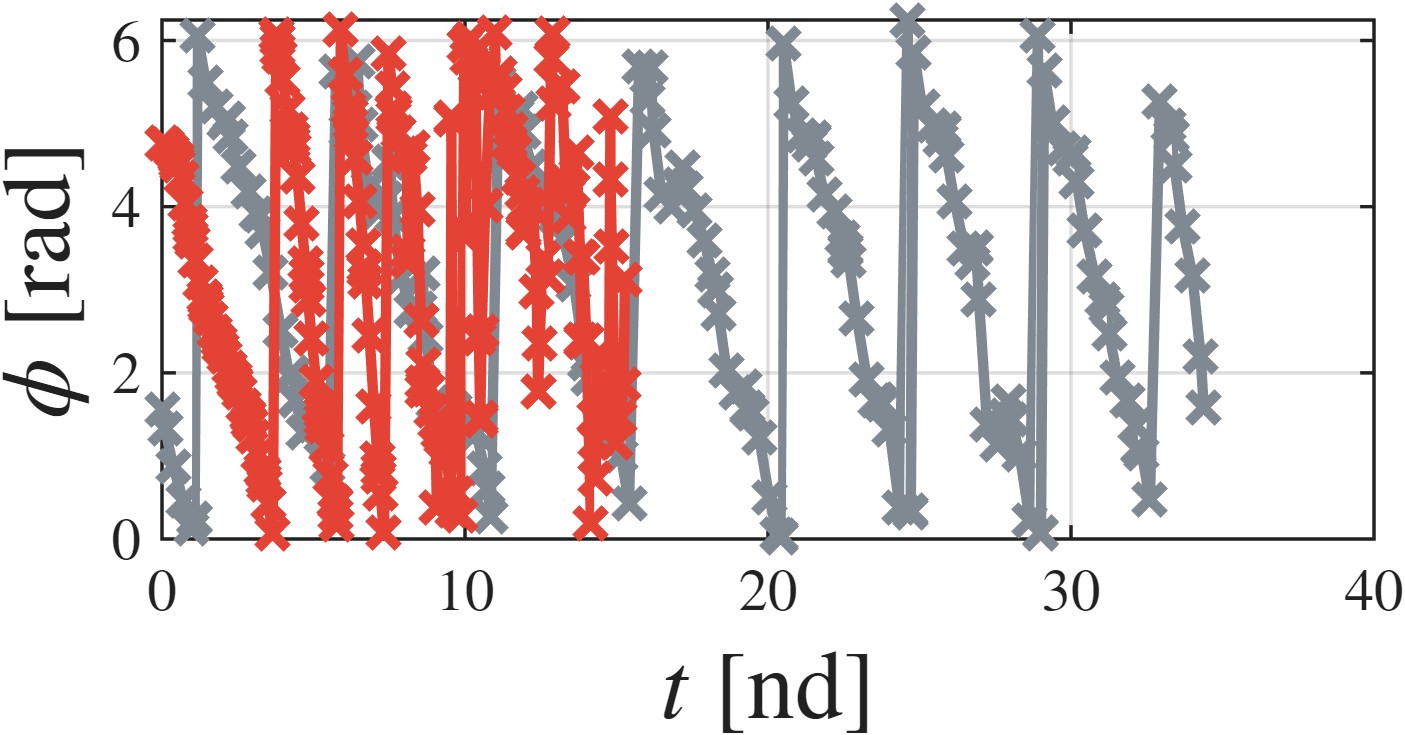}
         \caption{$u_{\max}=0.03$}
     \end{subfigure}
     \caption{Phase-free minimum time solutions.}
     \label{fig:mintime}
\end{figure}

\begin{figure}[htbp!]
     \centering
     \begin{subfigure}{0.24\textwidth}
         \centering
         \includegraphics[width=1\textwidth]{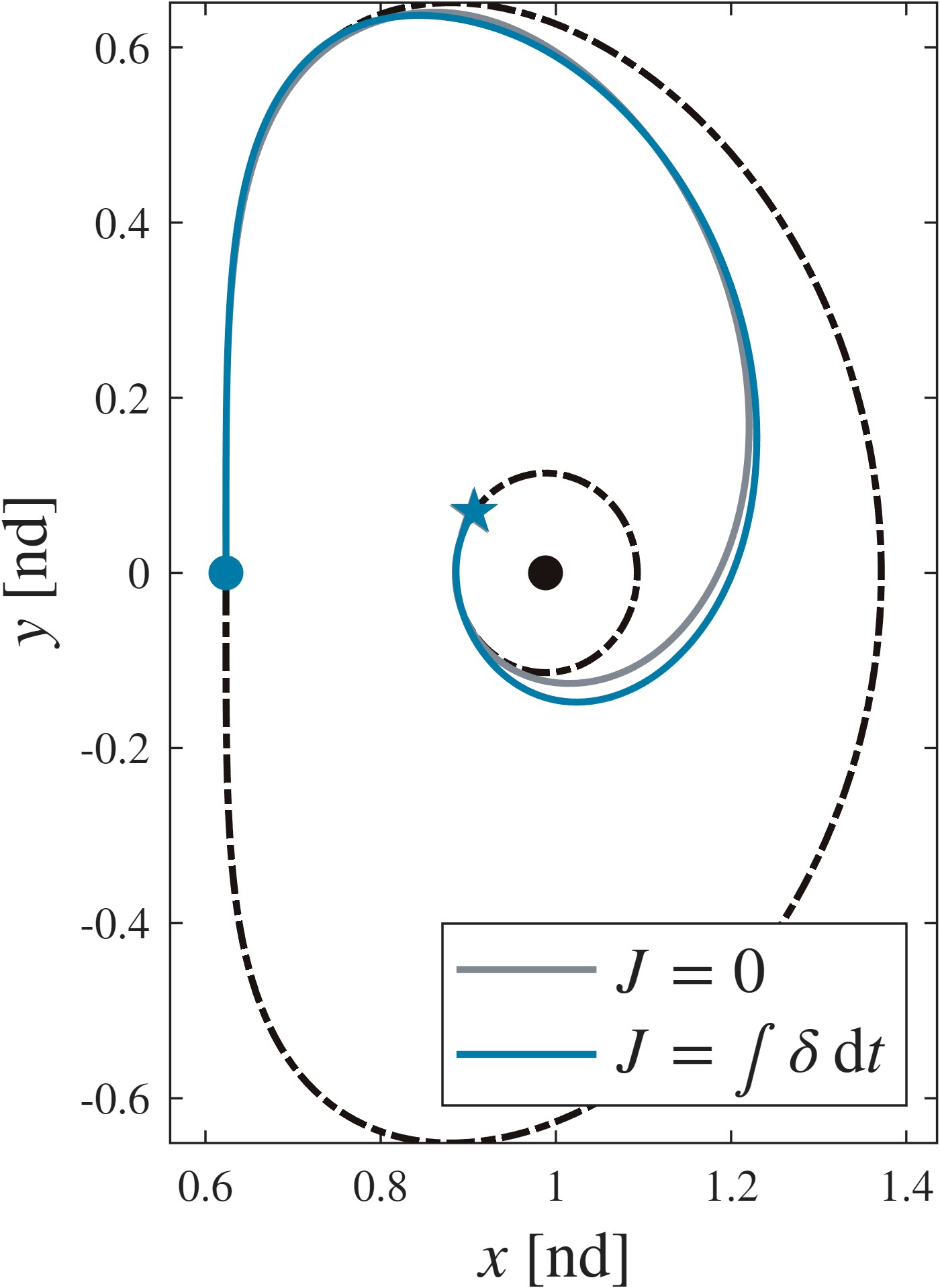}\\
         \includegraphics[width=1\textwidth]{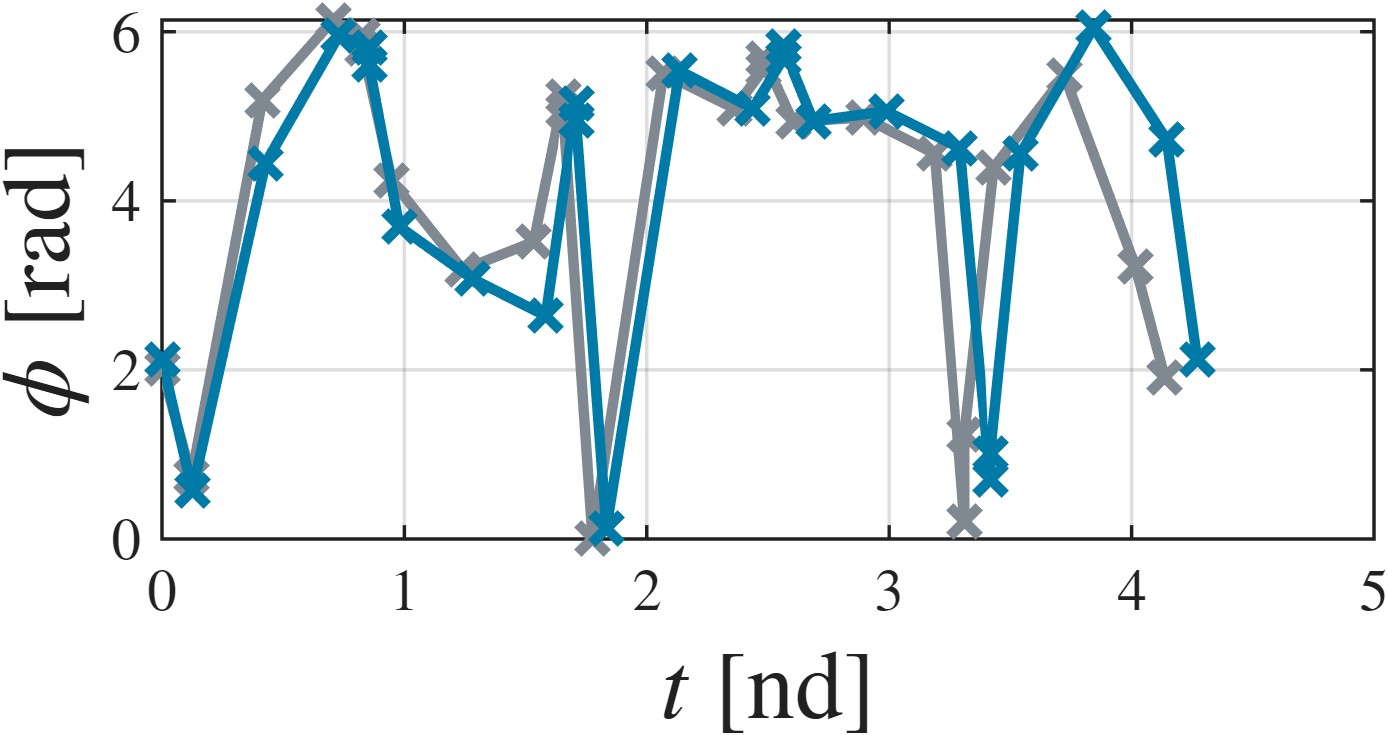}\\
         \includegraphics[width=1\textwidth]{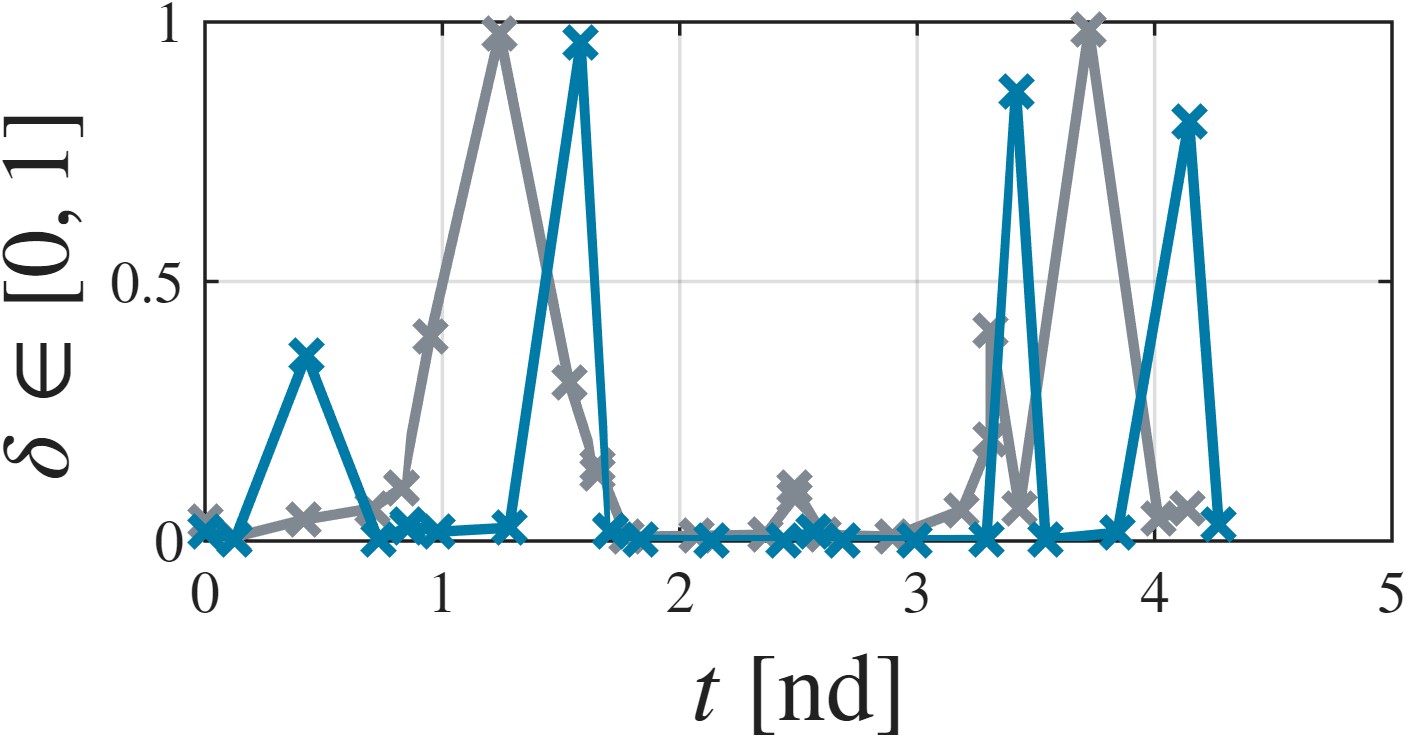}
         \caption{$u_{\max}=0.3$}
     \end{subfigure}\hfill
     \begin{subfigure}{0.24\textwidth}
         \centering
         \includegraphics[width=1\textwidth]{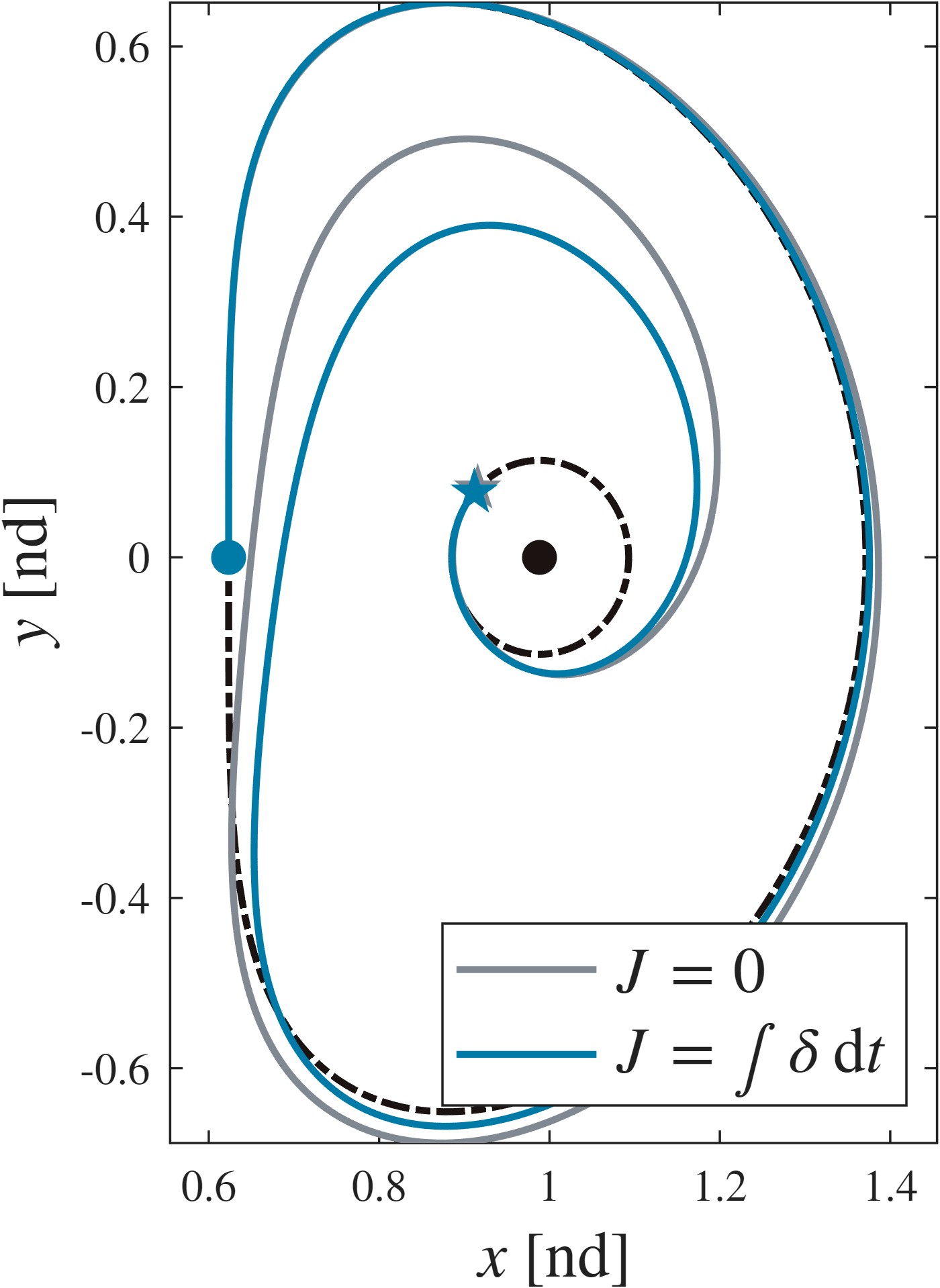}\\
         \includegraphics[width=1\textwidth]{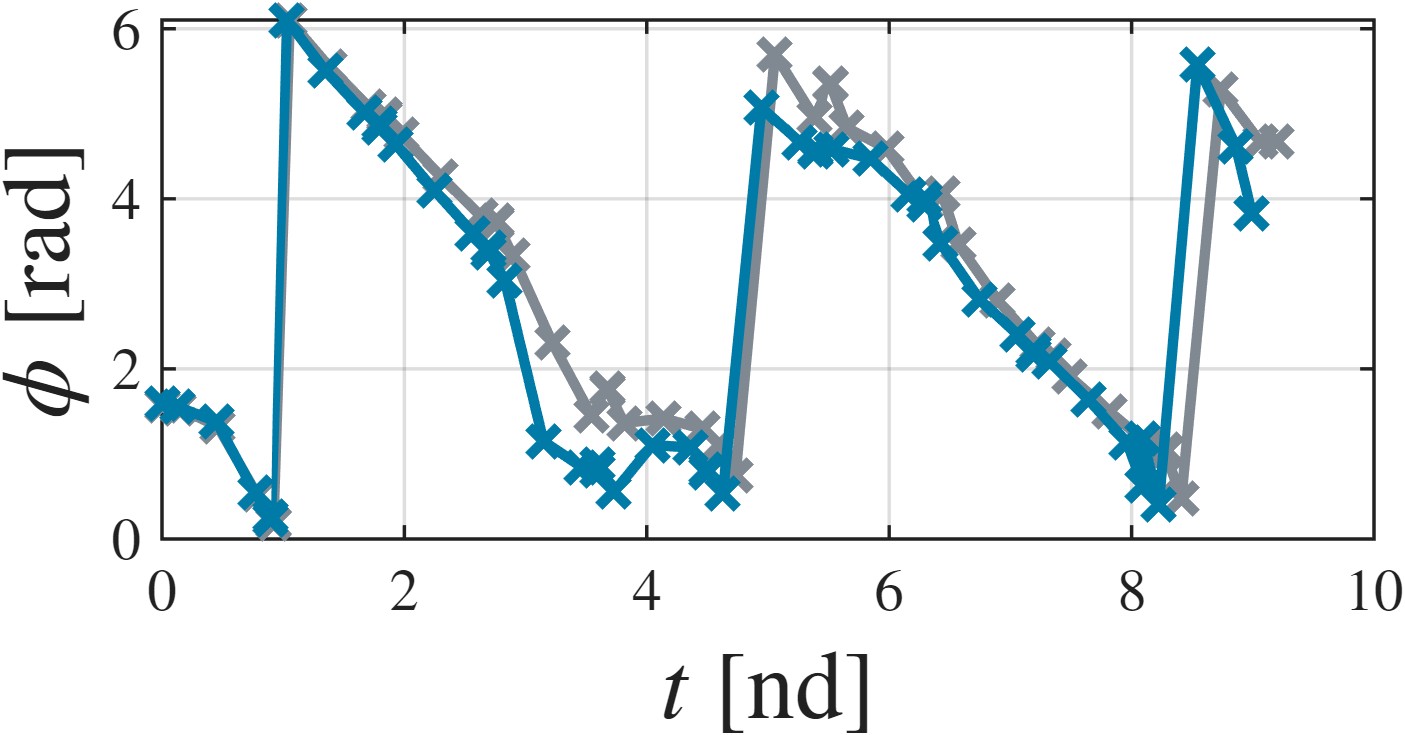}\\
         \includegraphics[width=1\textwidth]{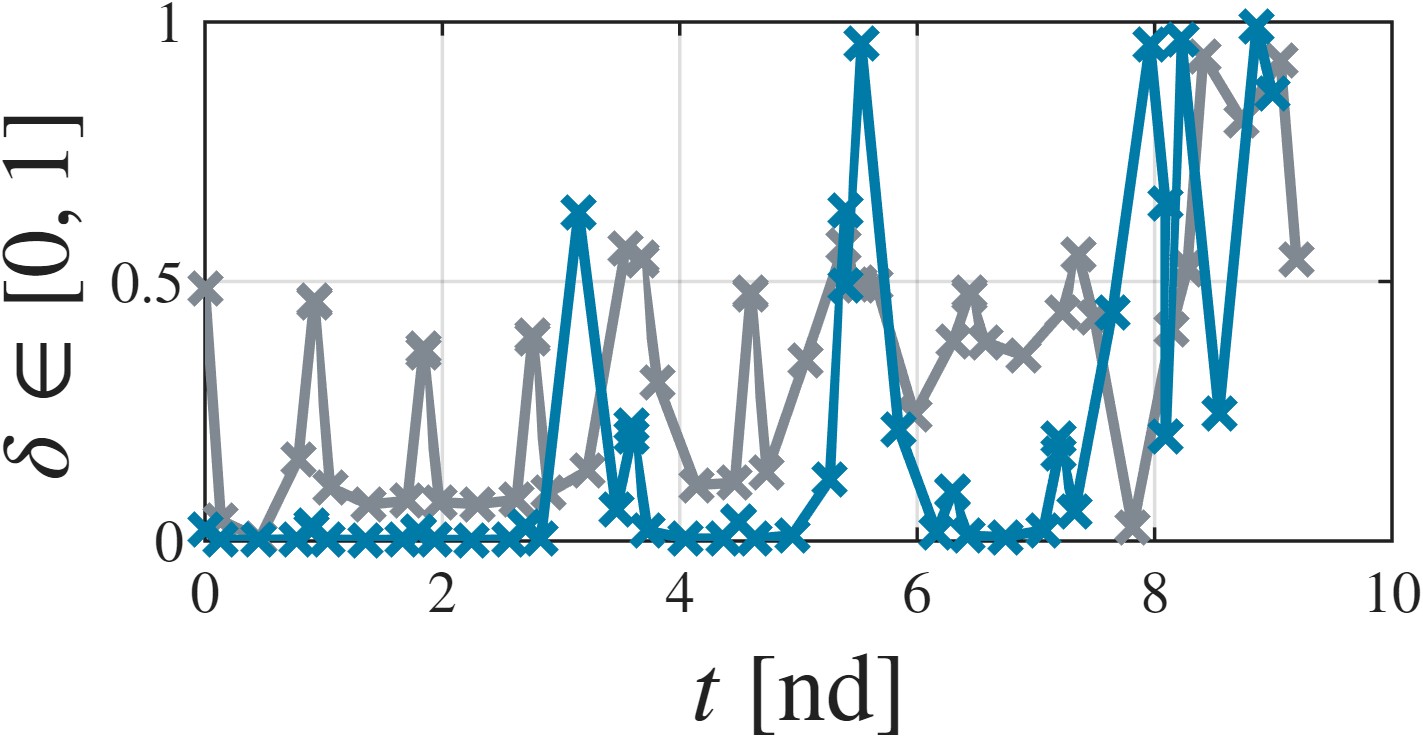}
         \caption{$u_{\max}=0.09$}
     \end{subfigure}\hfill
     \begin{subfigure}{0.24\textwidth}
         \centering
         \includegraphics[width=1\textwidth]{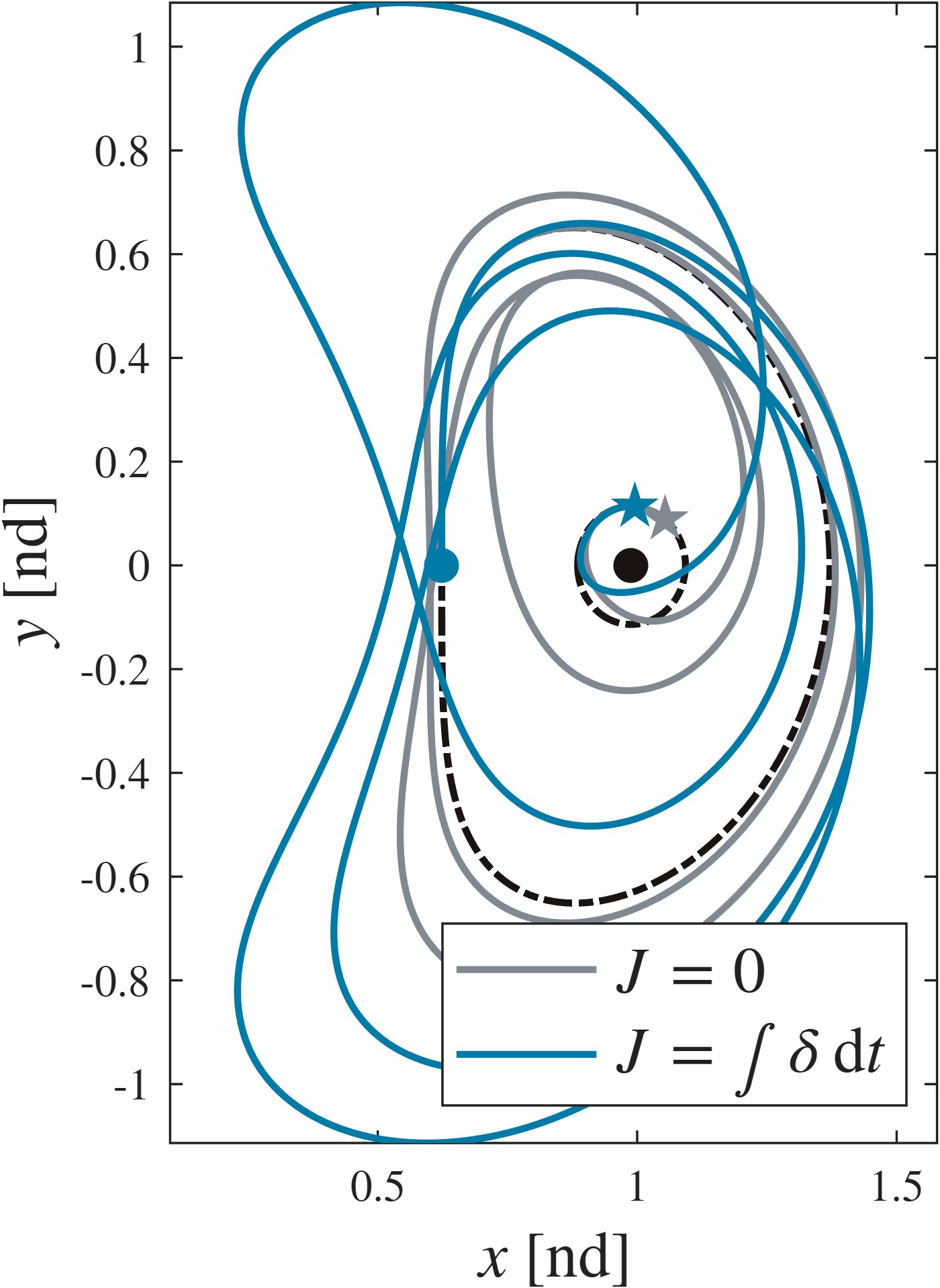}\\
         \includegraphics[width=1\textwidth]{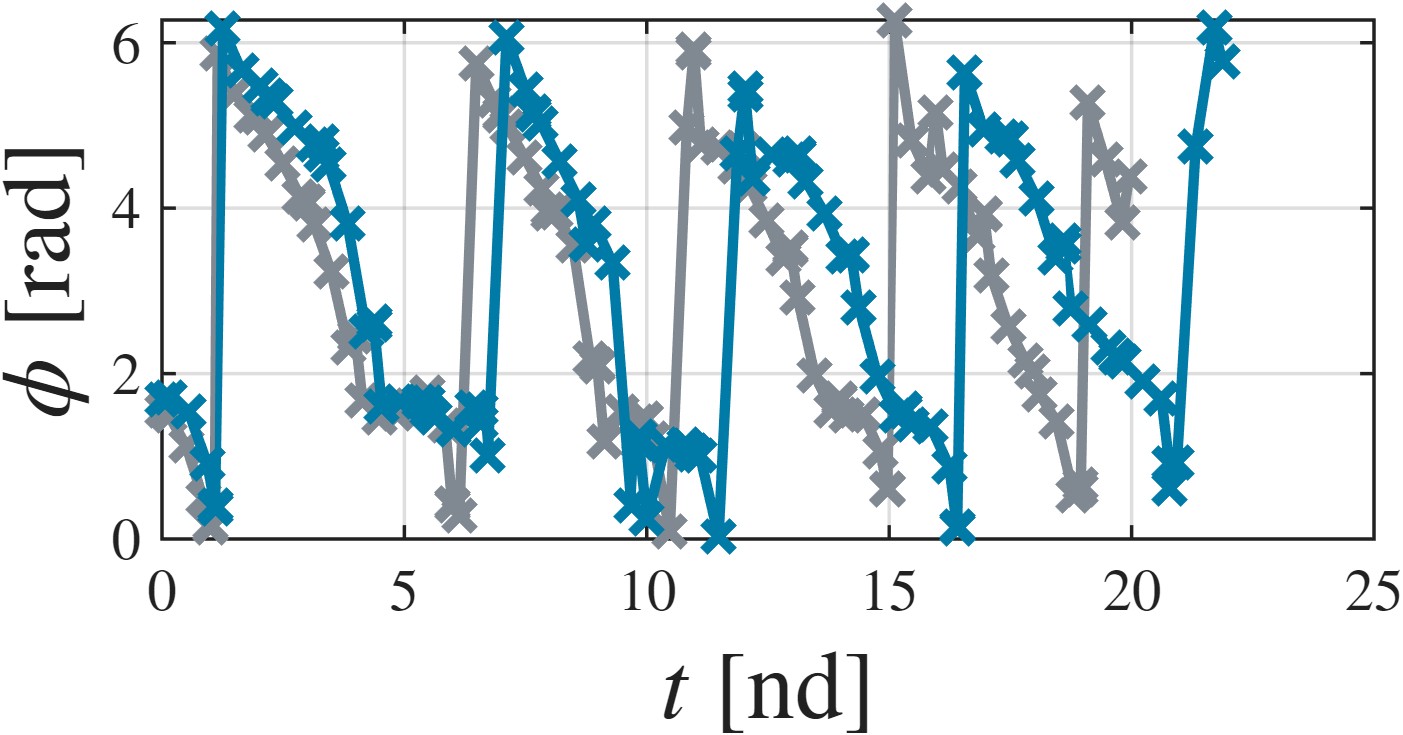}\\
         \includegraphics[width=1\textwidth]{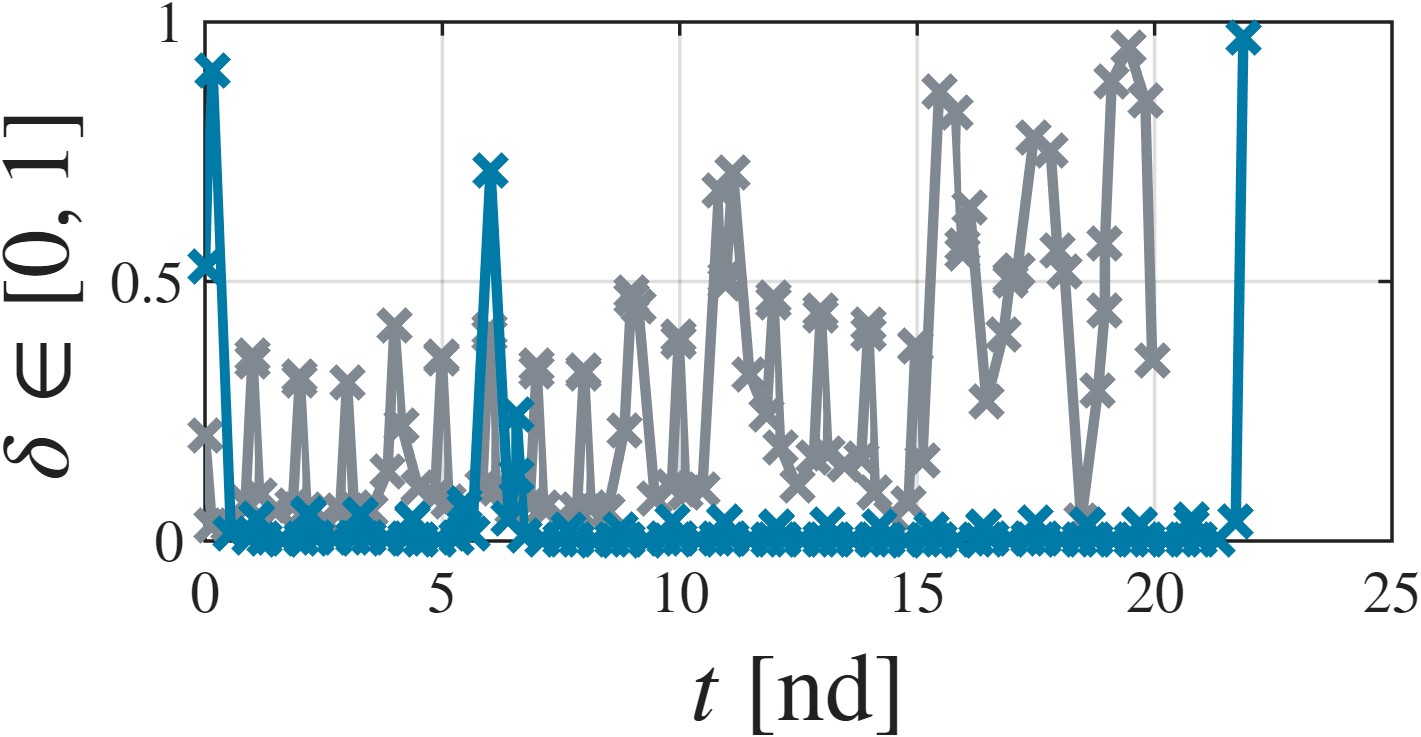}
         \caption{$u_{\max}=0.05$}
     \end{subfigure}\hfill
     \begin{subfigure}{0.24\textwidth}
         \centering
         \includegraphics[width=1\textwidth]{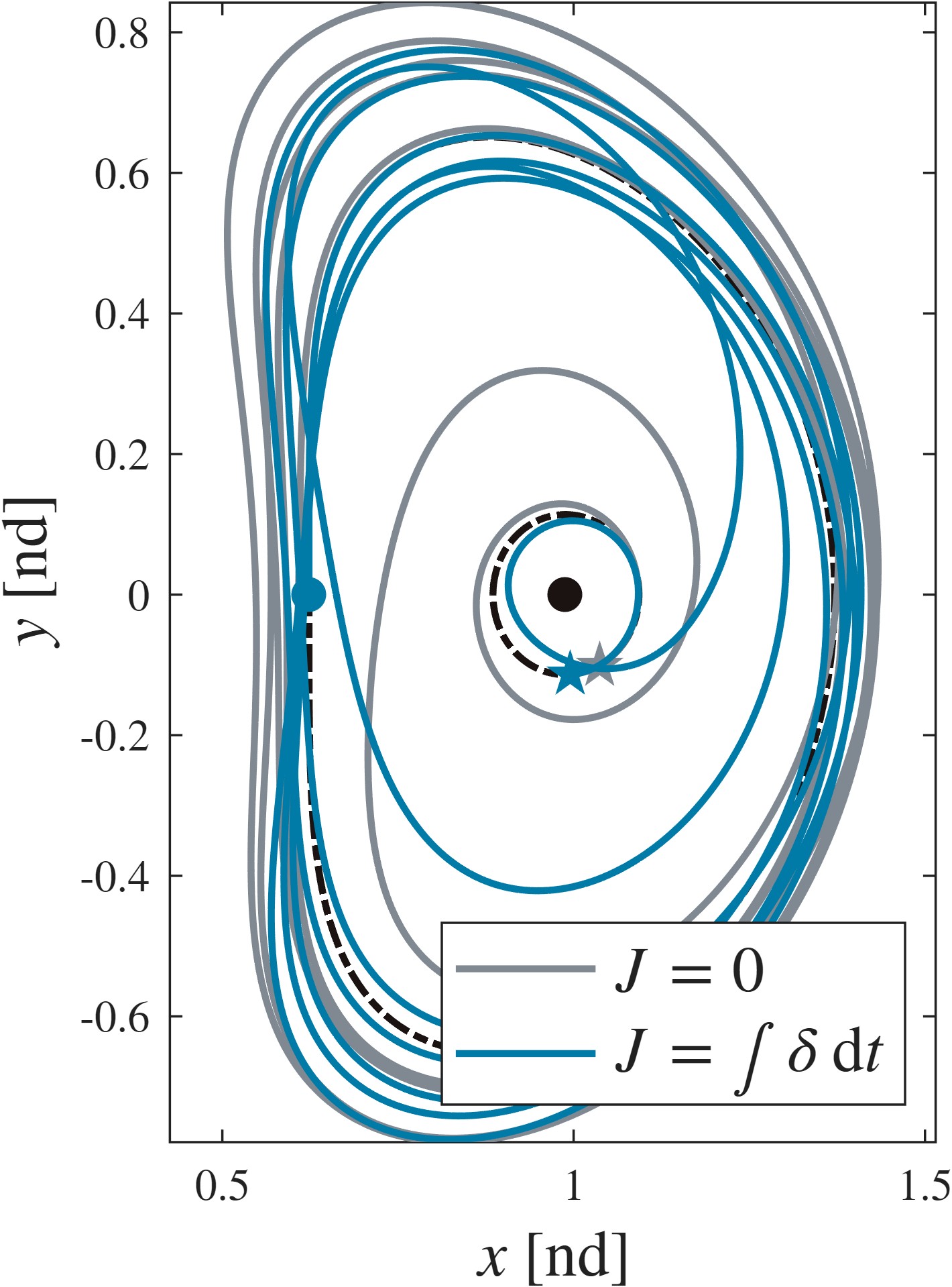}\\
         \includegraphics[width=1\textwidth]{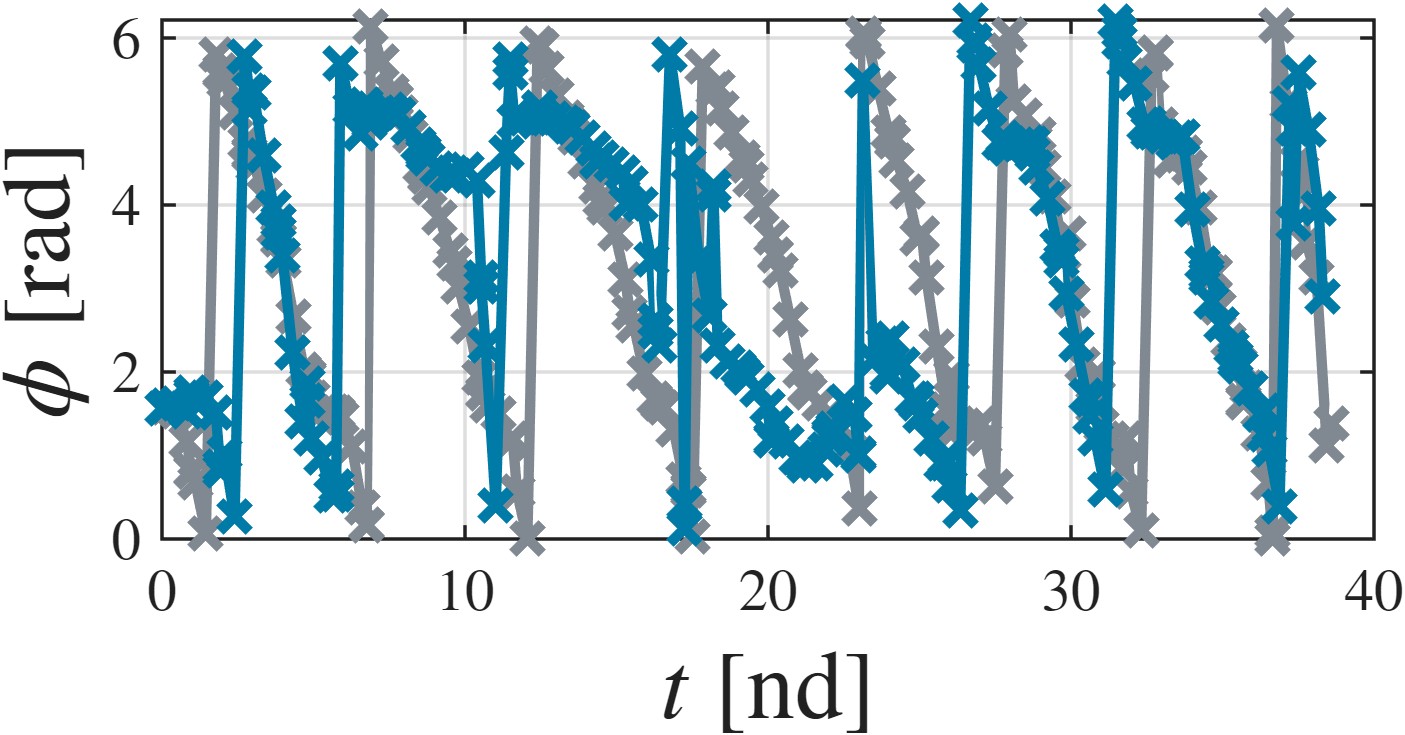}\\
         \includegraphics[width=1\textwidth]{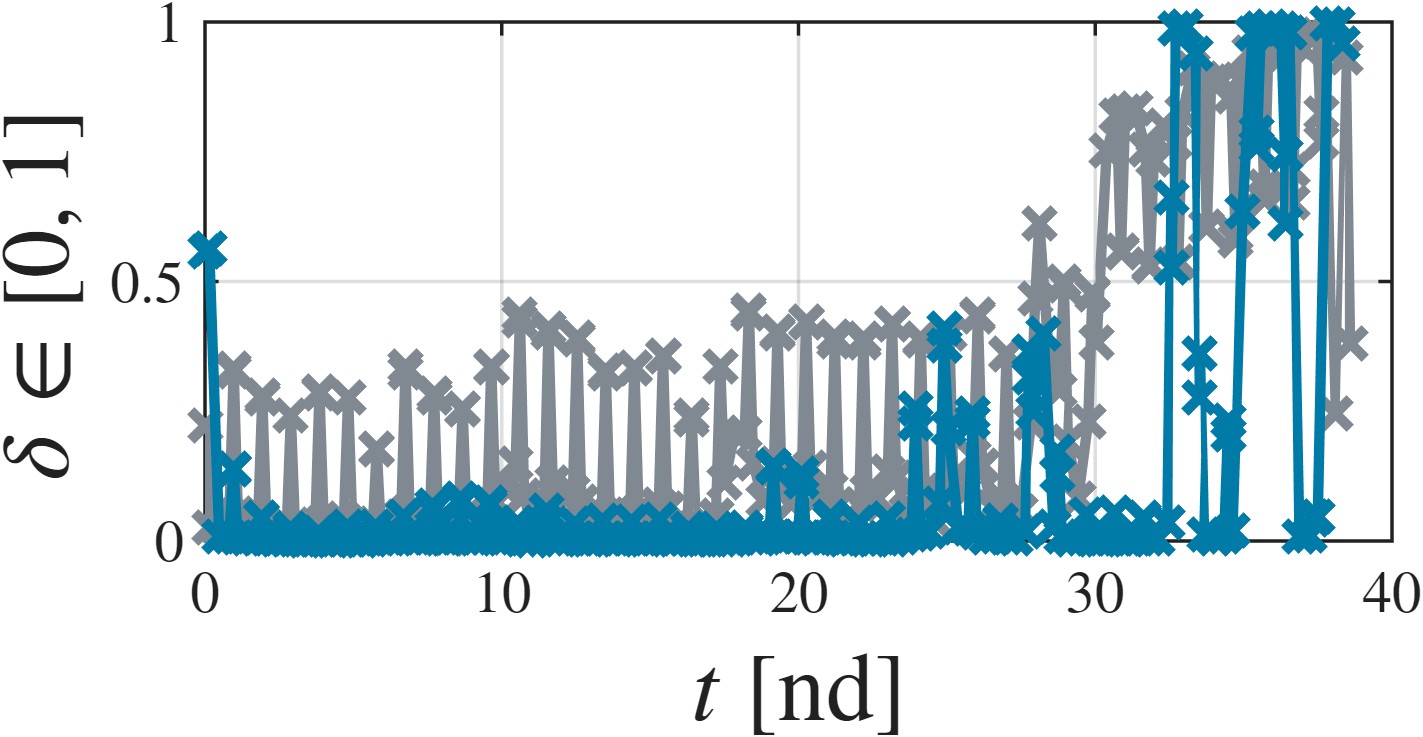}
         \caption{$u_{\max}=0.03$}
     \end{subfigure}
     \caption{Phase-free minimum fuel solutions.}
     \label{fig:minfuel}
\end{figure}

For the array of maximum control magnitudes, minimum time solutions are shown in Fig. \ref{fig:mintime}, minimum fuel solutions are shown in Fig. \ref{fig:minfuel}, and numerical results are displayed in Table \ref{tab:mintimefuel} with dimensional quantities. A feasibility problem with zero cost is first solved using the family space initialization, shown as gray trajectories, followed by the respective transcribed optimal control problem.

For both problems, the feasibility solution stays near the initialized total time of flight. Upon examining the minimum time solutions, the family space value for minimum time is always conservative -- a positive quality for direct methods because walking down towards a minimum time with a feasible solution is much easier than walking up towards the minimum time with an infeasible solution. In all cases with $u_{\max}<0.1$, the phase space converged minimum time was nearly half the predicted value. The control angle time histories are particularly smooth for $u_{\max}=0.3$ and $u_{\max}=0.03$, while all minimum time solutions show acceptable smoothness for discrete control points. This is a positive indication of well-behaved local minima. All minimum time solutions for $u_{\max}<0.1$ leverage one or more close passes of the Moon in the converged solution. When looking at the overall difference between family space solution, feasibility solution, and optimal solution, the integral collocation method is robust enough to dramatically change geometry to arrive at local optima; all trajectories shown in this work converged onto the locally optimal solutions without special treatment. Despite the geometric differences apparent in all three solution phases, the minimum time solutions still shared the same number of revolutions as the family space trajectories even with dramatically different final times. While this could be a product of the initialization trajectory staying with it's initialized revolution count basin, and no claims of global optimality are made in this work, the control angle behavior and overall dramatic time reduction are strong signs that the family space trajectory leads to a useful minimum time solution in the phase space while accurately predicting the number of revolutions. 

Without time in the cost function, the minimum fuel solutions vary very little from the family space initialized quantity, which makes the number of revolutions match largely meaningless. Still, the transition method is able to locate locally optimal, useful transfers indicated by the bang-off-bang throttle history. More segments with less nodes are used for the minimum fuel transcription to provide the solver more opportunities to resolve this control structure. The control angle history matches almost identically with the feasible solution despite  very different throttle histories because the control angle ultimately doesn't affect feasibility when throttle is turned off. With more time given (not minimized), the minimum fuel solution is allowed to find its way onto dynamical structures outside the family space. This is exemplified by the $u_{\max}=0.05$ minimum fuel trajectory, whose cost is significantly lower than any other transfer analyzed. This highlights the abundance of local minima for low-thrust trajectory design in cislunar space.

\begin{table}[htbp]
    \centering
    \caption{Phase-free minimum time and minimum fuel solution data.}
    \begin{tabular}{|l|c|c c c c|}
        \hline
        Solution & $u_{\max}$ &
        \makecell{0.300 [nd] \\ 0.817 [mm/s$^2$]} &
        \makecell{0.090 [nd] \\ 0.245 [mm/s$^2$]} &
        \makecell{0.050 [nd] \\ 0.136 [mm/s$^2$]} &
        \makecell{0.030 [nd] \\ 0.081 [mm/s$^2$]} \\
        \hline\hline
    
        \multirow{2}{*}{Family Space} & $N_\text{rev}$ & 1 & 2 & 4& 8\\ 
                                      & $t_f$ [days]   & 17.21& 52.43&95.51 & 164.6\\
        \hline
        \multirow{5}{*}{Minimum Time} & $N_\text{poly}$ & 5 & 5& 5& 5\\
                                      & $N_\text{seg}$ & 4 &8 &16 & 32\\
                                      & $N_\text{rev}$ &1 &2 & 4& 7\\
                                      & $t_f$ [days]   & 11.13& 16.74& 36.68& 66.43\\
                                      & $\Delta V$ [m/s] & 785.7& 354.6& 431.6 & 468.9\\
        \hline
        \multirow{5}{*}{Minimum Fuel} & $N_\text{poly}$ &4 &4 &4 &4 \\
                                      & $N_\text{seg}$ & 5& 10 & 20& 40\\
                                      & $N_\text{rev}$ &1 & 2& 4& 8\\
                                      & $t_f$ [days]   &18.57 &39.10 &95.09 & 166.5\\
                                      & $\Delta V$ [m/s] &183.9 &155.1 &40.13 & 127.5\\
        \hline
    \end{tabular}
    \label{tab:mintimefuel}
\end{table}

\subsection{Minimum Fuel Time-Varying Phase-Fixed}

\begin{figure}[htbp!]
     \centering
     \begin{subfigure}{0.2672\textwidth}
         \centering
         \includegraphics[width=1\textwidth]{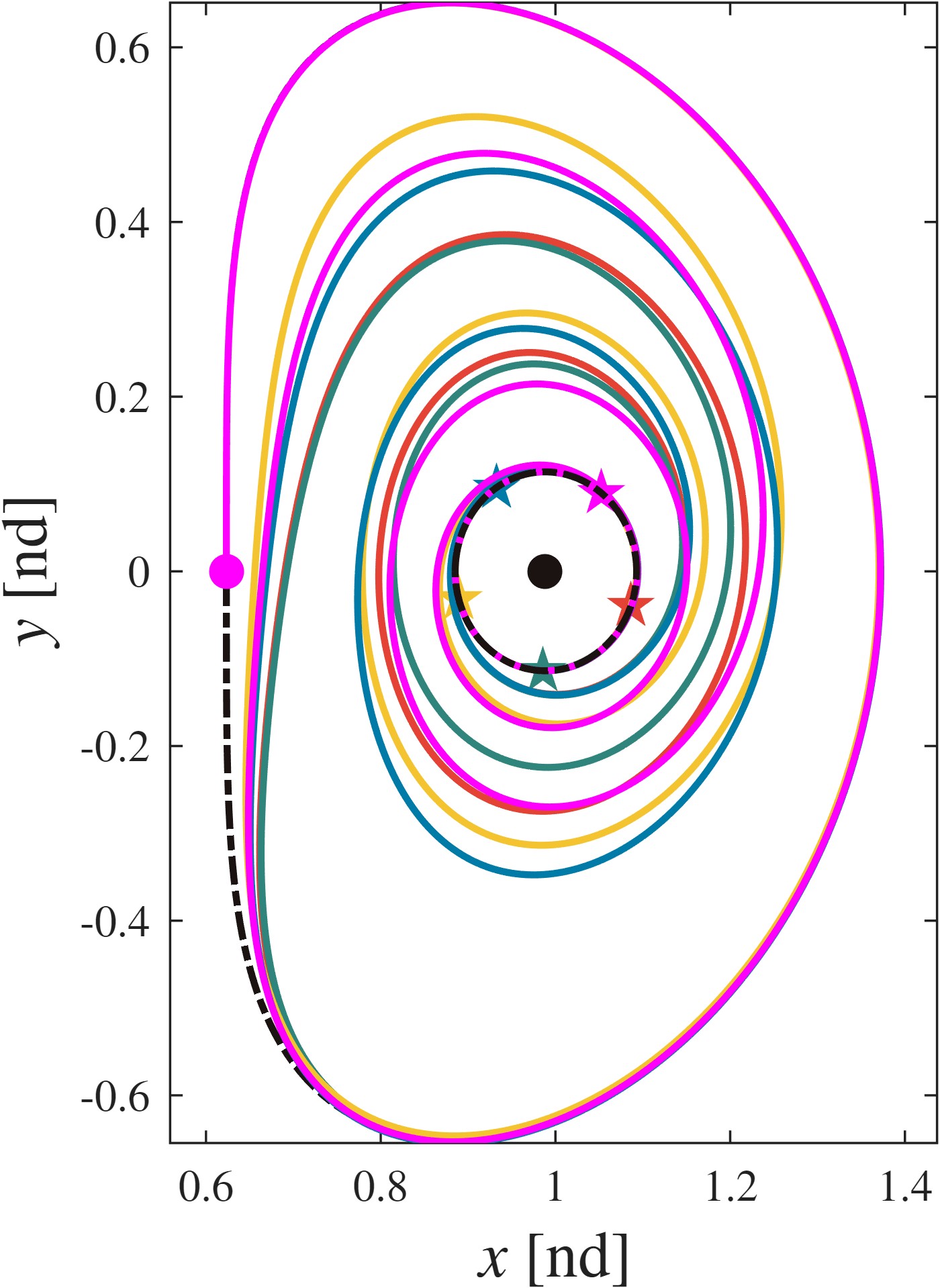}
         \caption{Earth-Moon Rotating}
     \end{subfigure}\hfill
     \begin{subfigure}{0.3563\textwidth}
         \centering
         \includegraphics[width=1\textwidth]{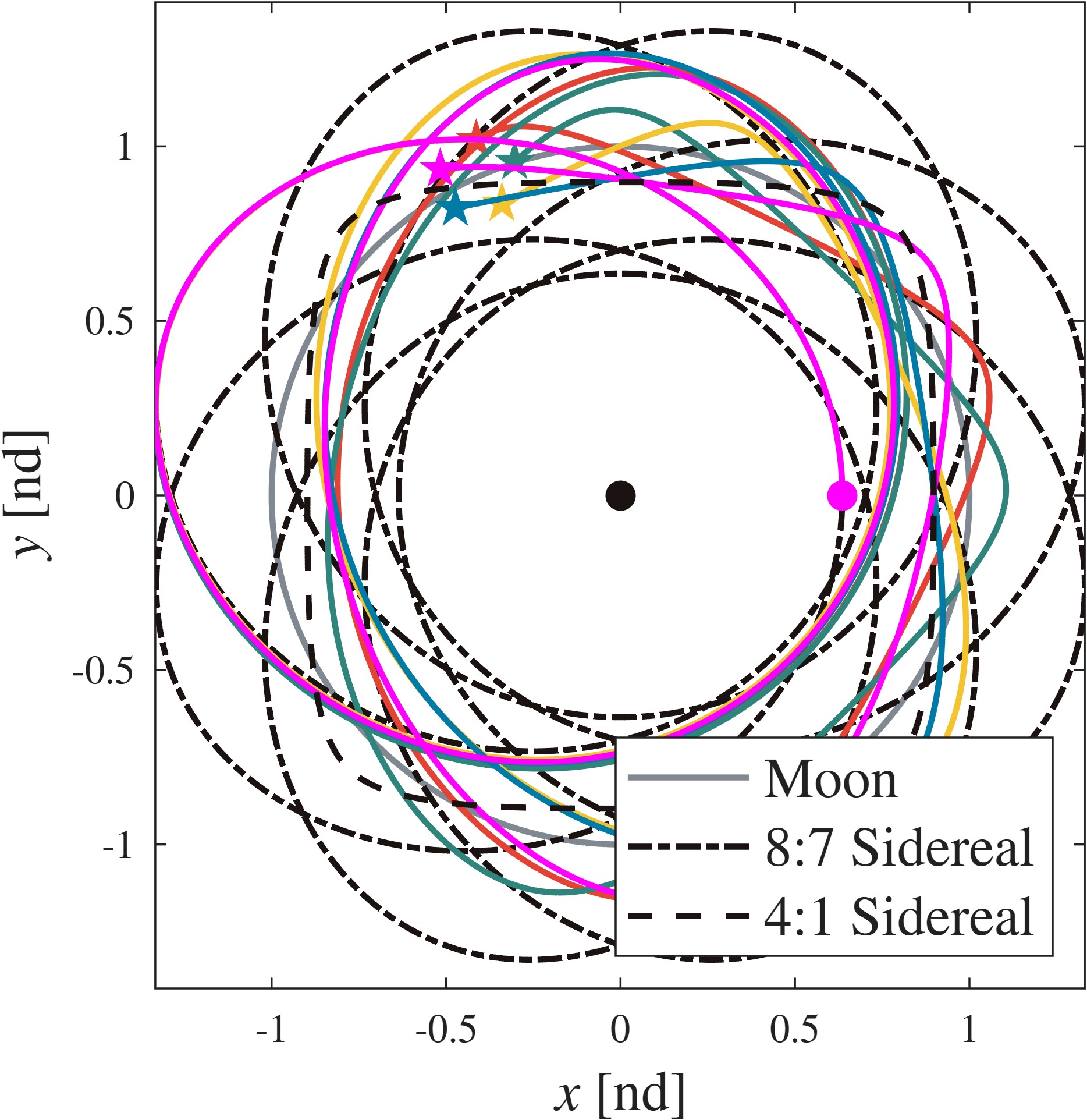}
         \caption{Earth-Centered Inertial}
     \end{subfigure}\hfill
     \begin{subfigure}{0.3563\textwidth}
         \centering
         \includegraphics[width=1\textwidth]{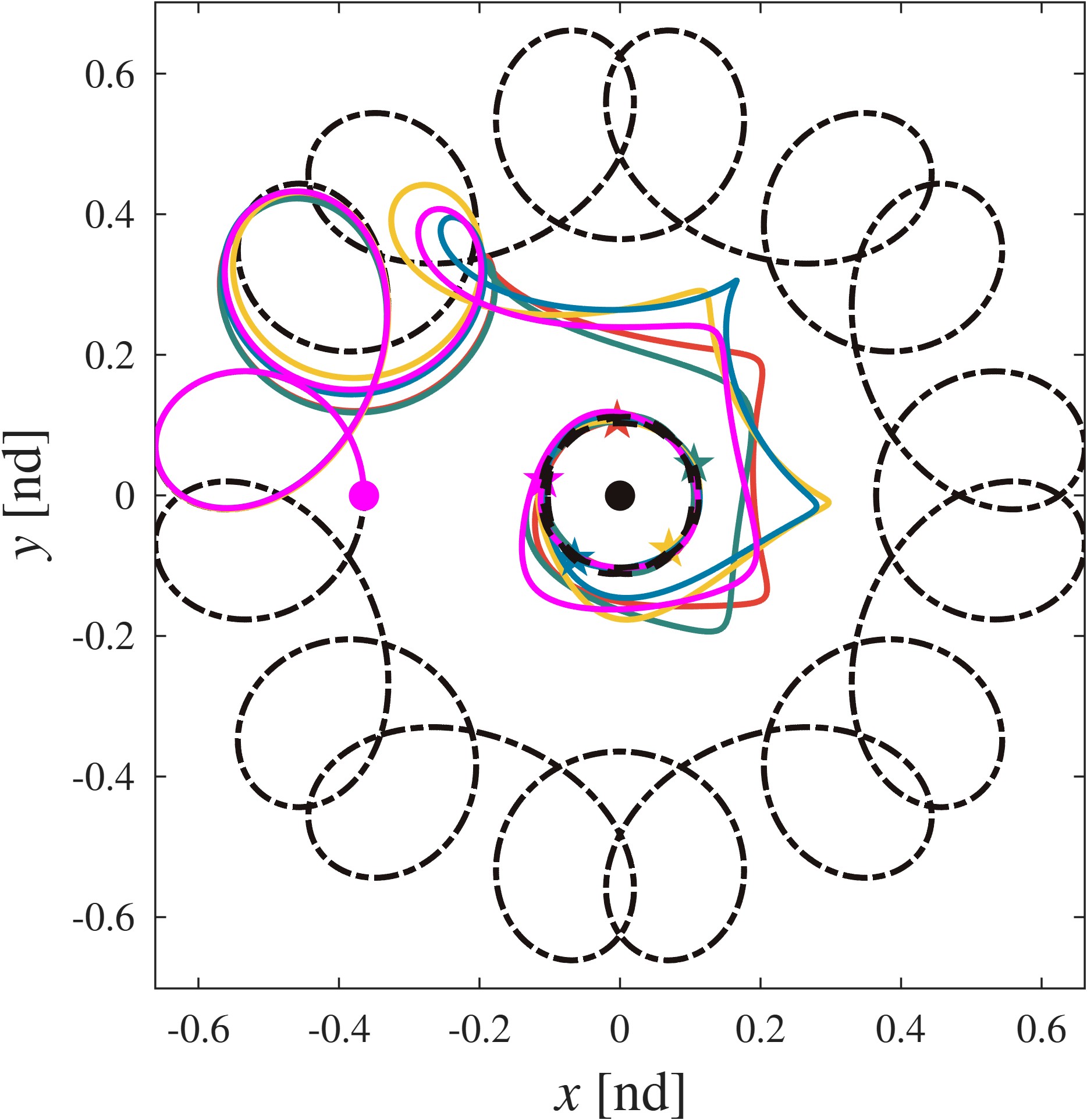}
         \caption{Moon-Centered Inertial}
     \end{subfigure}\\
     \begin{subfigure}{0.12\textwidth}
         \centering
         \includegraphics[width=1\textwidth]{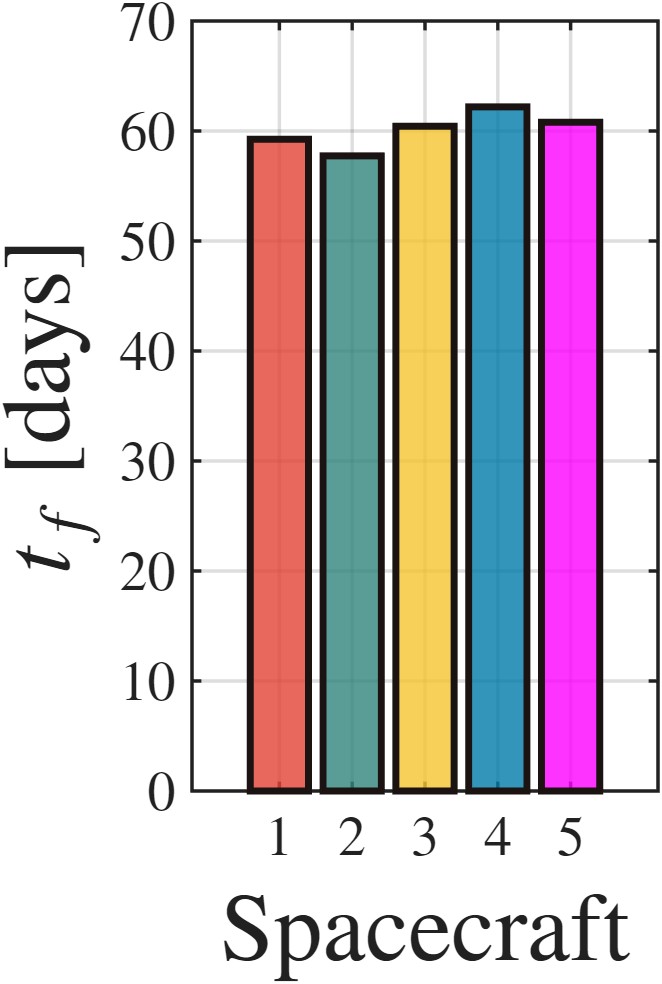}
         \caption{Time}
     \end{subfigure}\hfill
     \begin{subfigure}{0.12\textwidth}
         \centering
         \includegraphics[width=1\textwidth]{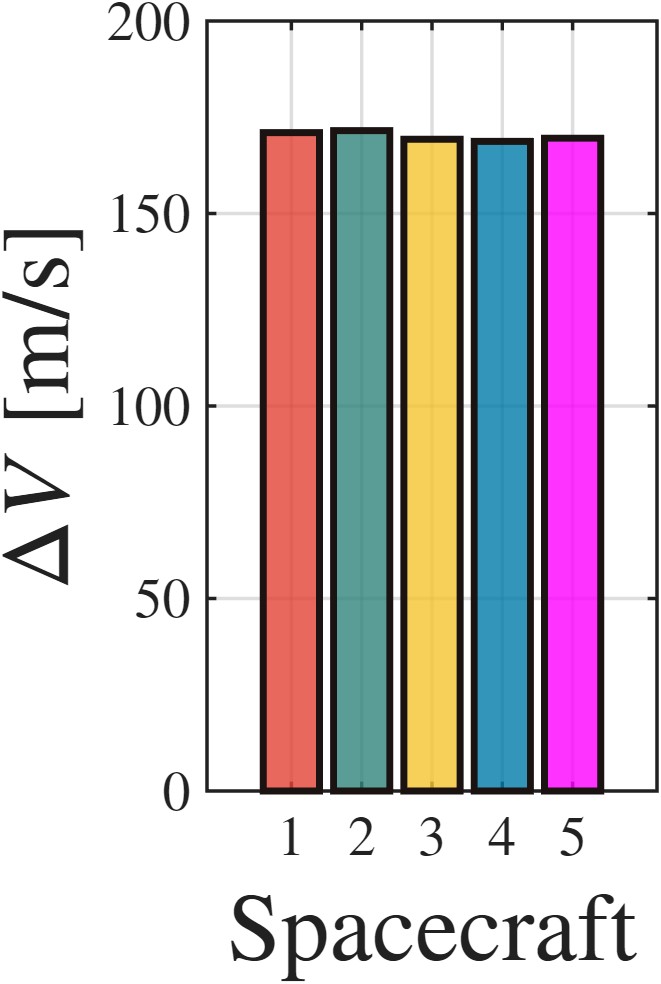}
         \caption{Cost}
     \end{subfigure}\hfill
     \begin{subfigure}{0.365\textwidth}
         \centering
         \includegraphics[width=1\textwidth]{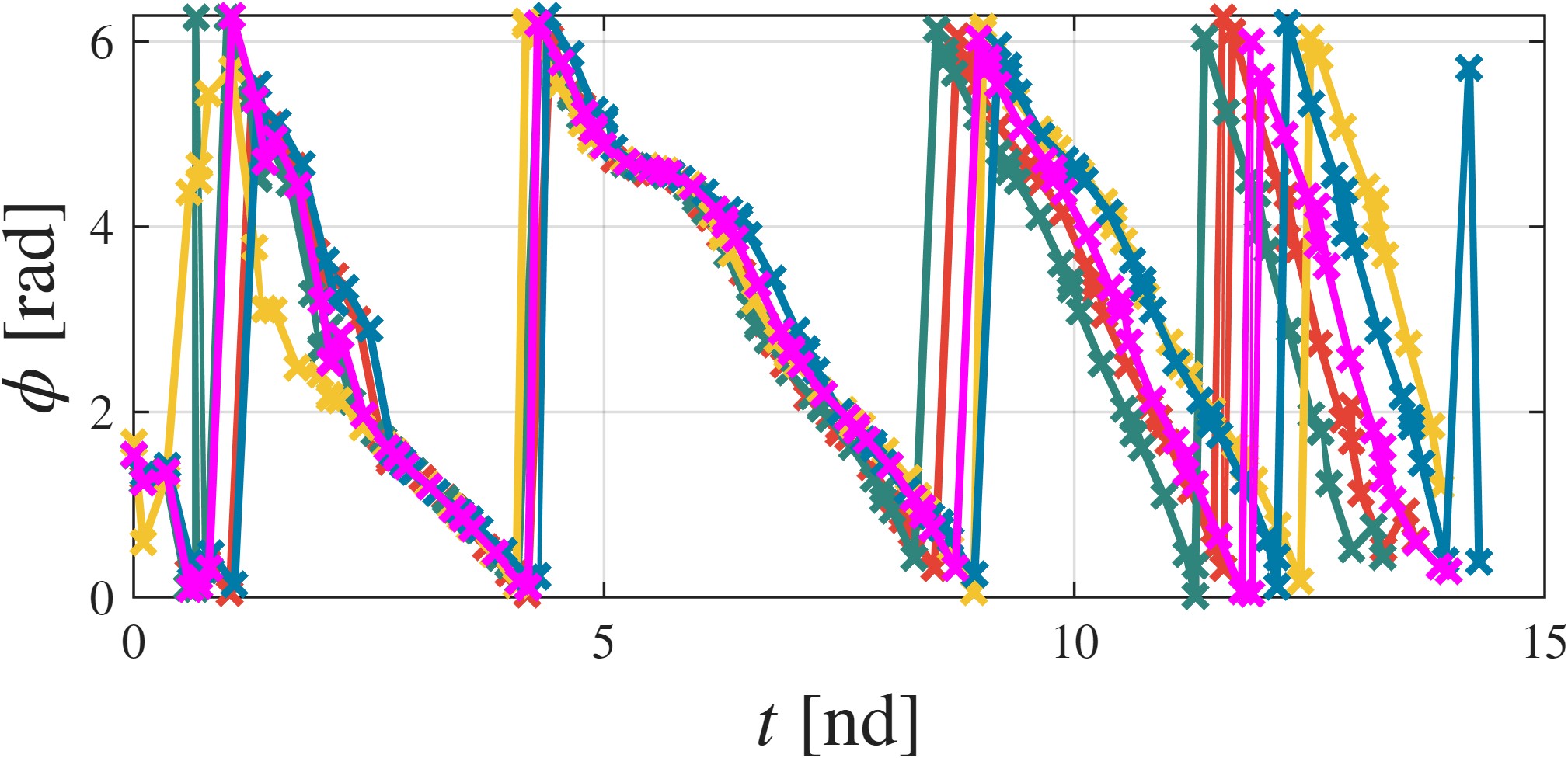}
         \caption{Angle Time History}
     \end{subfigure}\hfill
     \begin{subfigure}{0.37\textwidth}
         \centering
         \includegraphics[width=1\textwidth]{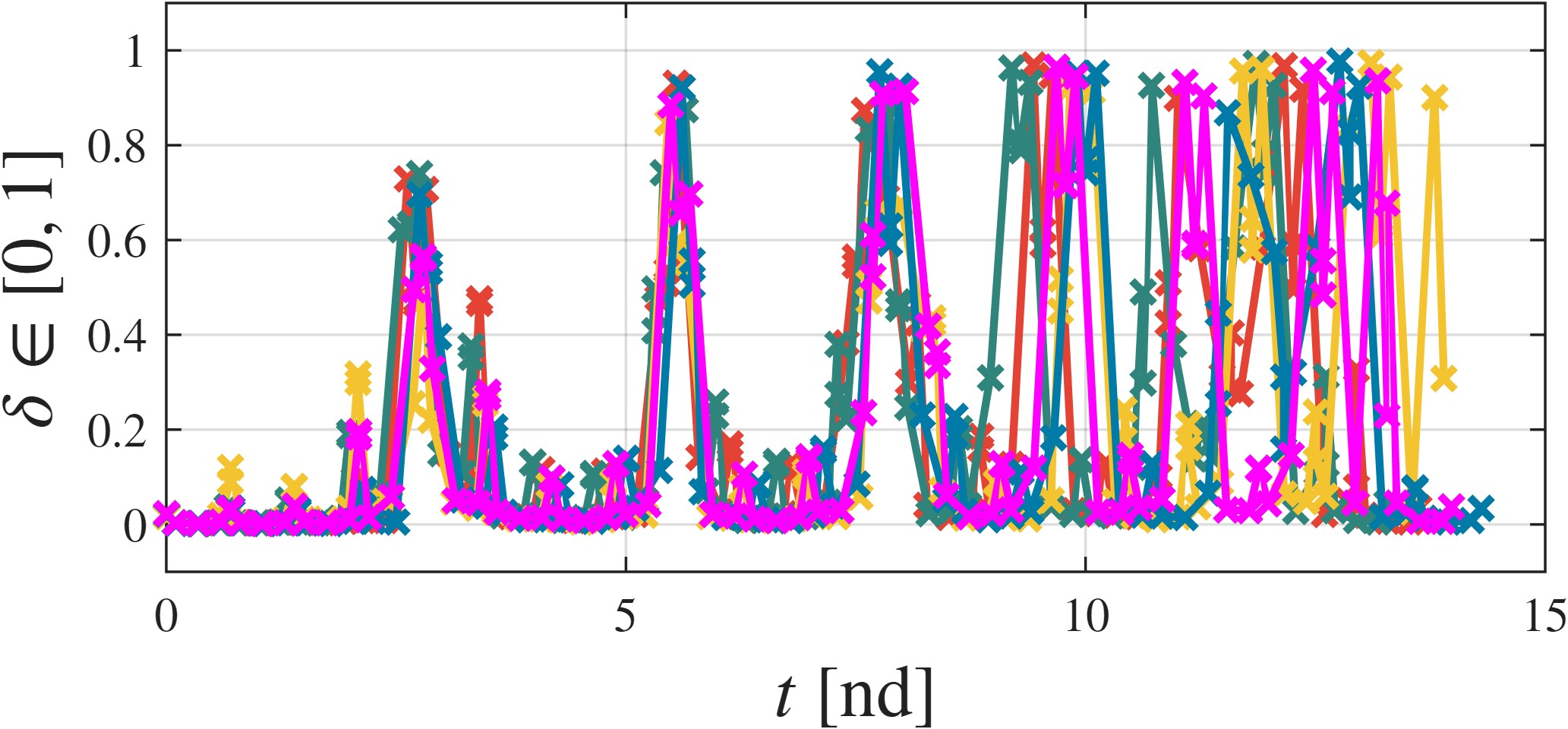}
         \caption{Throttle Time History}
     \end{subfigure}
     \caption{Minimum fuel time-varying fixed phase constellation deployment with $u_{\max}=0.06\,\,\text{[nd]}\approx0.1634$ mm/s$^2$.}
     \label{fig:const}
\end{figure}

Finally, a minimum fuel time-varying phase-fixed solution using the family space initialization scheme is now presented in the form of low-thrust constellation deployment. Let $\theta_i(t_0)$ be the desired phase of the $i$th spacecraft upon deployment from an orbital transfer vehicle, or launch vehicle upper stage. Then the terminal boundary condition for the transfer is fully parametrized by the free final time as
\begin{equation}
    \boldsymbol{x}(t_f) = \boldsymbol{\Gamma}(p_f,\omega(p_f)(t_{f,i}-t_0)+\theta_i(t_0))
\end{equation}
An equal angle distribution of five spacecraft are used to demonstrate this problem type. The results are shown in Fig. \ref{fig:const}.

The final insertion time for each spacecraft is varied, but stays around 60 days. Similarly, the fuel consumption is nearly equal between the spacecraft. The well-structured bang-off-bang throttle history shows that for the first three thrust arcs, the spacecraft are nearly in unison to depart and adjust their trajectories. Starting with the fourth thrust arc, the spacecraft begin to vary their burn times in accordance with their respective arrival phase. Motion is plotted in both the Earth-centered and Moon-centered inertial frames to highlight the non-Keplerian nature of transfers among DRO family members. Starting from family space initialization, through feasibility solution, phase-free solution, and time-varying phase-fixed solution for all spacecraft, this particular constellation deployment trajectory optimization problem took 6.532 minutes to run using six cores run in parallel for finite-differencing the nonlinear program Jacobians on a personal laptop with an Intel(R) Core(TM) Ultra 7 suite of processors. $N_\text{poly}=4$ and $N_\text{seg}=15$ were used for this $N_\text{rev}=3$ transfer. This showcases the computational efficiency of the family space to integral collocation family transfer framework developed in this work.

\section{Modifications for Time Regularization}

The distant retrograde orbit family was selected for the development above because it doesn't demand any of the following machinery. Its orbital frequency is monotonic across the family and can therefore serve directly as the family parameter, and its geometry in the rotating frame is uniformly well behaved, so a uniform sampling in time resolves every member equally well. Many families of practical interest fail one or both of these conditions. The $L_2$ halo family is the canonical counterexample: an NRHO has a dramatically larger dynamic range than the bifurcating Lyapunov-halo orbit, thus a uniform sampling in time for approximation does not work well for both of them; time regularization must be applied to approximate the family using a tractable number of data points, which can be dynamically placed for more coverage in fast regions.

\subsection{Regularized Dynamics and Family Function}

Following the Sundman-type transformations used for many-revolution transfer design\cite{aziz2018low} and for regularized family
interpolation\cite{bourjeili2025approximations}, a strictly positive scalar regularization function $h(\boldsymbol{x})$ is introduced such that
\begin{equation}
    \frac{\mathrm{d}t}{\mathrm{d}s} = h(\boldsymbol{x}) > 0
    \label{eq:reg:sundman}
\end{equation}
where $s$ is a new independent variable (fictitious time). Applying the chain rule, the dynamics are transformed
\begin{equation}
    \bm{x}' \;=\; \frac{\mathrm{d}\bm{x}}{\mathrm{d}s}
    \;=\; h(\bm{x})\Big[\bm{f}(\bm{x}) + \bm{g}(\bm{x})\bm{u}\Big]
    \label{eq:reg:dyn}
\end{equation}
where $(\cdot)'$ denotes differentiation with respect to $s$. The multi-body regularization function used in this work selected by trial and error to fit the $L_2$ halo family is
\begin{equation}
    h(\bm{x}) \;=\; d_1\cdot d_2^{3/2}
    \label{eq:reg:h}
\end{equation}
The exponent $3/2$ on the secondary distance is the classical Sundman choice for two-body-like periapsis passage, and it is what equalizes the sampling through perilune. The additional $d_1$ factor does not affect the sampling behavior but improves the conditioning of the Chebyshev system in Eq. \ref{eq:fit_system} by roughly a factor of $33$ at equal cost. Seven candidate functions of $(d_1, d_2)$ were compared in their invariance error residuals to select Eq. \ref{eq:reg:h}.

While collecting family data, a periodic orbit is now sampled at $m_\theta$ points equally spaced in $s$. Thus, a different family function is generated from the fit, as the torus angles advance according to their respective regularized orbital frequencies at each parameter (orbit) $p$.
\begin{equation}
    T_\text{reg}( p) = \oint \frac{\mathrm{d}t}{h(\bm{x})},
    \quad
    \omega_\text{reg}( p) = \frac{2\pi}{T_\text{reg}( p)}
    \label{eq:reg:omega}
\end{equation}
The regularized family dynamics become
\begin{equation}
    \bm{y}' = \bm{l}(\bm{y}) = \begin{bmatrix} 0 \\ \omega_\text{reg}( p)\end{bmatrix},
    \qquad
    \theta( p, s) = \omega_\text{reg}( p)(s - s_0) + \theta(s_0)
    \label{eq:reg:angle}
\end{equation}
and the reduced family invariance condition becomes
\begin{equation}
    \left[\frac{\partial \bm{\Gamma}}{\partial \theta}\right]\omega_\text{reg}( p)
    \;=\; h\big(\bm{\Gamma}(\bm{y})\big)\,\bm{f}\big(\bm{\Gamma}(\bm{y})\big)
    \label{eq:reg:invariance}
\end{equation}
The right-hand side is the only structural change: the field being matched is the regularized field $h\bm{f}$ rather than $\bm{f}$. The approximation error measure is likewise evaluated against the regularized field,
\begin{equation}
    \epsilon = \frac{1}{m_s m_p m_\theta}\sum_{k=1}^{m_s}\sum_{j=1}^{m_p}\sum_{i=1}^{m_\theta}
    \left\lVert \left[\frac{\partial \bm{\Gamma}}{\partial\tilde\theta_i}\right]
    \omega_\text{reg}\big( p_k(\tilde\tau_j)\big)
    - h\Big(\bm{\Gamma}\big( p_k(\tilde\tau_j),\tilde\theta_i\big)\Big)
      \bm{f}\Big(\bm{\Gamma}\big( p_k(\tilde\tau_j),\tilde\theta_i\big)\Big)\right\rVert
    \label{eq:reg:err}
\end{equation}

\subsection{Control Mapping and Initialization Under Regularization}
Assuming the family parameter's units have been normalized (if needed), the family function sensitivity matrices $\bm{T}$ and $\dot{\bm{T}}$ are built entirely from partial derivatives of the family function with respect to its own arguments. They are properties of the geometry of $\bm{\Gamma}$ and carry no dependence on the independent variable used to traverse it. The velocity slip penalty function $\alpha(\bm{y})$ is therefore \emph{unchanged} by regularization. The only modification required anywhere in the control mapping is the
conversion into the variable actually being integrated. Since
$\mathrm{d}t = h\,\mathrm{d}s$,
\begin{equation}
    \nu_{s,\max}(\bm{y}) = h\big(\bm{\Gamma}(\bm{y})\big)\frac{\alpha(\bm{y})\,u_{\max}}
    {\sqrt{\hat{\bm{\nu}}^\mathrm{T}\dot{\bm{T}}^\mathrm{T}\dot{\bm{T}}\hat{\bm{\nu}}}}
    \label{eq:reg:numax}
\end{equation}

Physical time is no longer the independent variable, so it must be carried explicitly as an augmented state in the initialization scheme to seed the phase space optimization problem. With the control direction $\hat{\bm{\nu}} = [\,\text{sign}(\bar p_f - \bar p_0),\; 0\,]^\mathrm{T}$, the family space initialization dynamics become
\begin{equation}
    \bm{y}' = \begin{bmatrix} 0 \\ \omega_\text{reg}(\bar p)\end{bmatrix}
    + \hat{\bm{\nu}}\,\nu_{s,\max}(\bm{y}),
    \qquad
    t' = h\big(\bm{\Gamma}(\bm{y})\big)
    \label{eq:reg:guessode}
\end{equation}
and the targeting problem is posed in $s$:
\begin{equation}
    \text{find } s_f: \;\;
    \bar p_0 + \int_{s_0}^{s_f}\text{sign}(\bar p_f - \bar p_0)\,\nu_{s,\max}(\bm{y})\,\mathrm{d}s = \bar p_f,
    \qquad
    t_f = \int_{s_0}^{s_f} h\big(\bm{\Gamma}(\bm{y})\big)\,\mathrm{d}s
    \label{eq:reg:target}
\end{equation}

\subsection{$L_2$ Halo Family}

\begin{figure}[htpb!]
    \centering
    \begin{subfigure}{0.357\textwidth}
         \centering
         \includegraphics[width=1\textwidth]{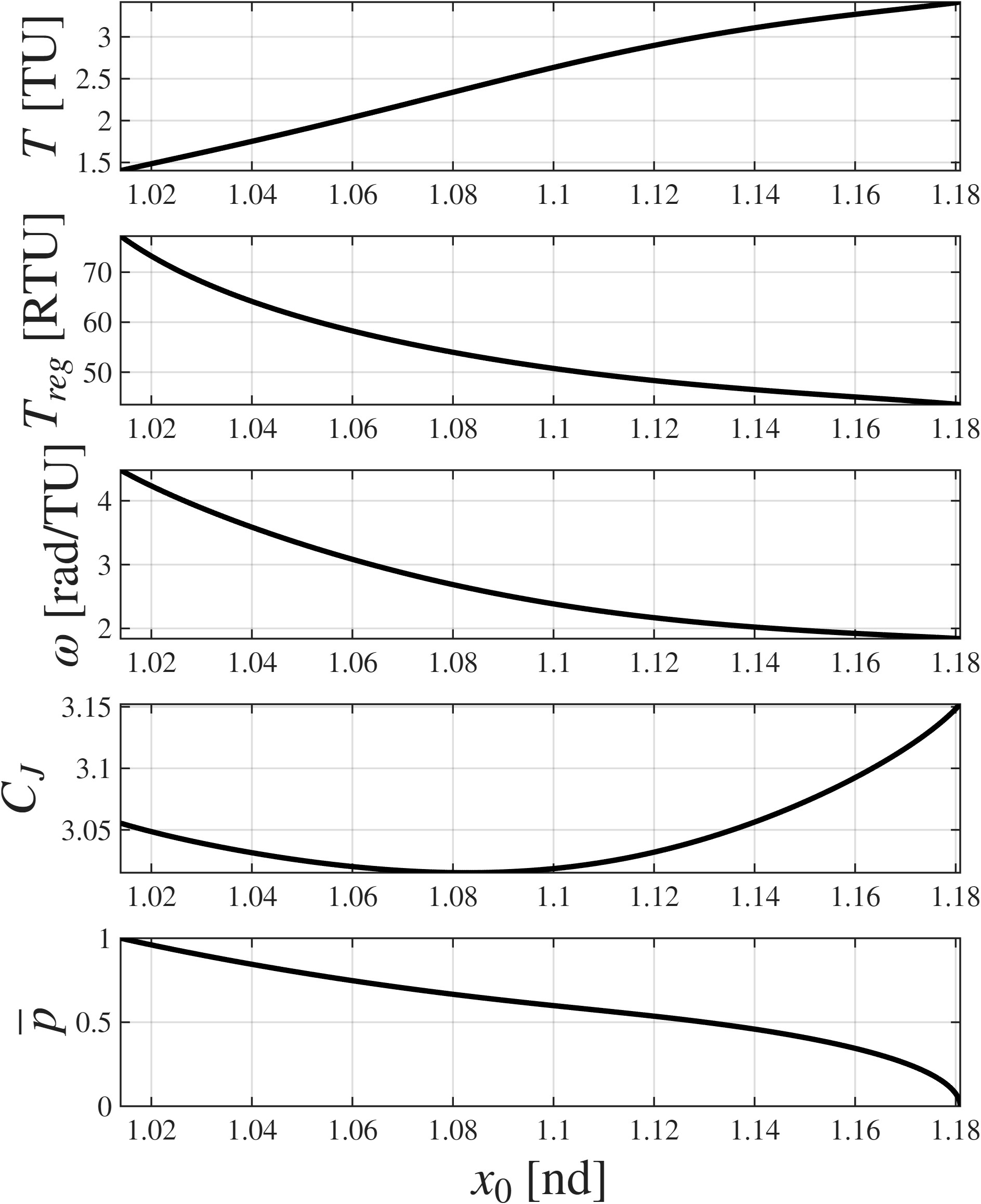}
         \caption{Family Data}
     \end{subfigure}\hfill
     \begin{subfigure}{0.233\textwidth}
         \centering
         \includegraphics[width=1\textwidth]{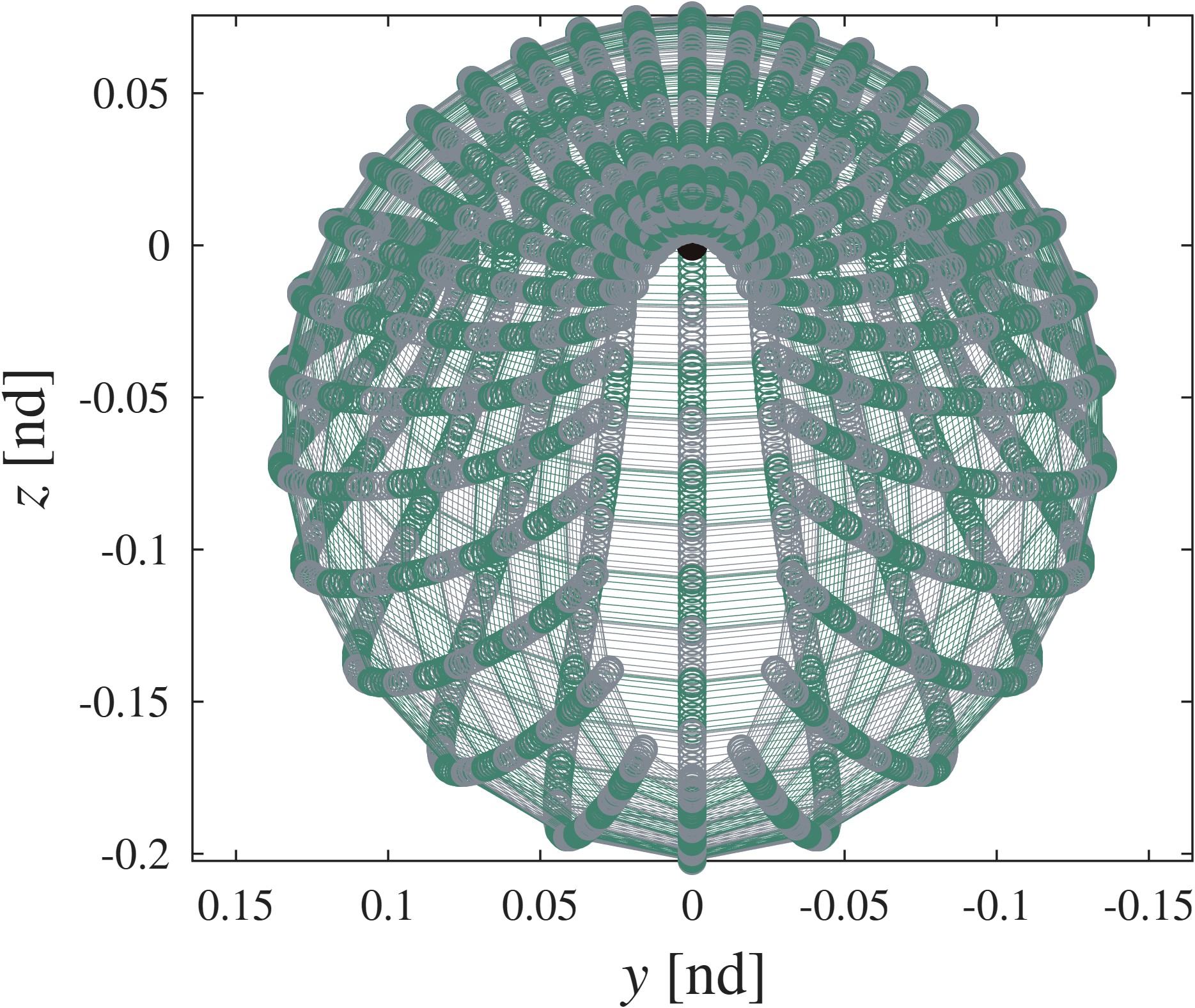}
         \includegraphics[width=1\textwidth]{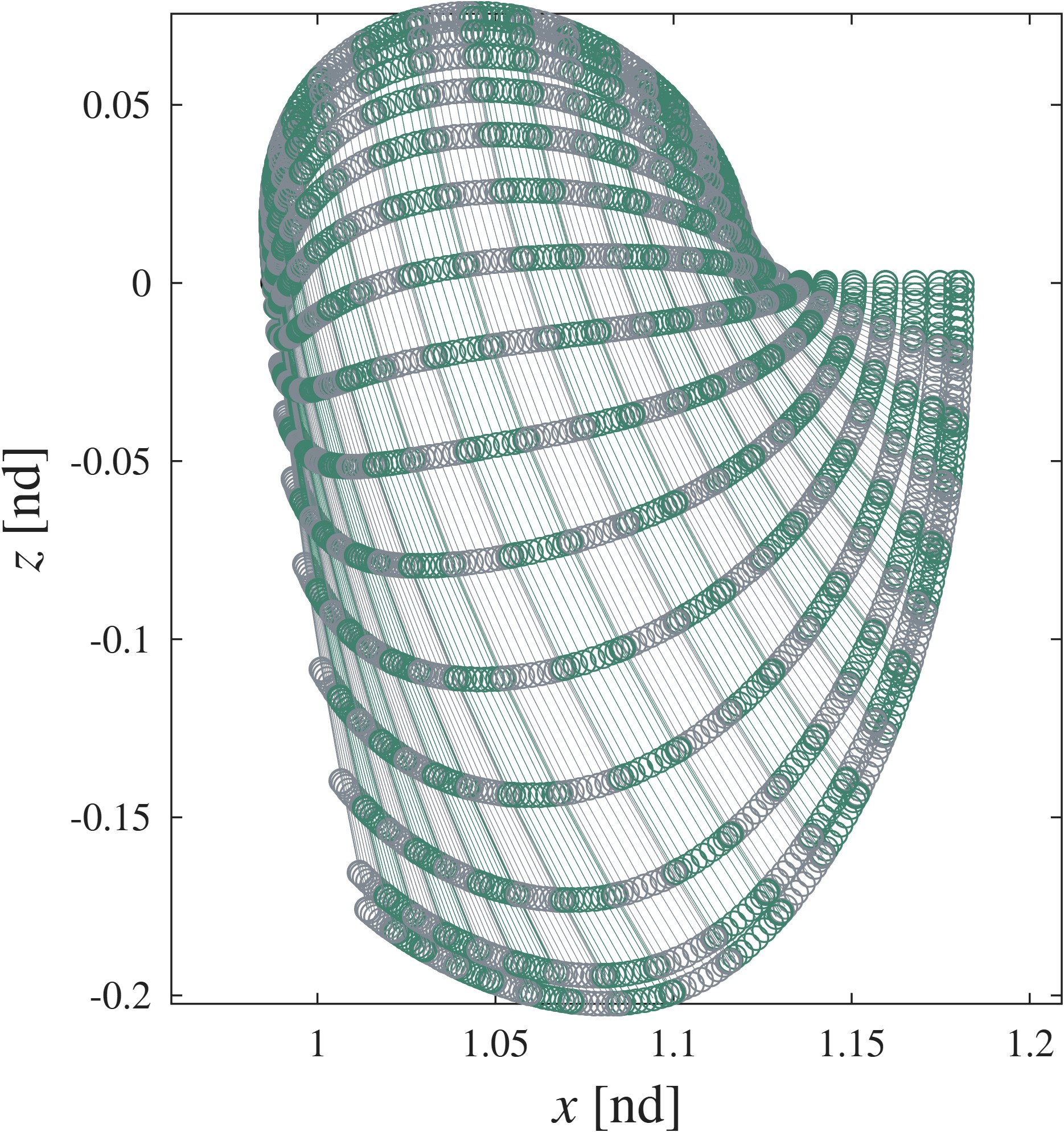}
         \caption{Family Fit}
     \end{subfigure}\hfill
     \begin{subfigure}{0.357\textwidth}
         \centering
         \includegraphics[width=1\textwidth]{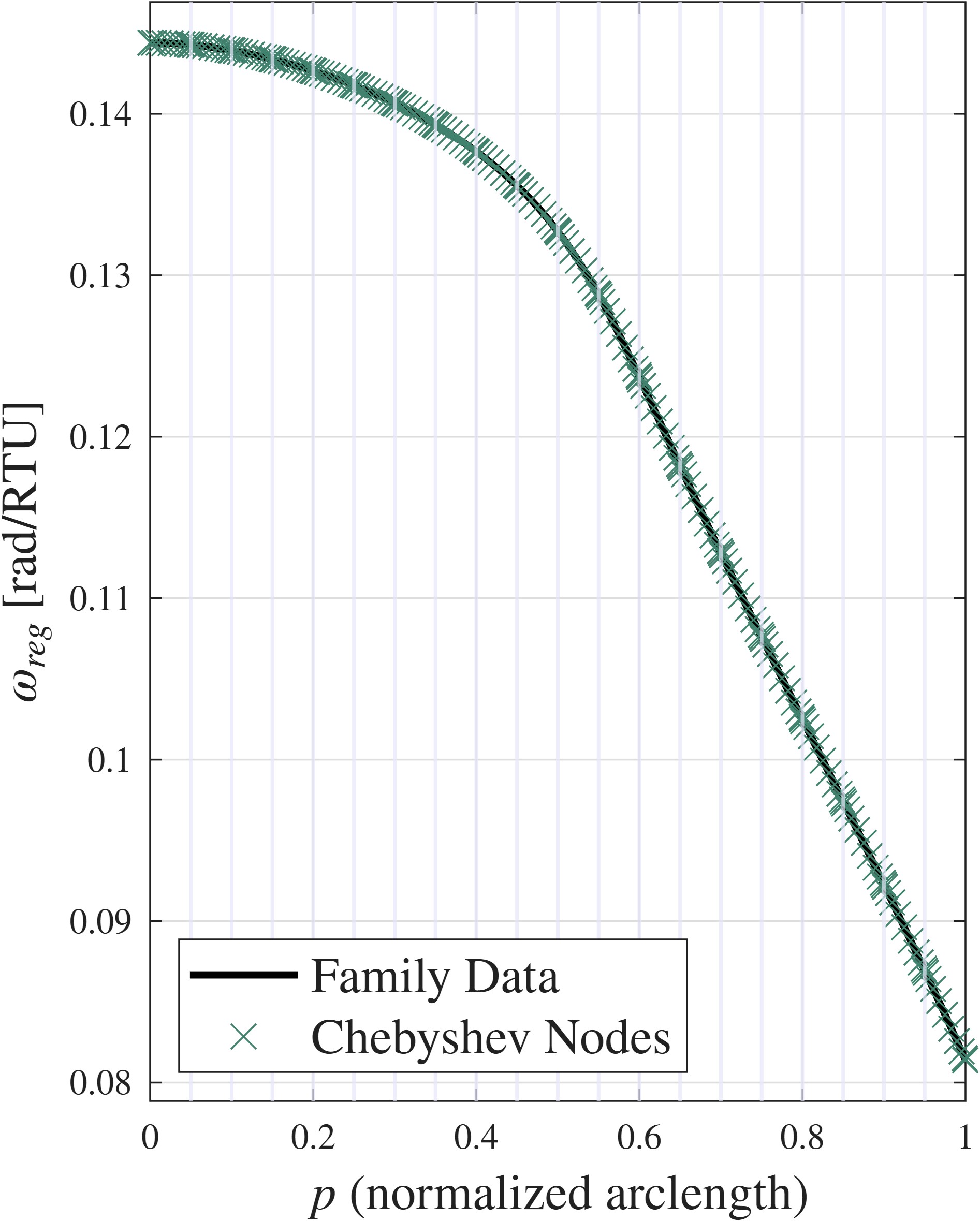}
         \caption{Regularized Frequency Fit}
     \end{subfigure}
    \caption{$L_2$ halo orbit family: characteristics and fit using arc-length as the parameter. Fit shown with $m_s=20$, $m_p=10$, $m_\theta=30$, and $\beta=1$, yielding a global invariance error $\epsilon=6.1055\times10^{-8}$ via time regularization. RTU = regularized time unit.}
    \label{fig:halo_fam_fit}
\end{figure}

The $L_2$ halo family is now used to demonstrate the regularized approximation approach. The normalized state space family arclength $\bar p$ is used as the family parameter with the recognition that parameter monotonicity is a necessary but not sufficient condition. At the $L_2$ halo pitchfork bifurcation, the out-of-plane amplitude grows as the square root of family distance, so any parameter even in that amplitude is stationary there and $\partial\bm\Gamma/\partial\bar p$ is unbounded. Since $d_1,d_2$ depend on $z$ only through $z^2$, this holds for every $\omega_\text{reg}$ generated by $h(d_1,d_2)$, hence arclength is needed. Let $\bm{x}_{0}(\sigma)$ denote the perpendicular crossing state as a function of any monotonic family index, then
\begin{equation}
    \bar p(\sigma) = \frac{1}{L}\int_{\sigma_0}^{\sigma}
    \left\lVert \frac{\mathrm{d}\bm{x}_{0}}{\mathrm{d}\sigma'}\right\rVert \mathrm{d}\sigma'\in[0,1],
    \quad
    L = \int_{\sigma_0}^{\sigma_f}
    \left\lVert \frac{\mathrm{d}\bm{x}_{0}}{\mathrm{d}\sigma'}\right\rVert \mathrm{d}\sigma'
    \label{eq:reg:arc}
\end{equation}
which is evaluated on the $n$ supplied family members as
\begin{equation}
    \bar p_k = \frac{1}{L}\sum_{i=1}^{k-1}\left\lVert \bm{x}_{0,i+1} - \bm{x}_{0,i}\right\rVert_2,
    \qquad
    L = \sum_{i=1}^{n-1}\left\lVert \bm{x}_{0,i+1} - \bm{x}_{0,i}\right\rVert_2
    \label{eq:reg:arcdisc}
\end{equation}
resulting in $\bar p = 0$ at the near-planar end and
$\bar p = 1$ at the NRHO end of the $L_2$ halo family. No aggressive clustering is required because arclength is selected, thus $\beta=1$.

The resulting fit is shown in Fig.~6 with $m_s = 20$, $m_p = 10$, $m_\theta = 30$ and $\beta = 1$, yielding a global invariance error of $\epsilon = 6.1055\times10^{-8}$ by Eq. \ref{eq:reg:err}. Figure \ref{fig:halo_fam_fit} shows the family characteristics including regularized frequency. Note that because a parameter other than frequency was chosen for the family, an additional fit of the frequency as a function of that parameter must be made in order to evaluate the initialization scheme in the family space.

\subsection*{Time-Free, Phase-Free Minimum Time and Fuel Solutions}

\begin{figure}[htbp!]
     \centering
     \begin{subfigure}{0.24\textwidth}
         \centering
         \includegraphics[width=1\textwidth]{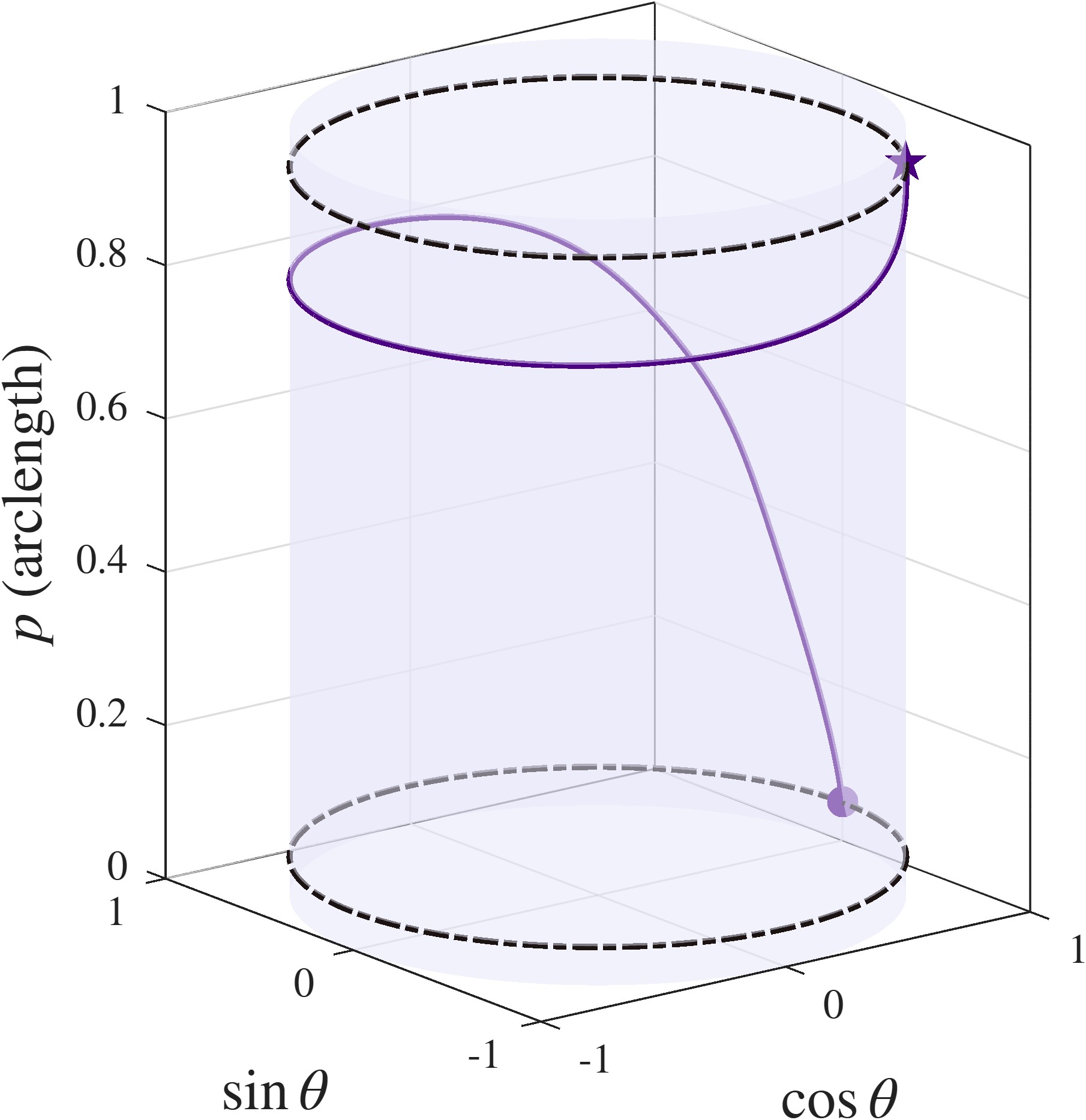}\\
         \includegraphics[width=1\textwidth]{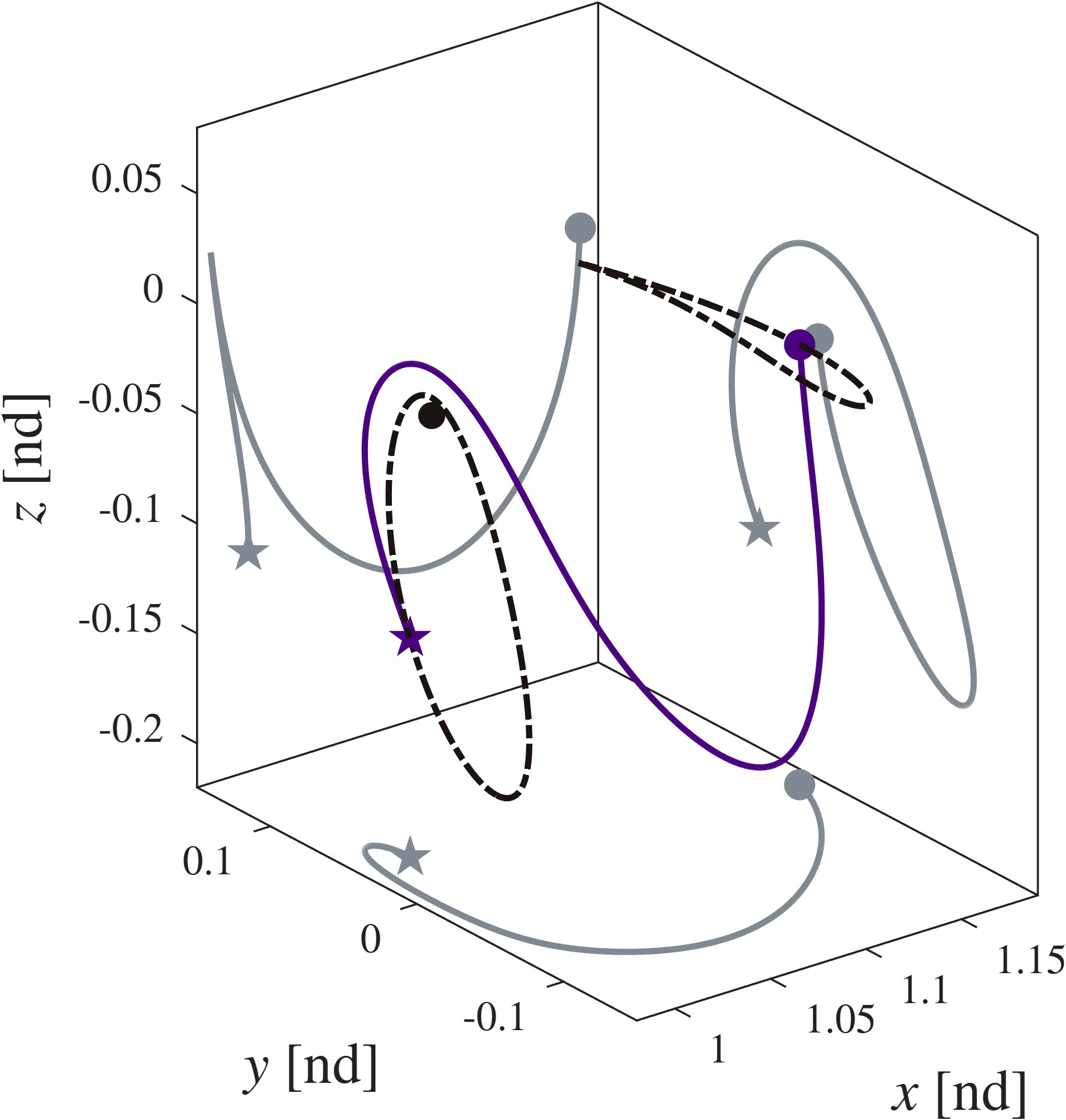}
         \caption{$u_{\max}=0.3$}
     \end{subfigure}\hfill
     \begin{subfigure}{0.24\textwidth}
         \centering
         \includegraphics[width=1\textwidth]{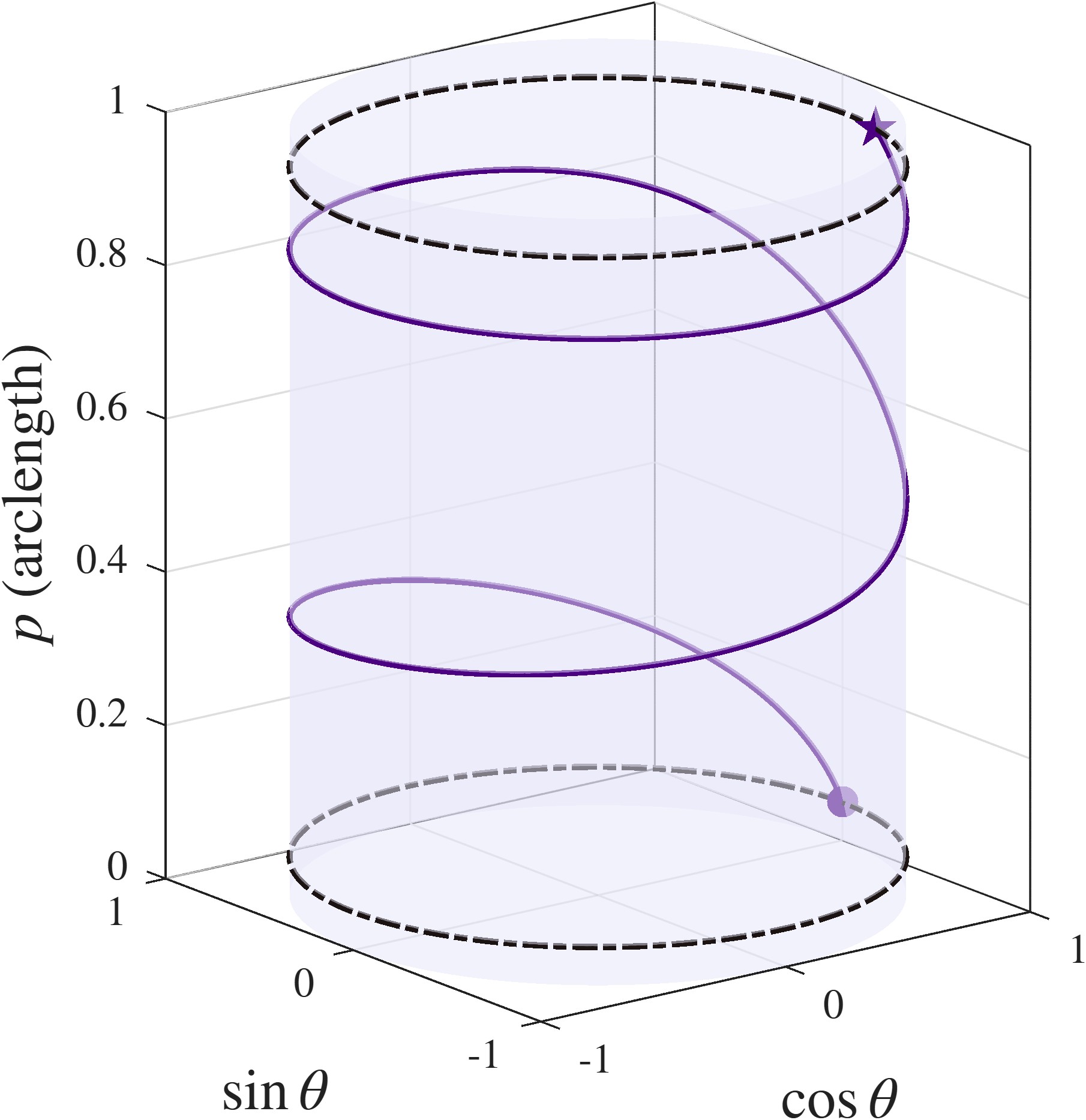}\\
         \includegraphics[width=1\textwidth]{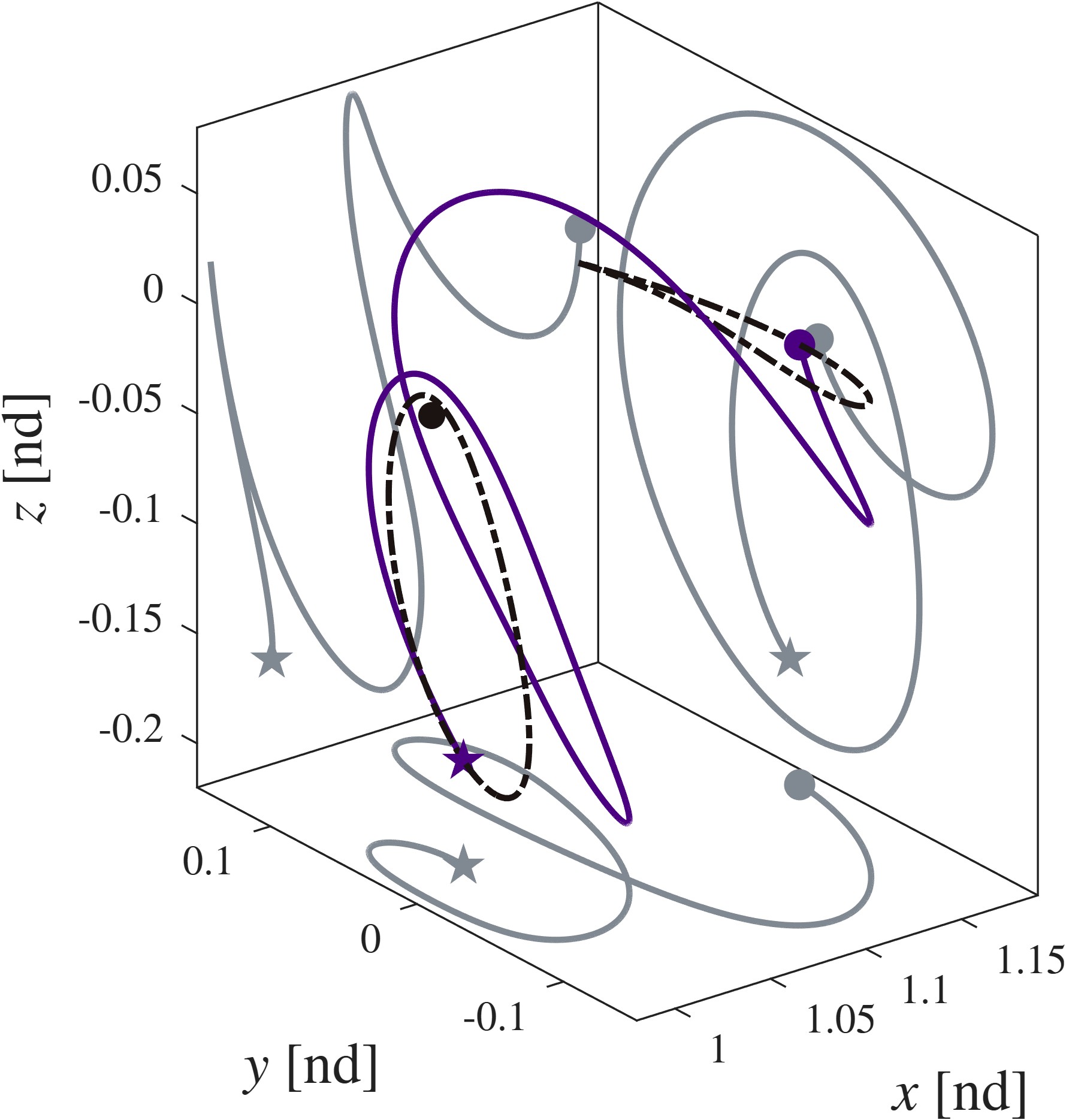}
         \caption{$u_{\max}=0.09$}
     \end{subfigure}\hfill
     \begin{subfigure}{0.24\textwidth}
         \centering
         \includegraphics[width=1\textwidth]{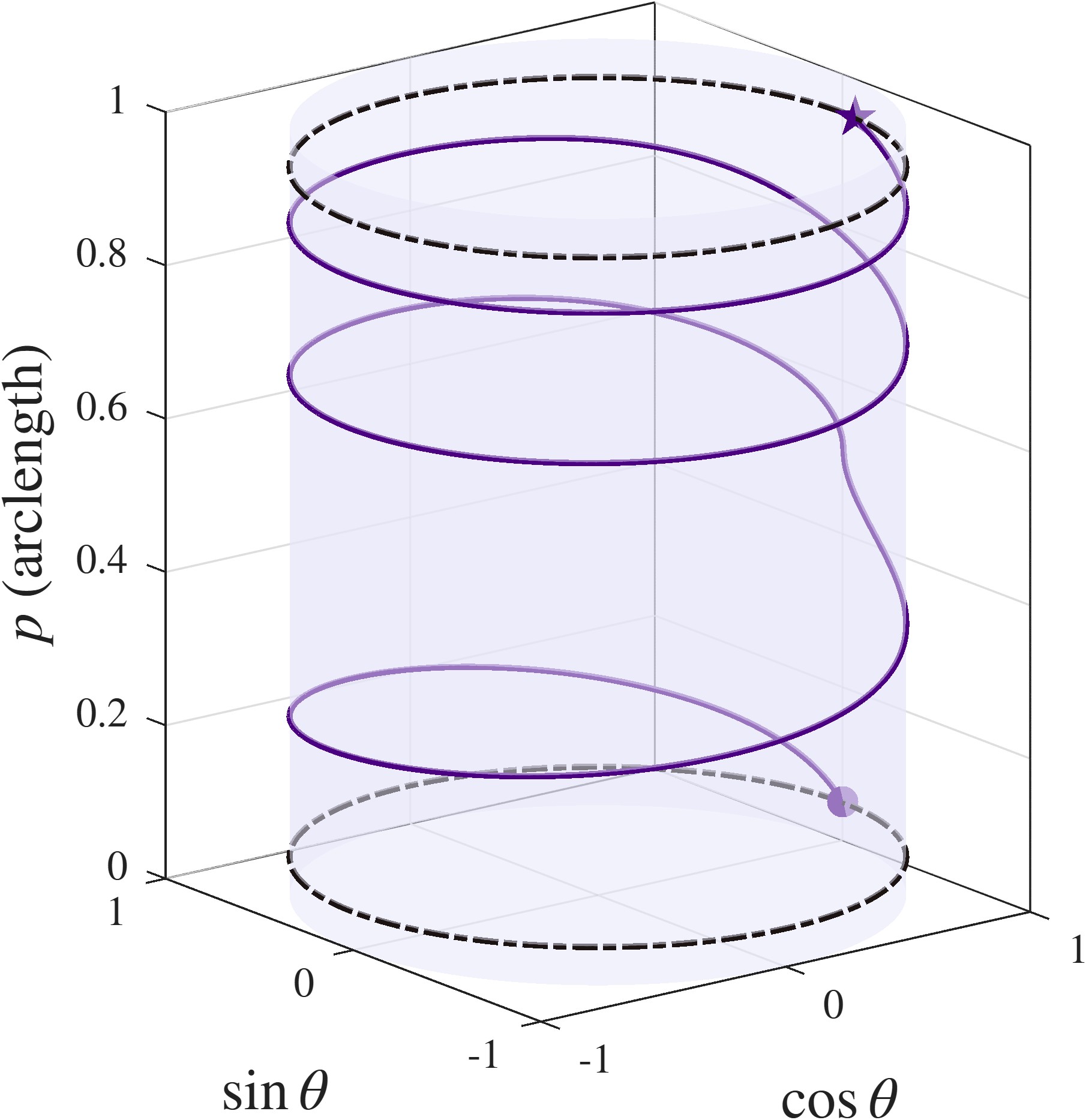}\\
         \includegraphics[width=1\textwidth]{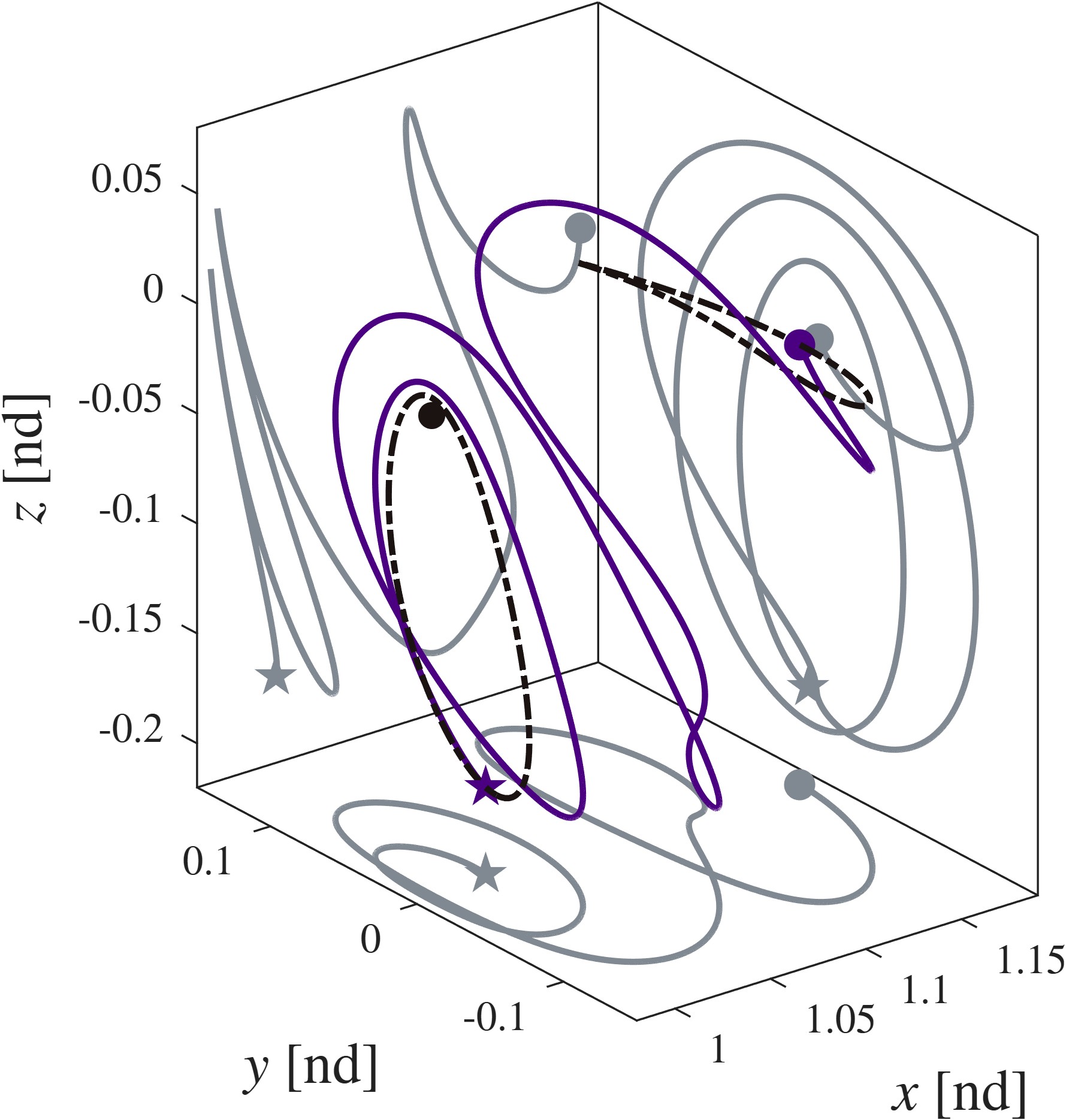}
         \caption{$u_{\max}=0.05$}
     \end{subfigure}\hfill
     \begin{subfigure}{0.24\textwidth}
         \centering
         \includegraphics[width=1\textwidth]{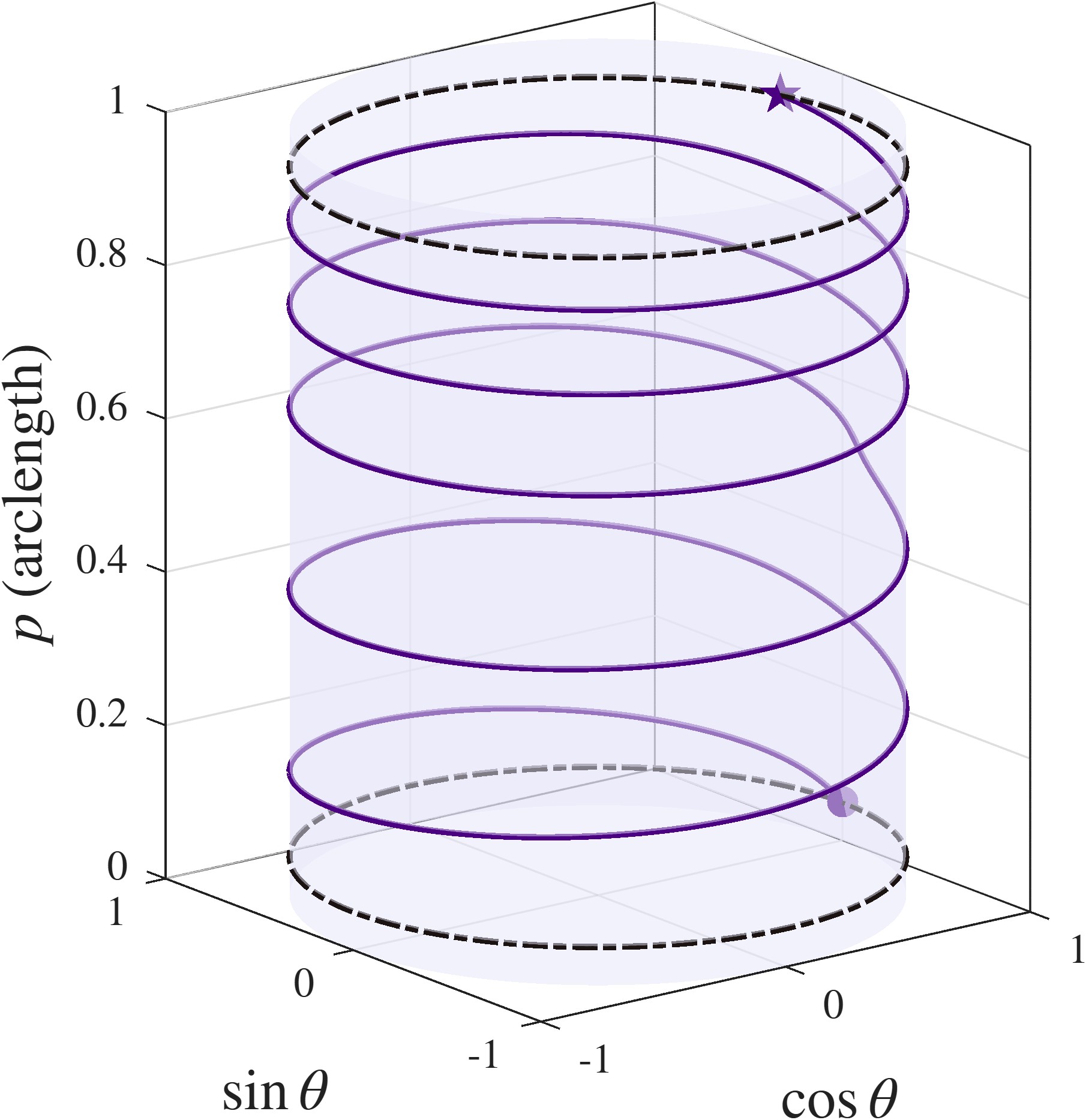}\\
         \includegraphics[width=1\textwidth]{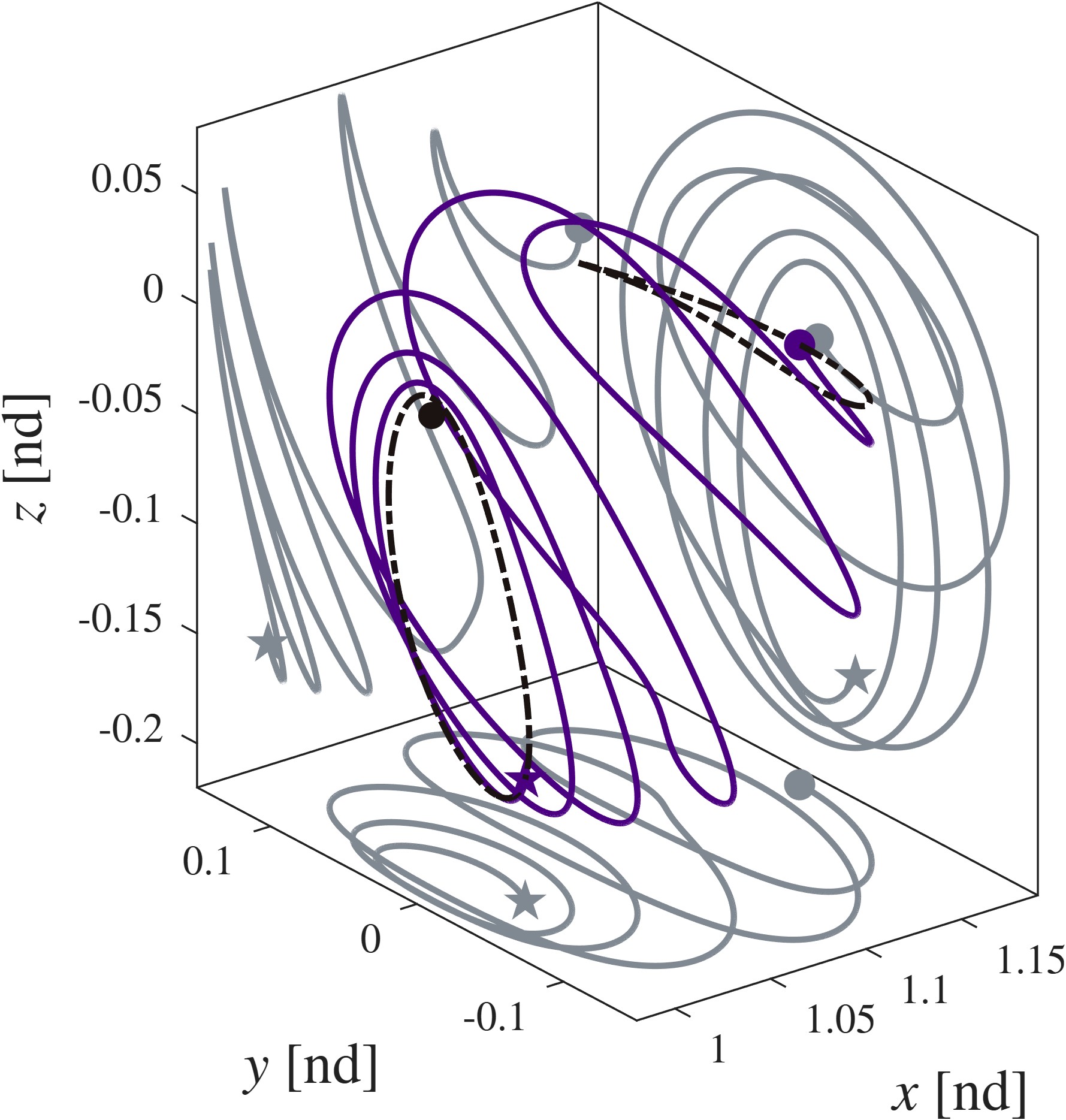}
         \caption{$u_{\max}=0.03$}
     \end{subfigure}
     \caption{Family space initialization trajectories for $L_2$ halo transfer.}
     \label{fig:halo_famspace}
\end{figure}

\begin{figure}[htbp!]
     \centering
     \begin{subfigure}{0.24\textwidth}
         \centering
         \includegraphics[width=1\textwidth]{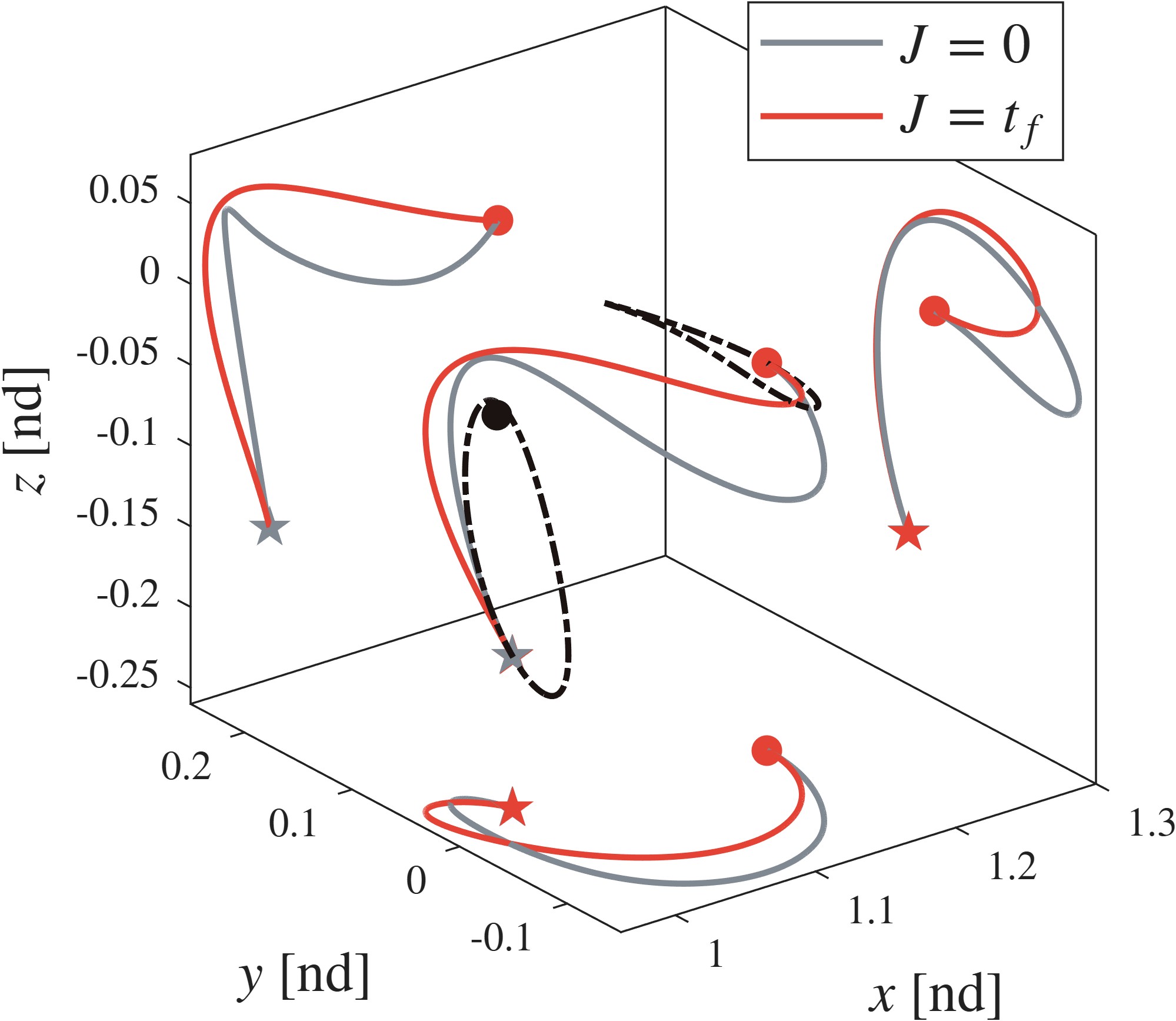}\\
         \includegraphics[width=1\textwidth]{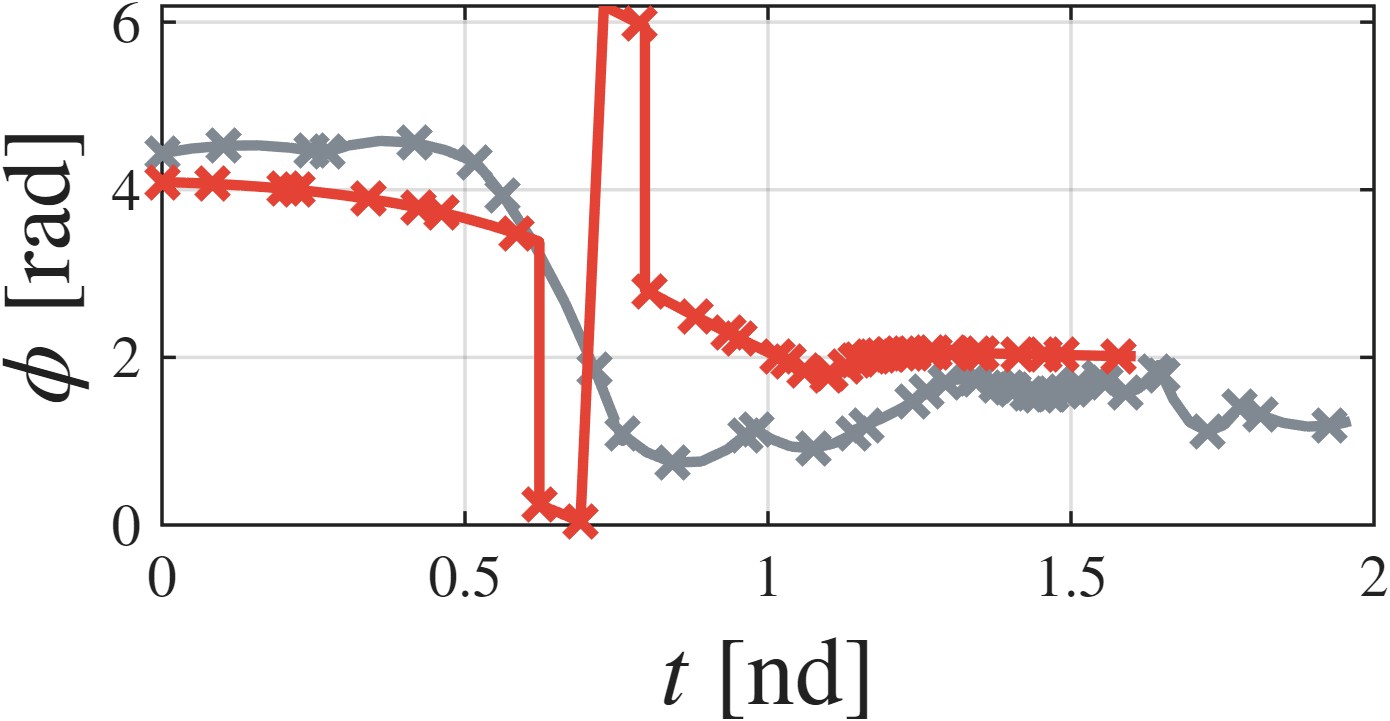}\\
         \includegraphics[width=1\textwidth]{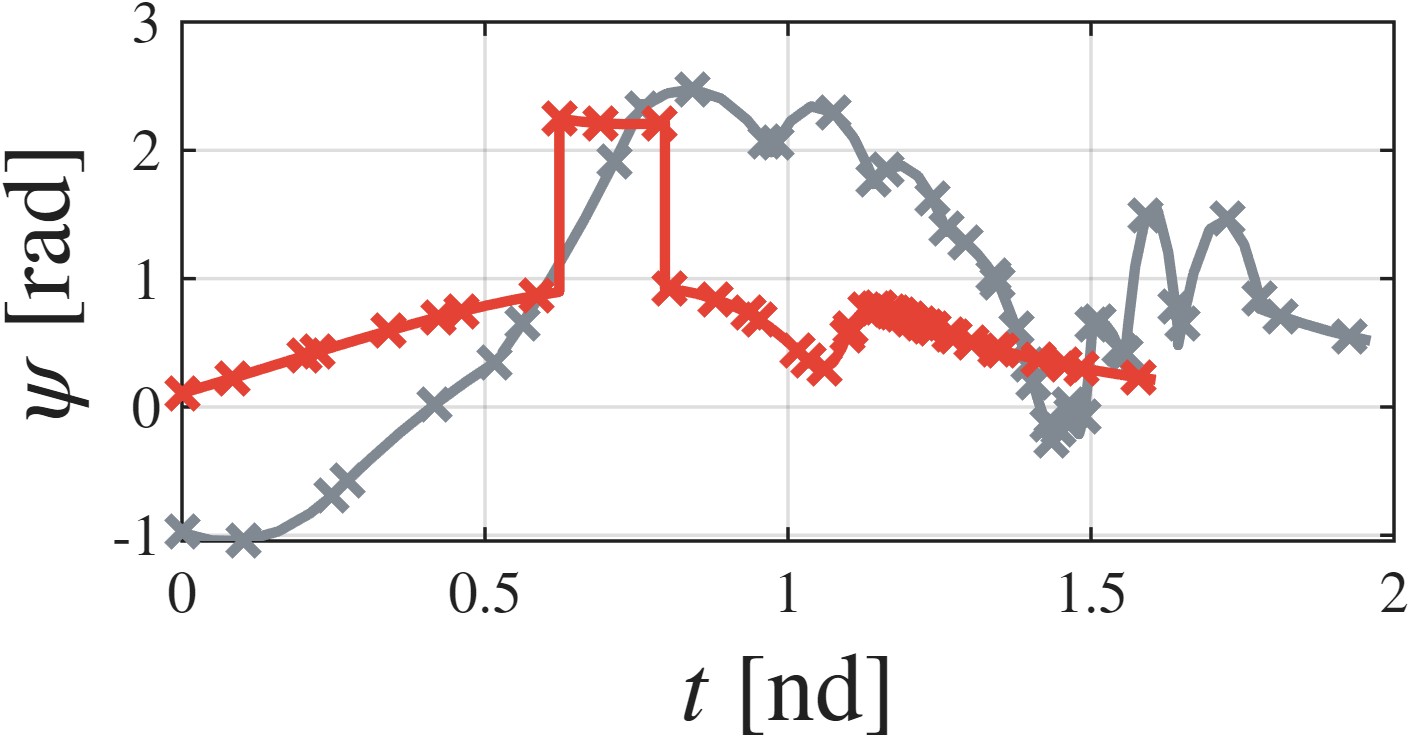}
         \caption{$u_{\max}=0.3$}
     \end{subfigure}\hfill
     \begin{subfigure}{0.24\textwidth}
         \centering
         \includegraphics[width=1\textwidth]{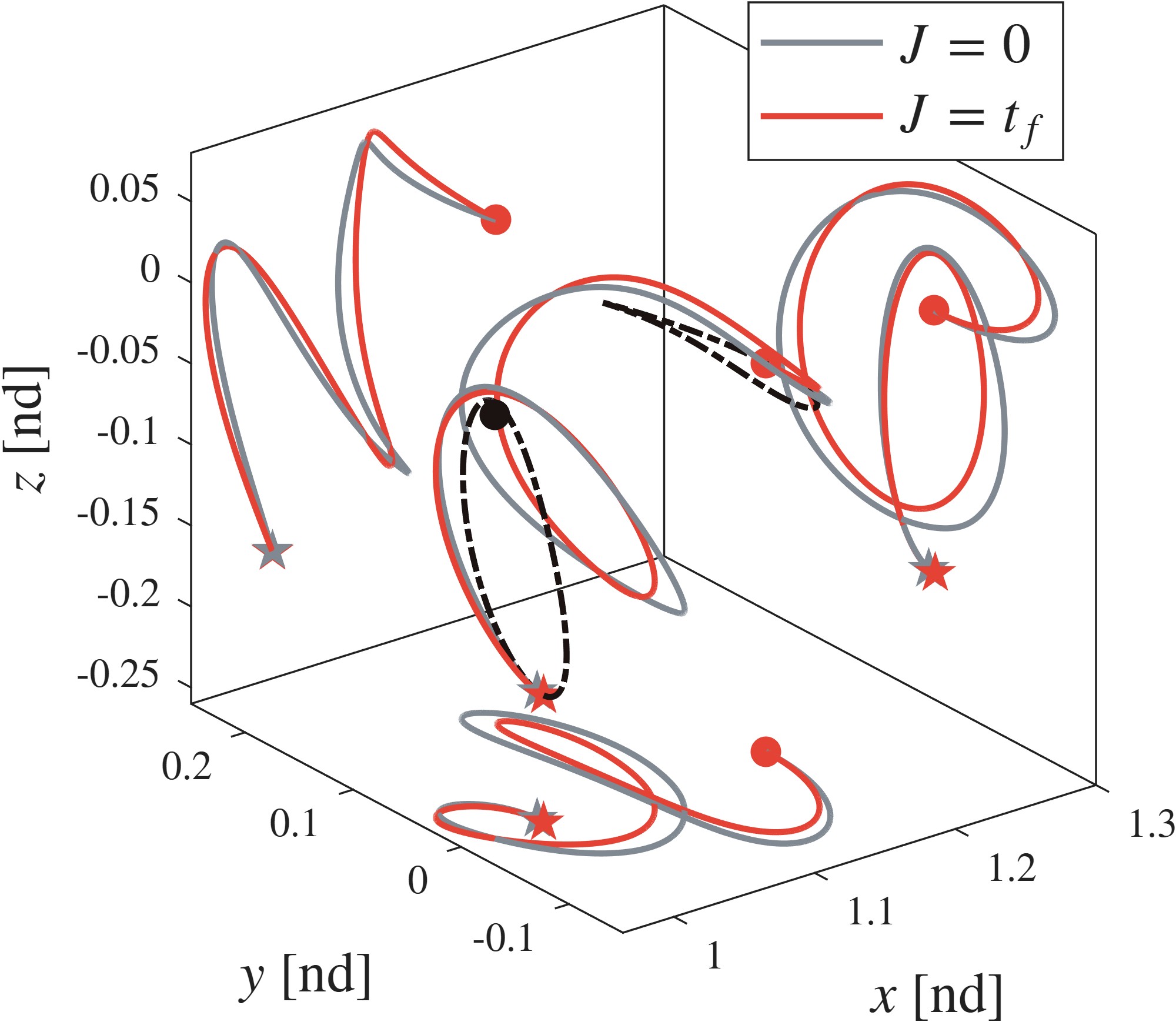}\\
         \includegraphics[width=1\textwidth]{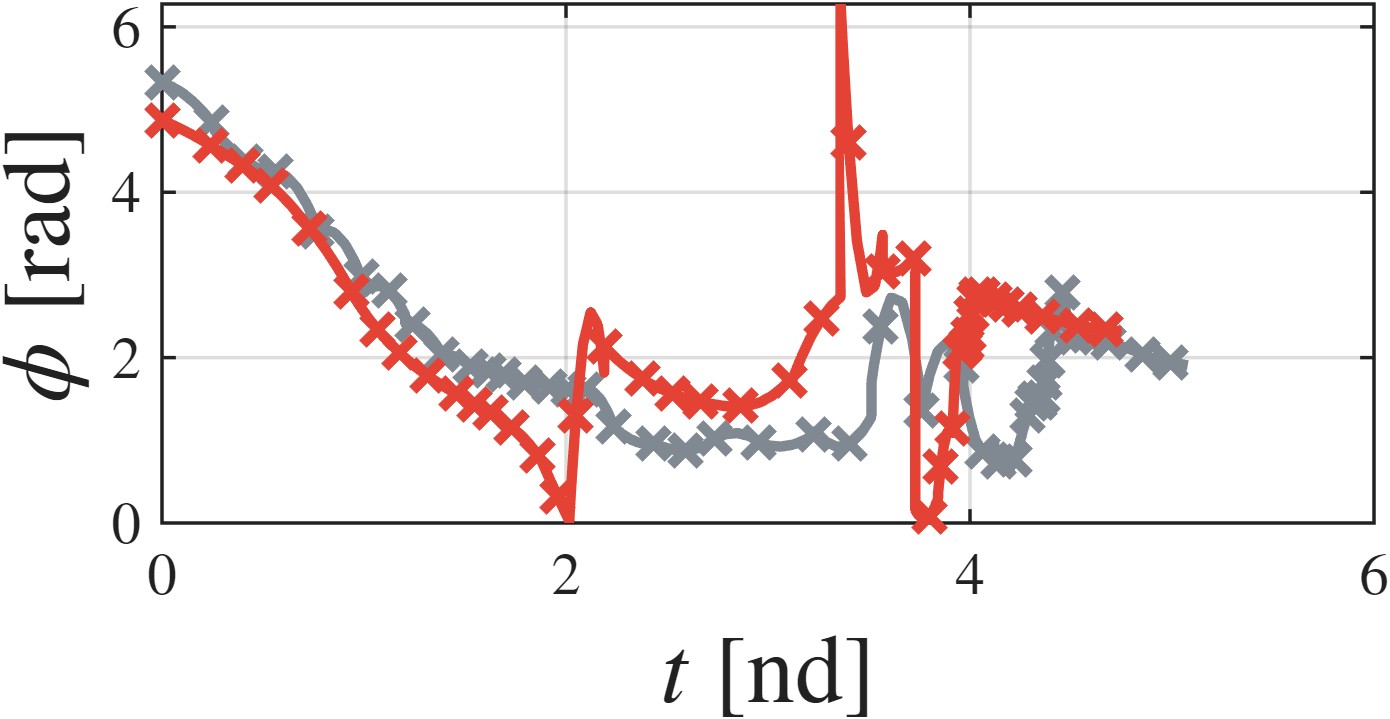}\\
         \includegraphics[width=1\textwidth]{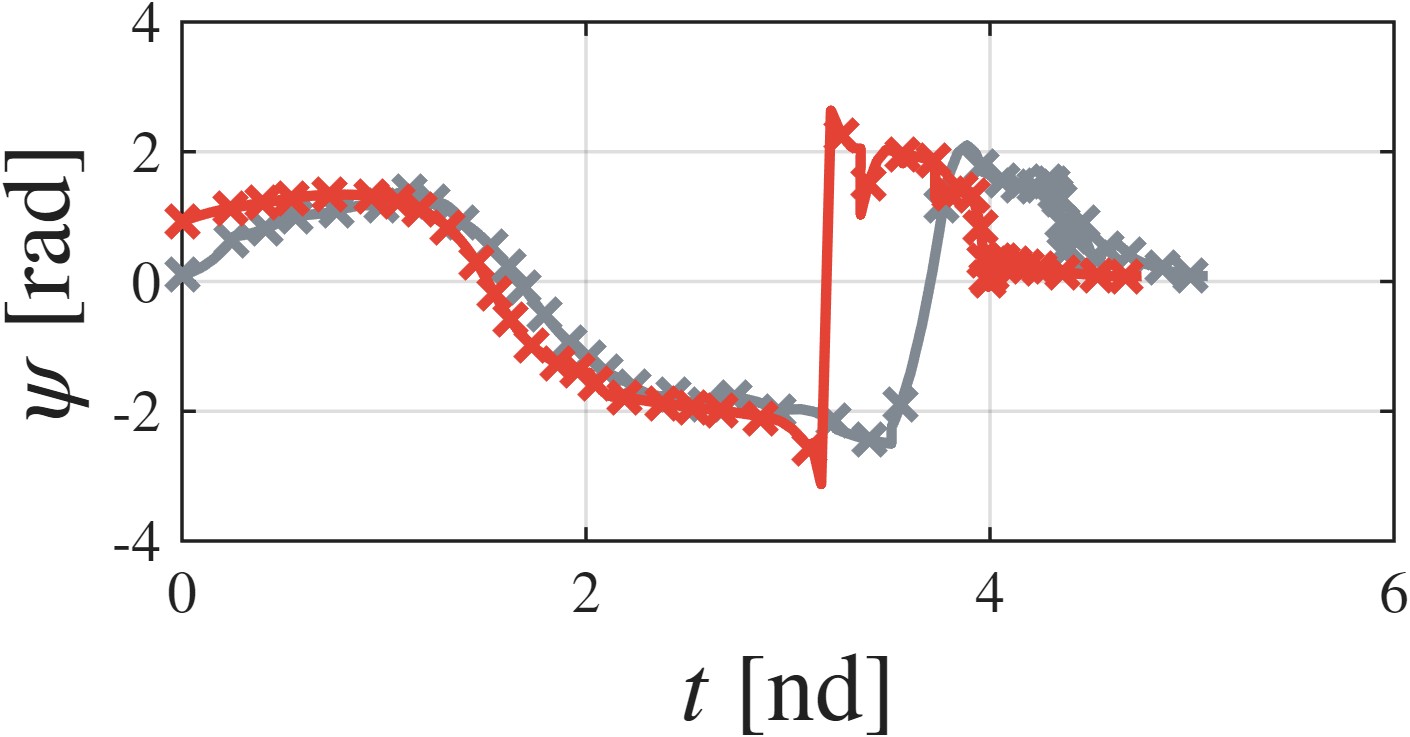}
         \caption{$u_{\max}=0.09$}
     \end{subfigure}\hfill
     \begin{subfigure}{0.24\textwidth}
         \centering
         \includegraphics[width=1\textwidth]{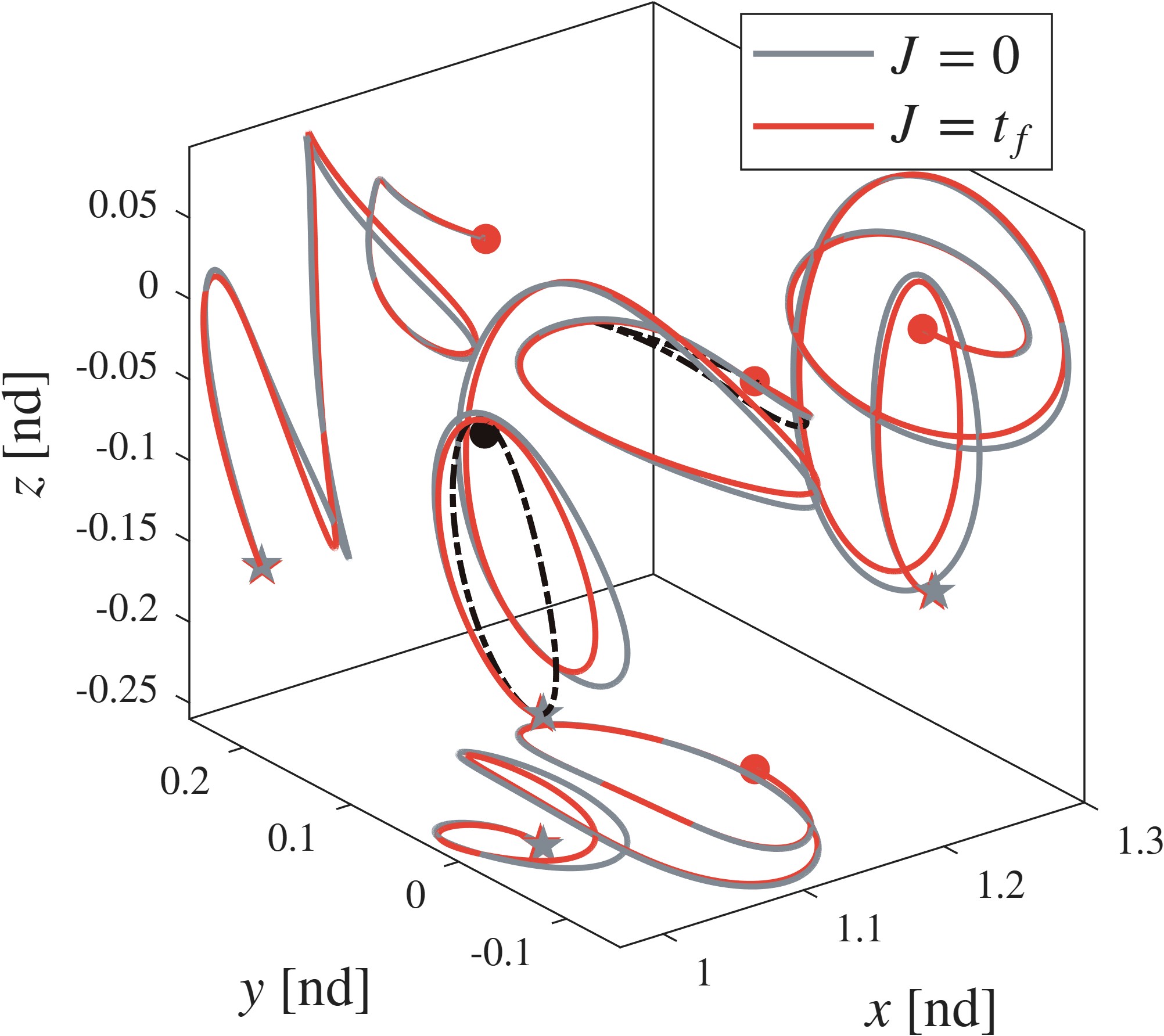}\\
         \includegraphics[width=1\textwidth]{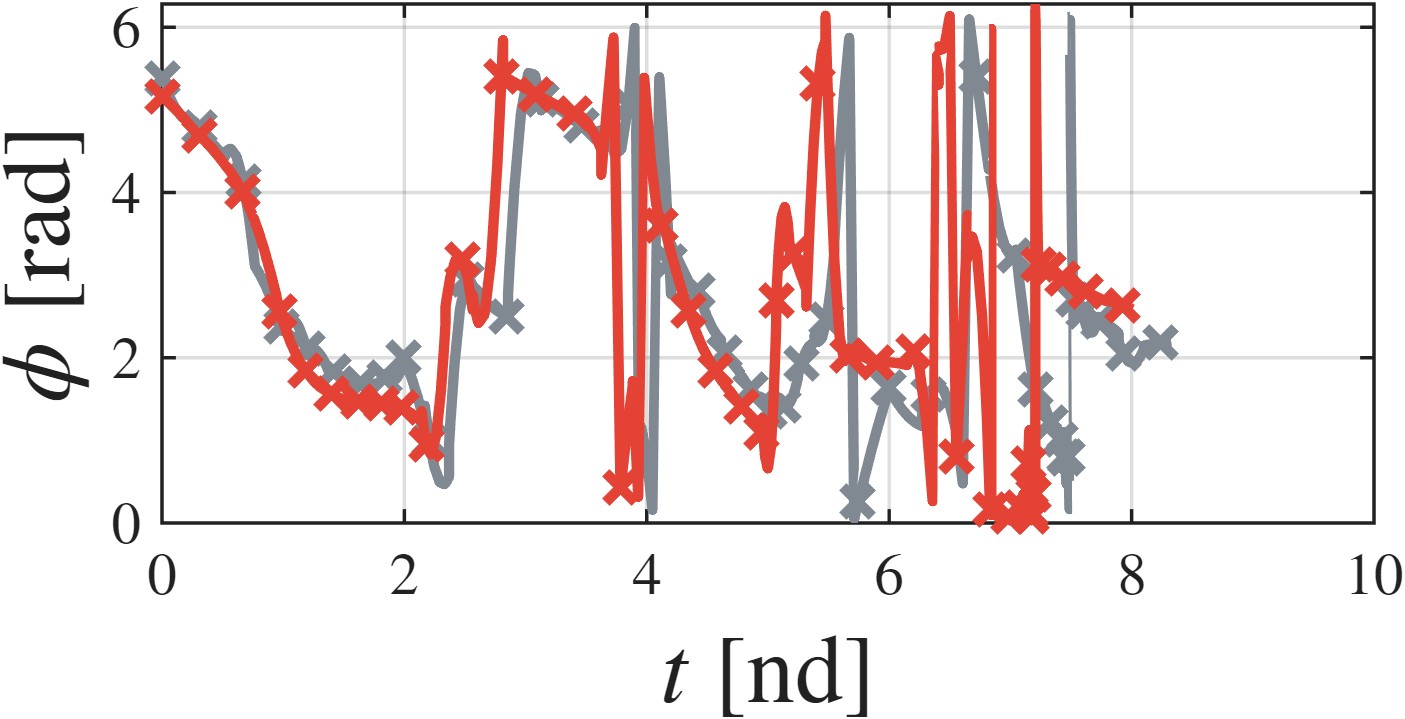}\\
         \includegraphics[width=1\textwidth]{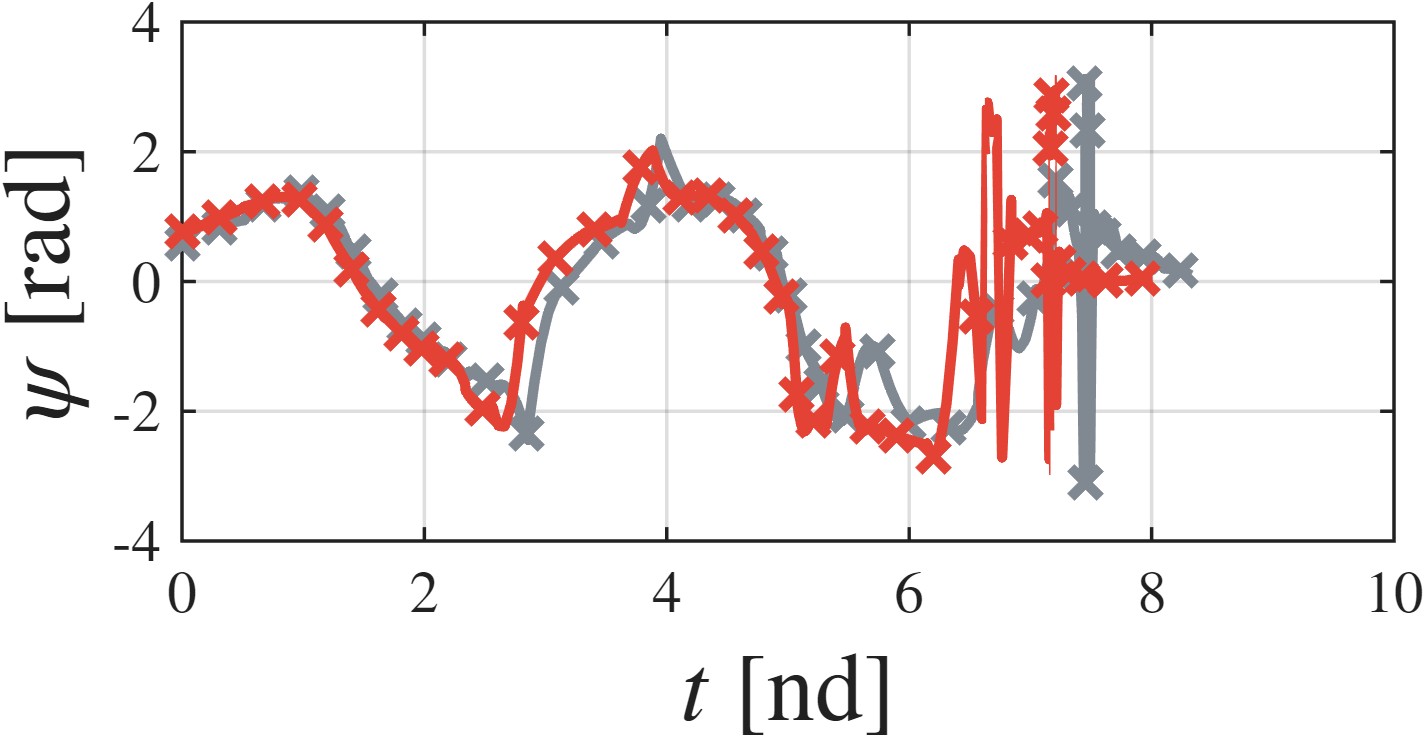}
         \caption{$u_{\max}=0.05$}
     \end{subfigure}\hfill
     \begin{subfigure}{0.24\textwidth}
         \centering
         \includegraphics[width=1\textwidth]{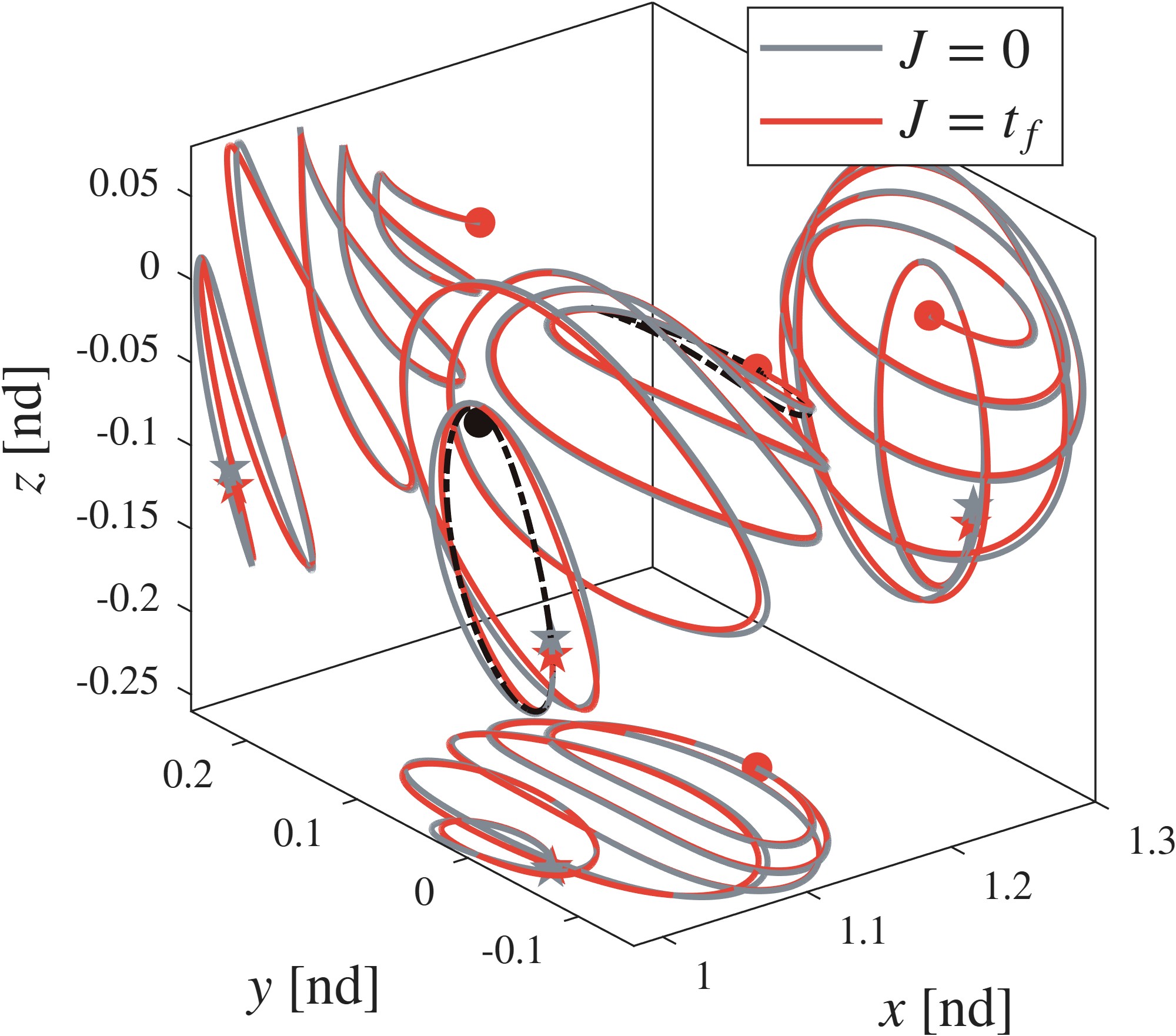}\\
         \includegraphics[width=1\textwidth]{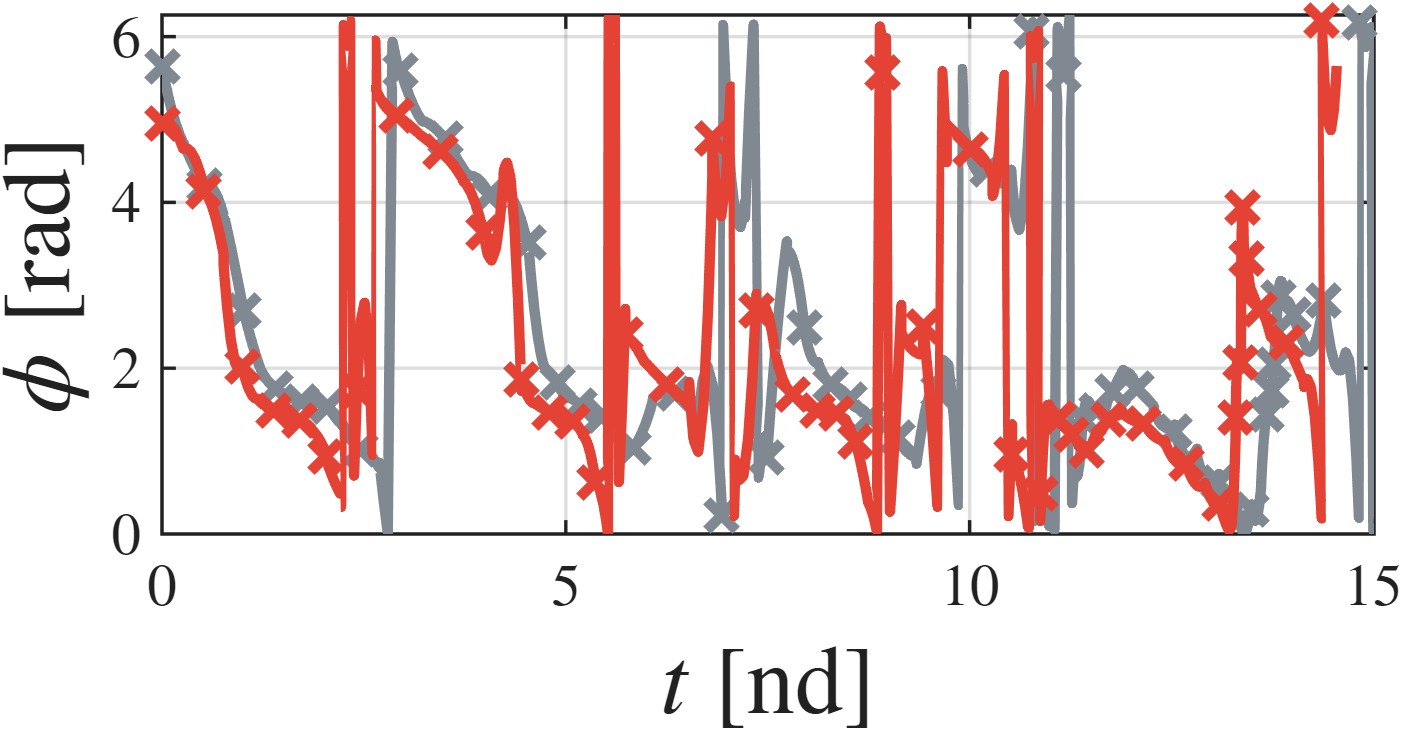}\\
         \includegraphics[width=1\textwidth]{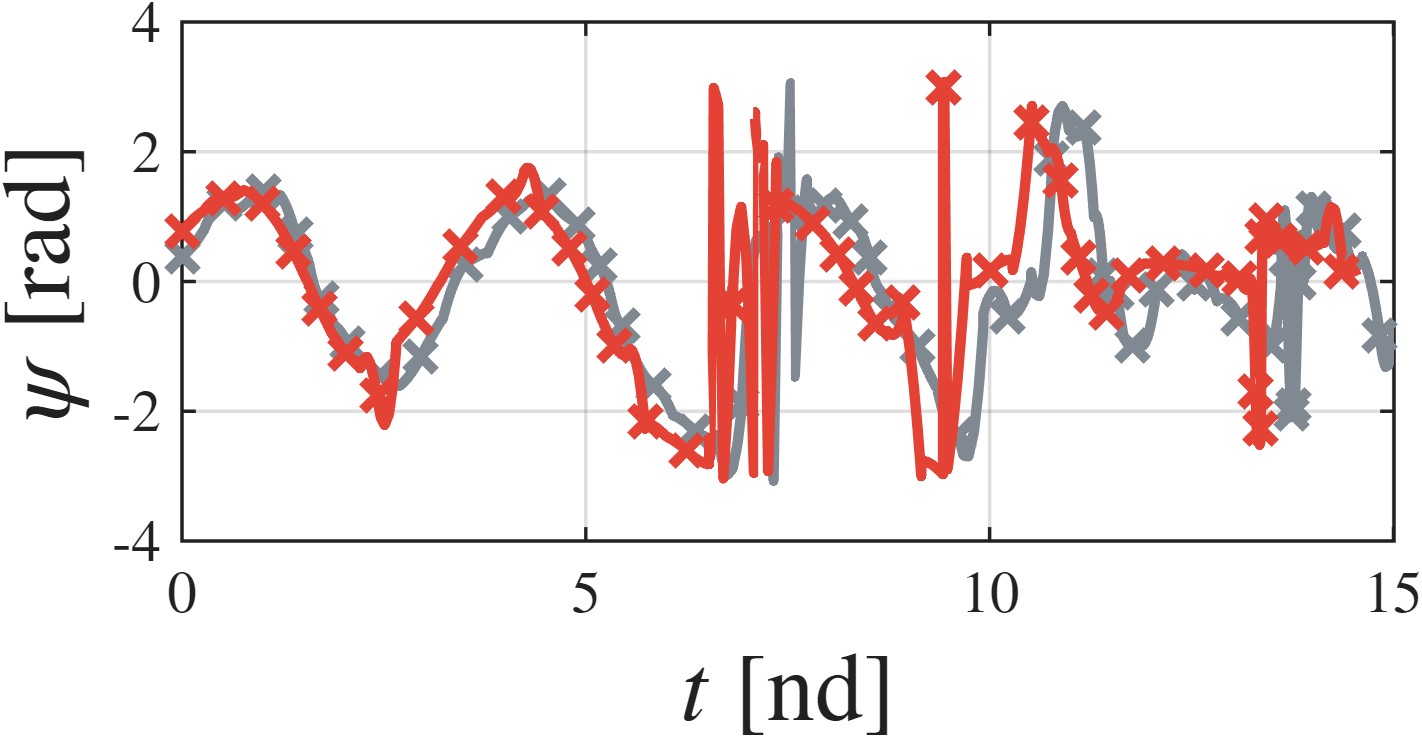}
         \caption{$u_{\max}=0.03$}
     \end{subfigure}
     \caption{Phase-free minimum time solutions for $L_2$ halo transfer.}
     \label{fig:halo_mintime}
\end{figure}

\begin{figure}[htbp!]
     \centering
     \begin{subfigure}{0.24\textwidth}
         \centering
         \includegraphics[width=1\textwidth]{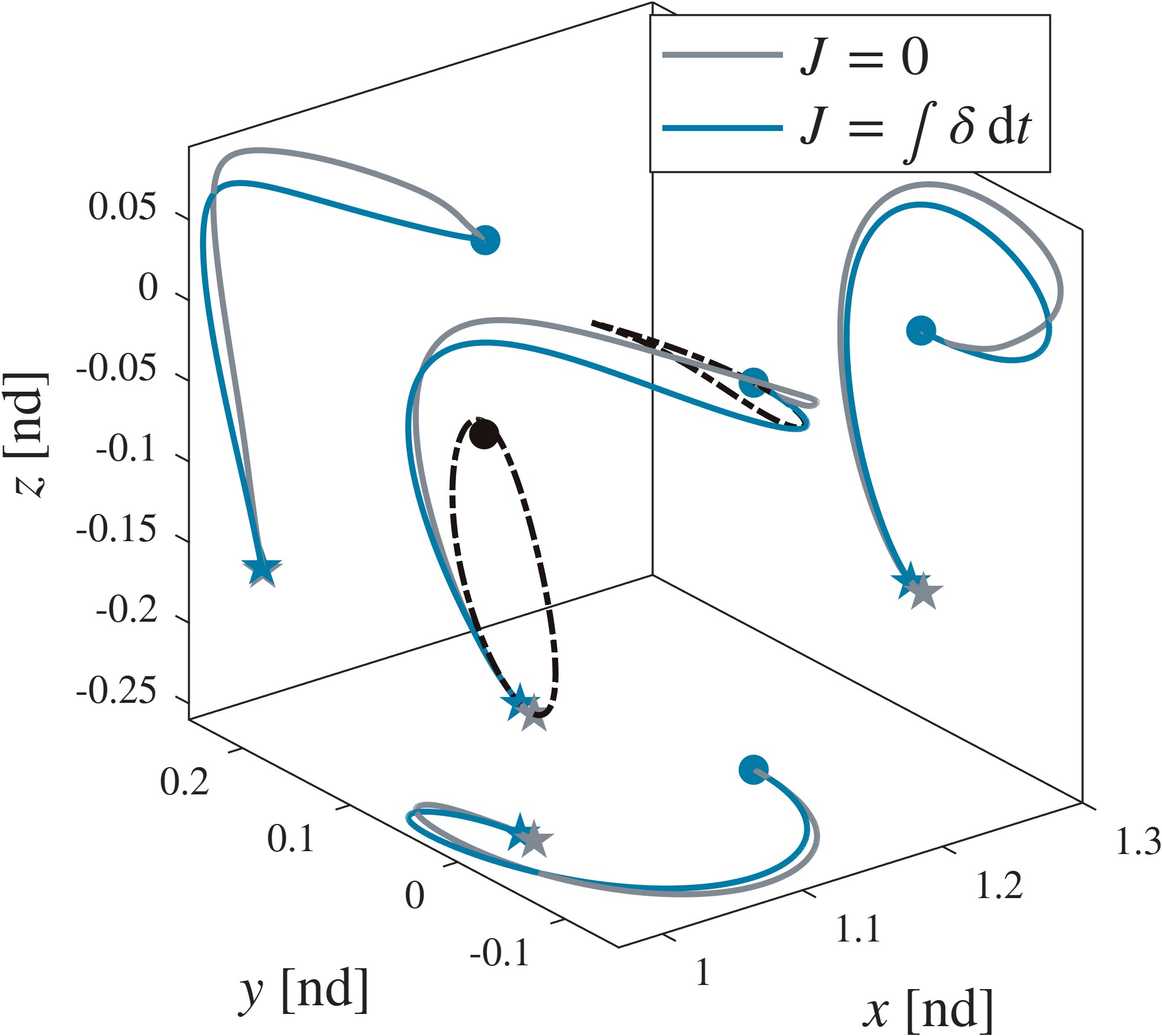}\\
         \includegraphics[width=1\textwidth]{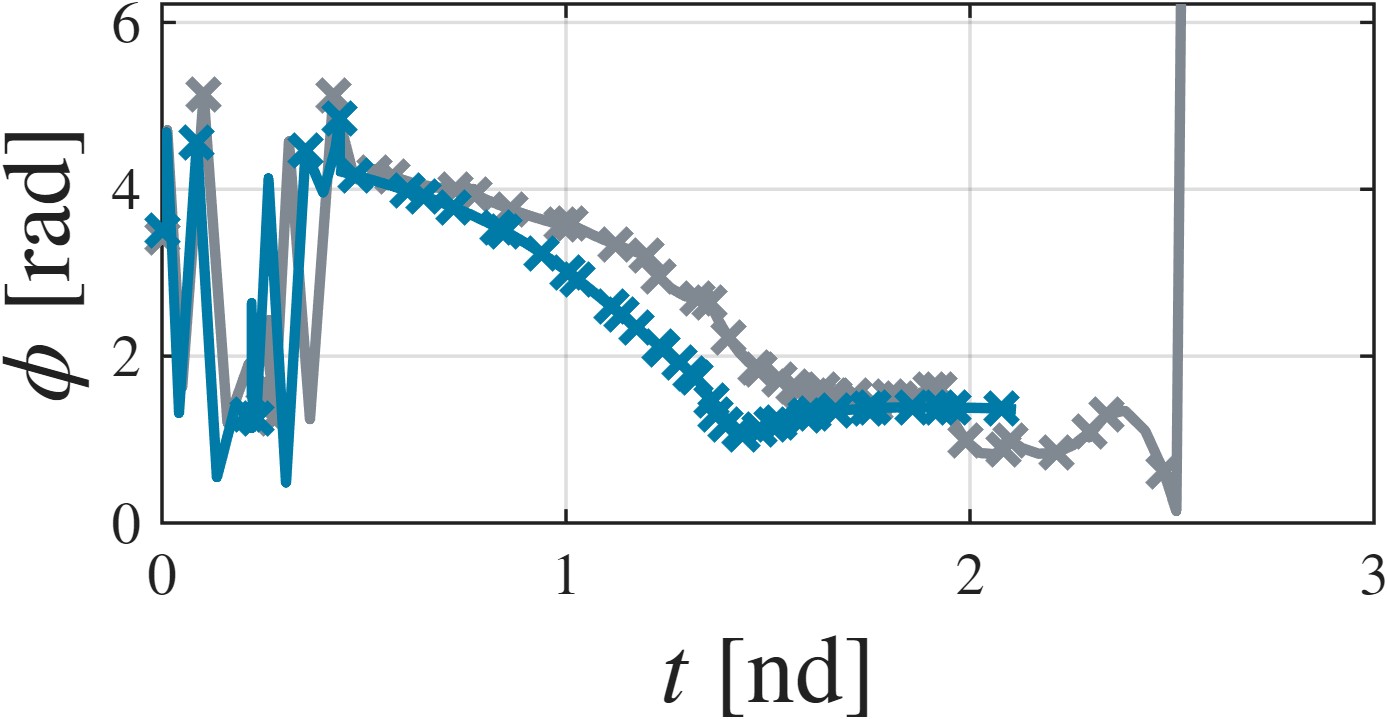}\\
         \includegraphics[width=1\textwidth]{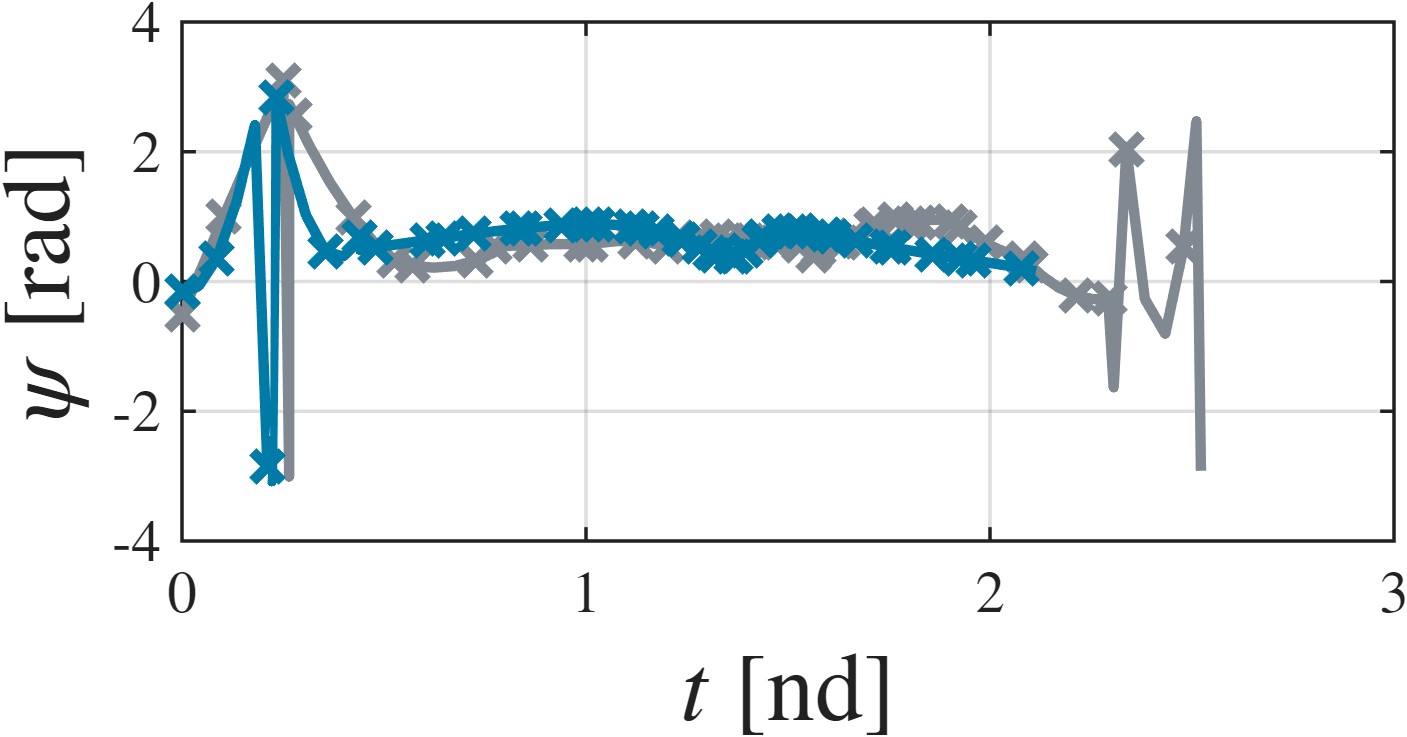}\\
         \includegraphics[width=1\textwidth]{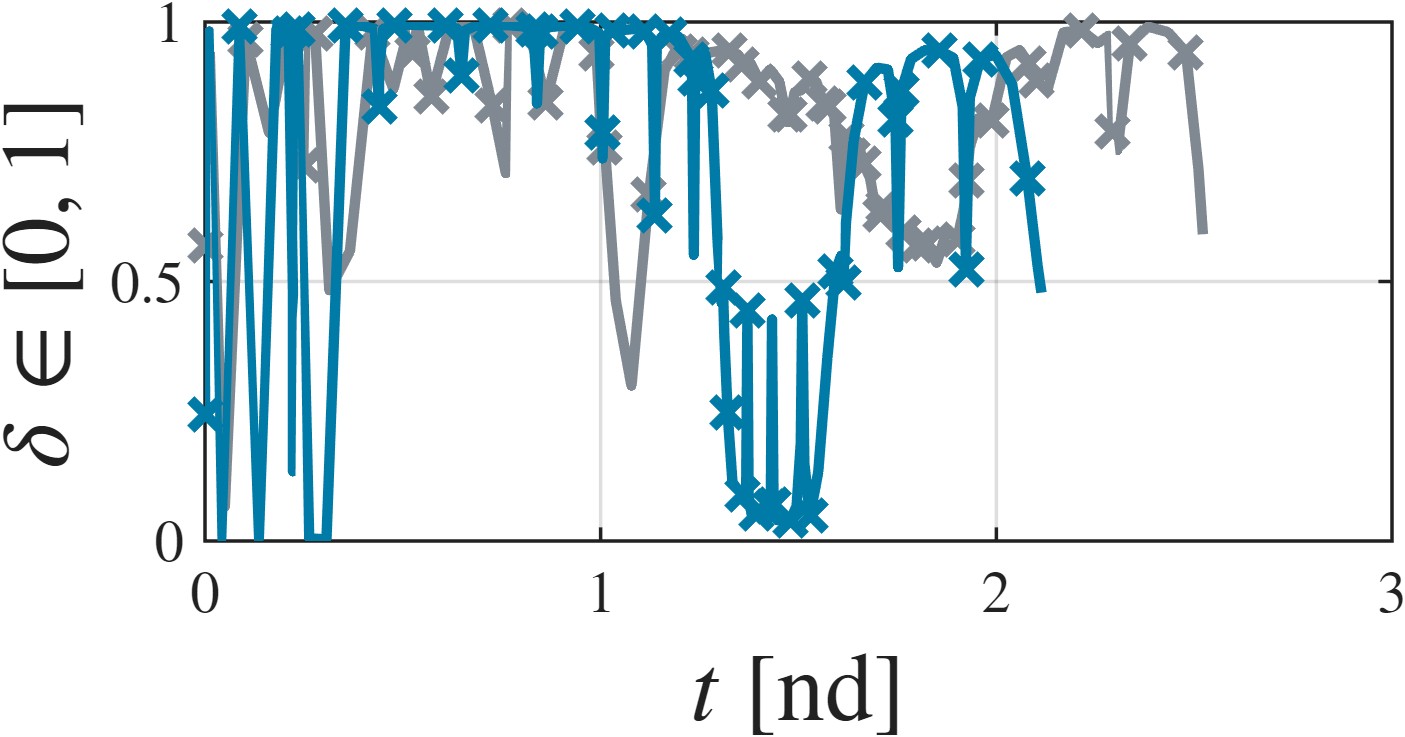}
         \caption{$u_{\max}=0.3$}
     \end{subfigure}\hfill
     \begin{subfigure}{0.24\textwidth}
         \centering
         \includegraphics[width=1\textwidth]{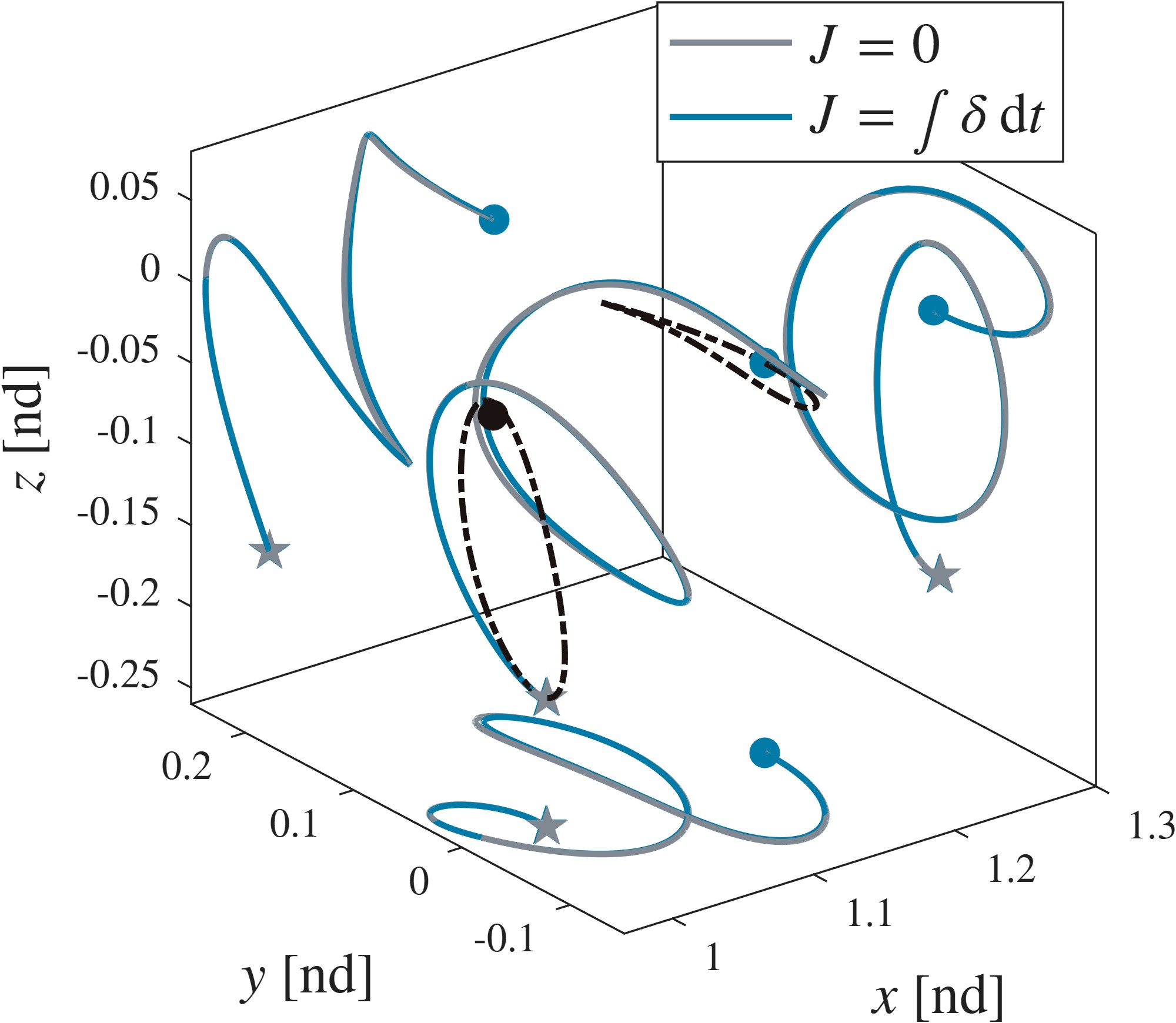}\\
         \includegraphics[width=1\textwidth]{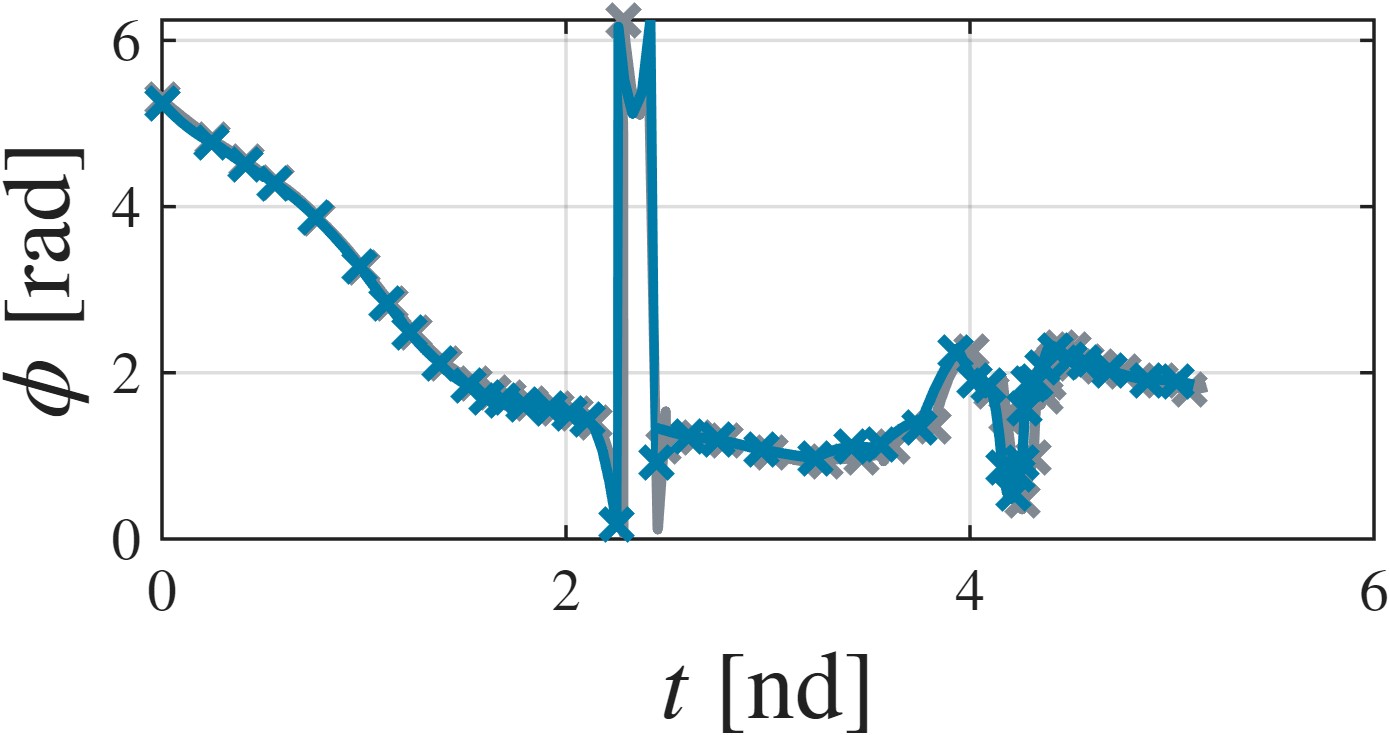}\\
         \includegraphics[width=1\textwidth]{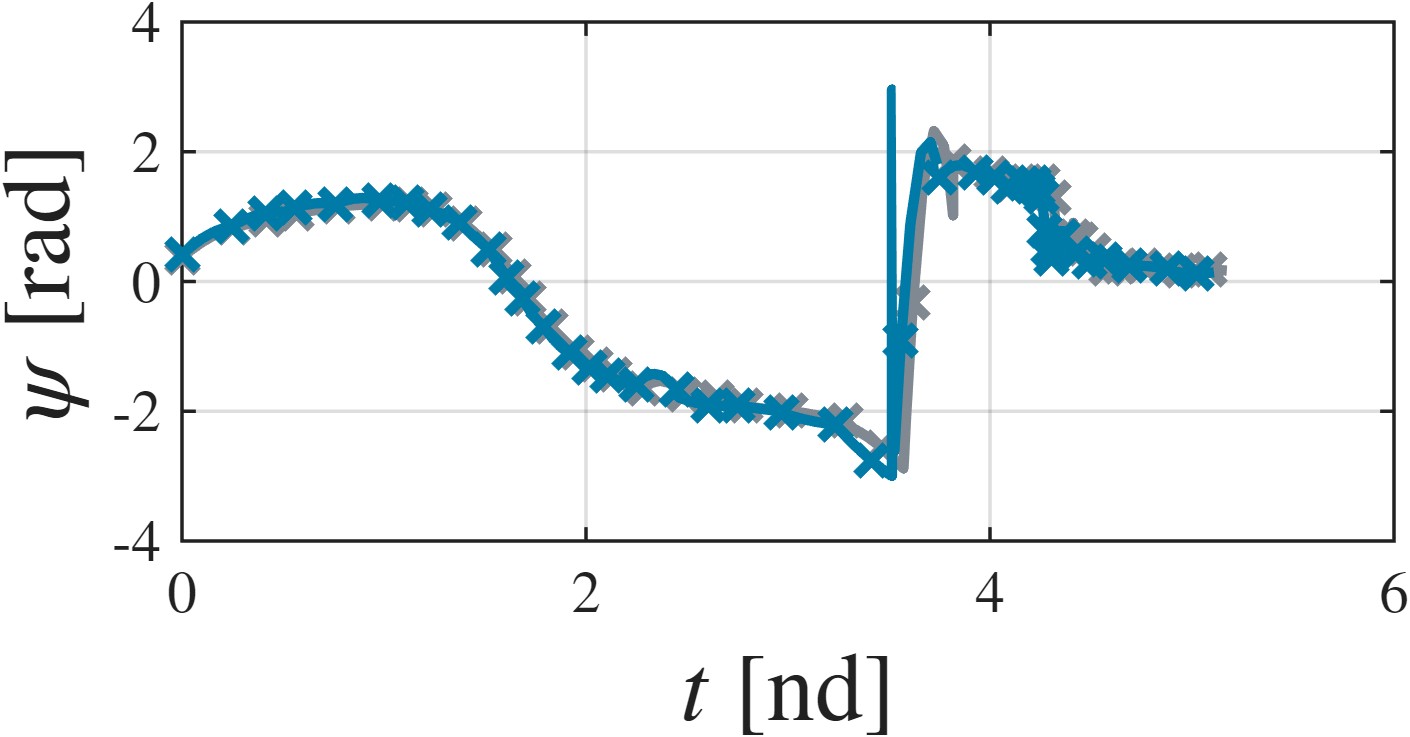}\\
         \includegraphics[width=1\textwidth]{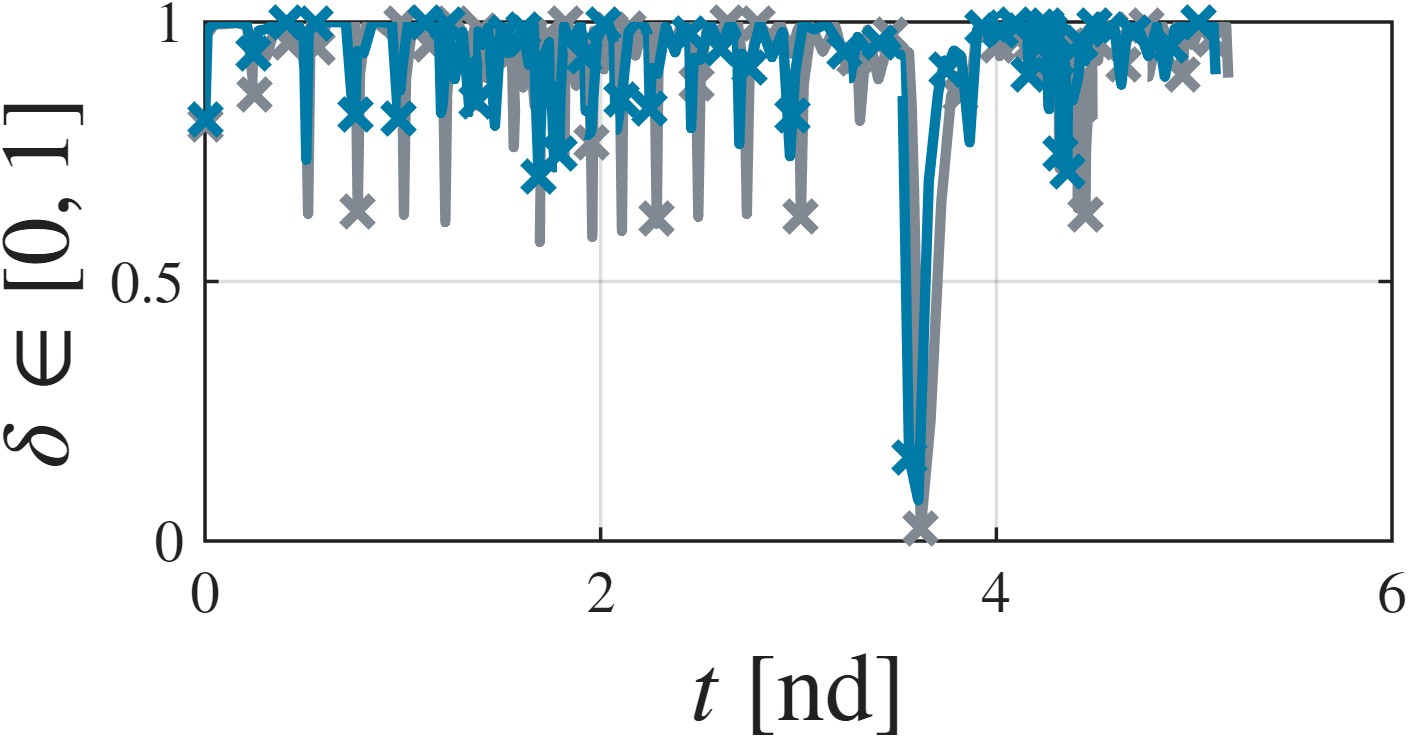}
         \caption{$u_{\max}=0.09$}
     \end{subfigure}\hfill
     \begin{subfigure}{0.24\textwidth}
         \centering
         \includegraphics[width=1\textwidth]{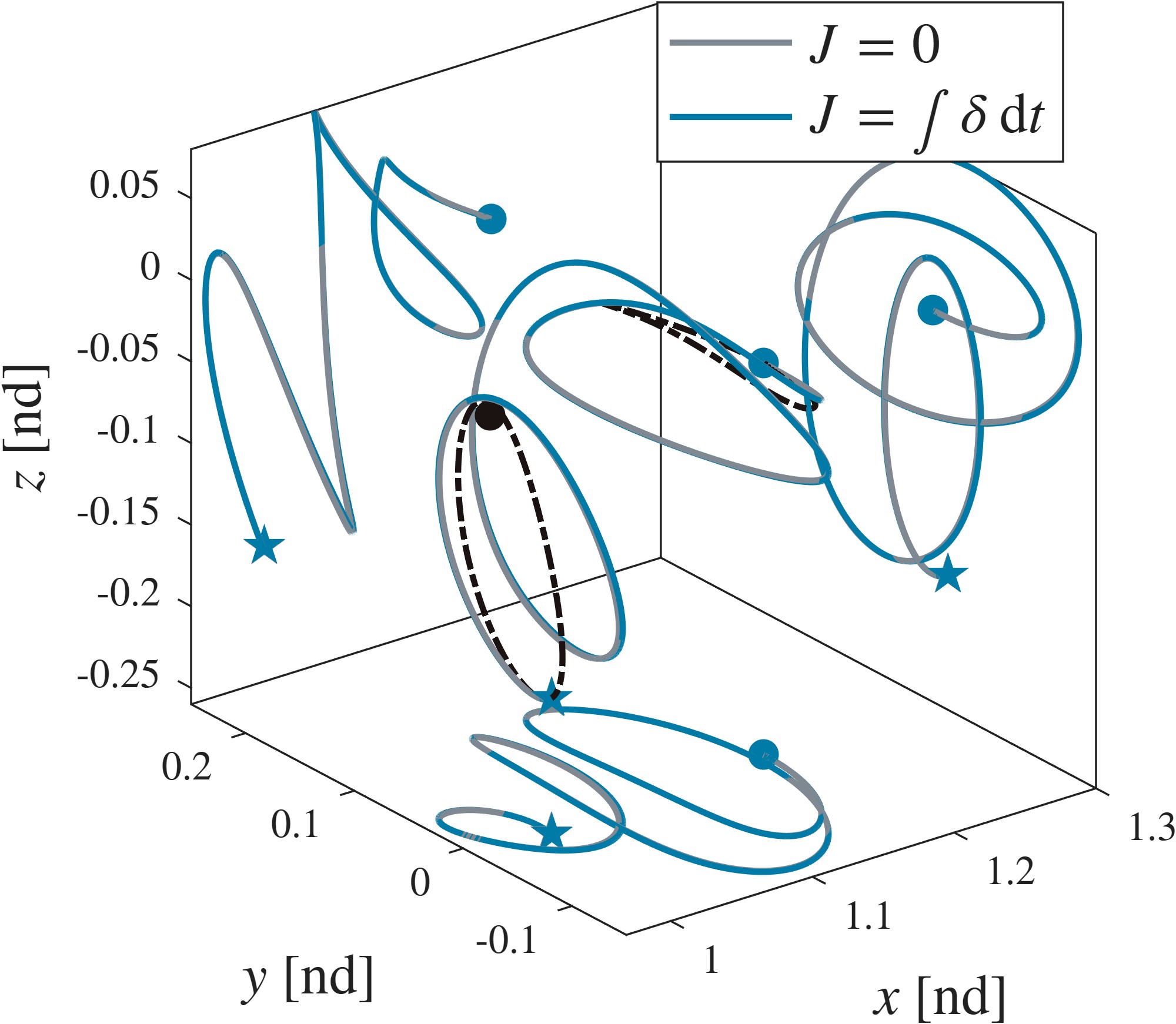}\\
         \includegraphics[width=1\textwidth]{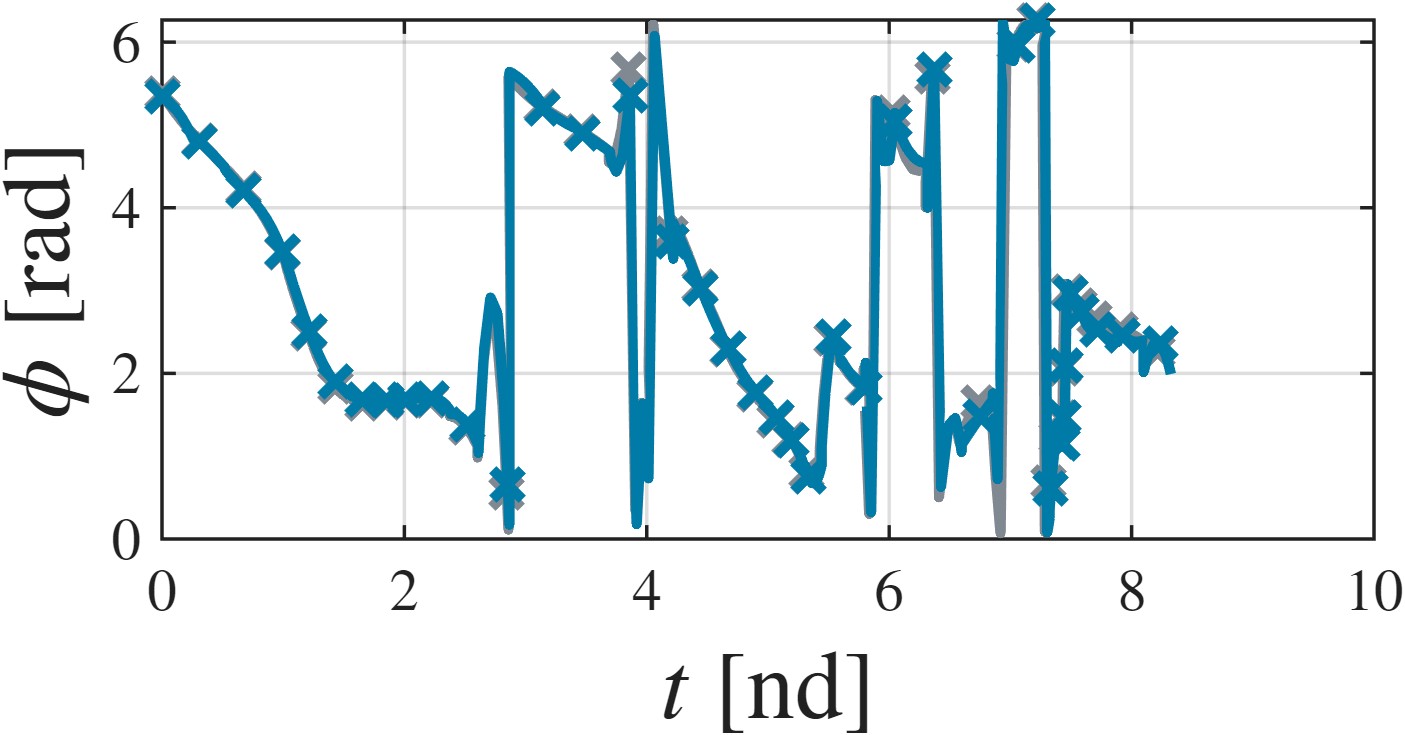}\\
         \includegraphics[width=1\textwidth]{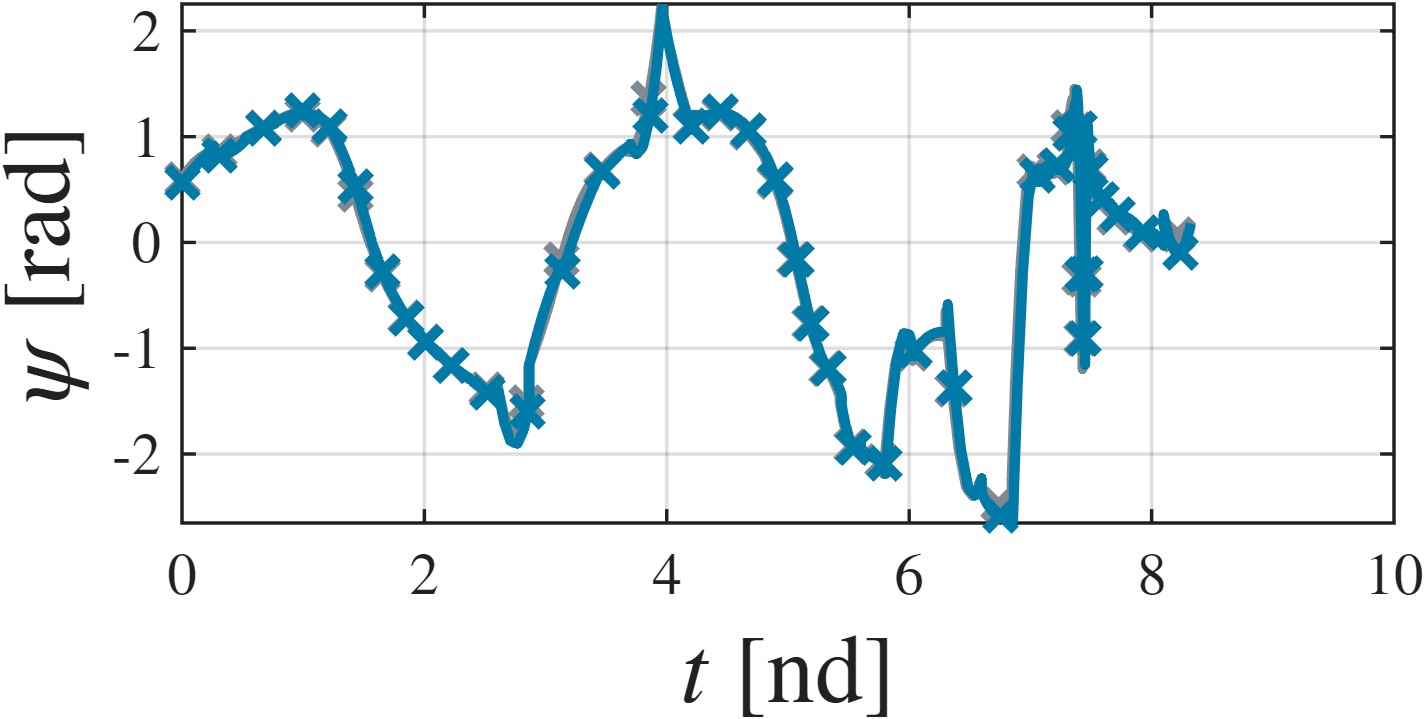}\\
         \includegraphics[width=1\textwidth]{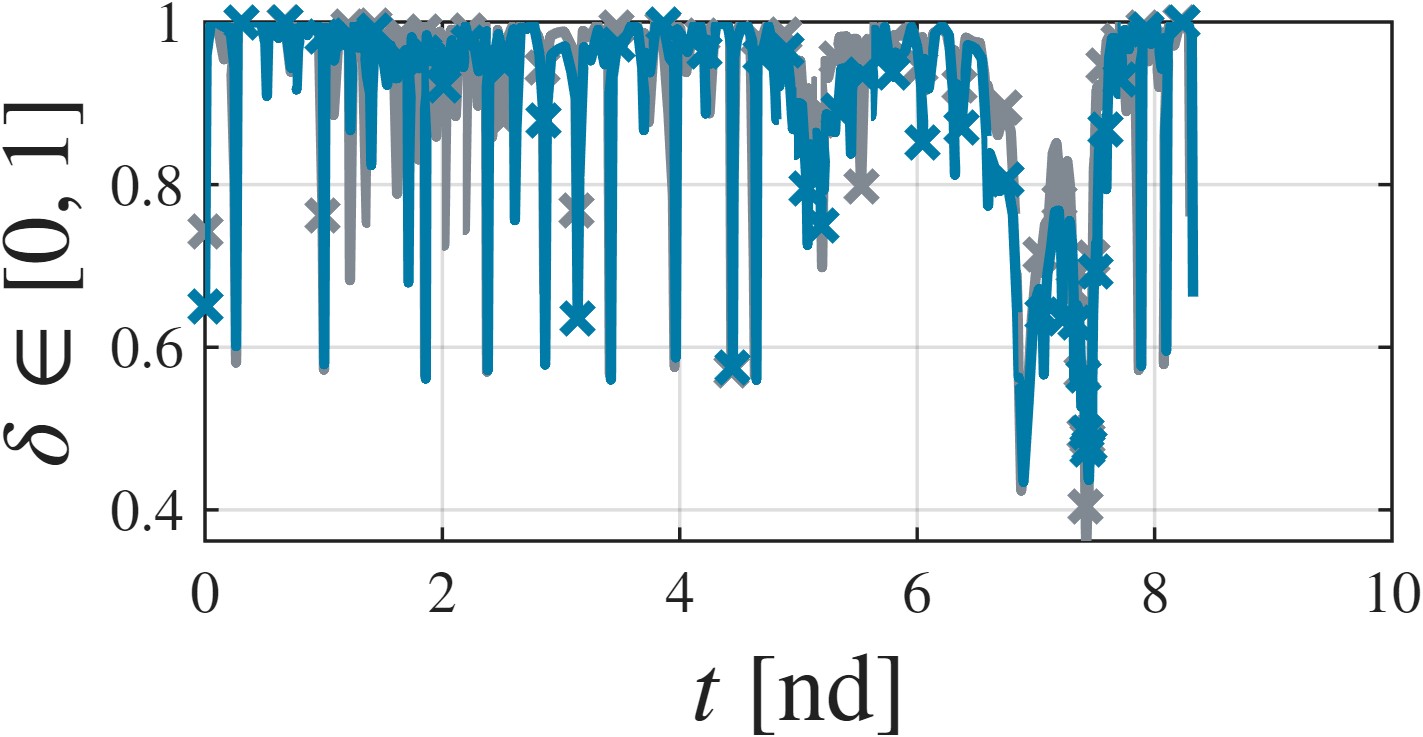}
         \caption{$u_{\max}=0.05$}
     \end{subfigure}\hfill
     \begin{subfigure}{0.24\textwidth}
         \centering
         \includegraphics[width=1\textwidth]{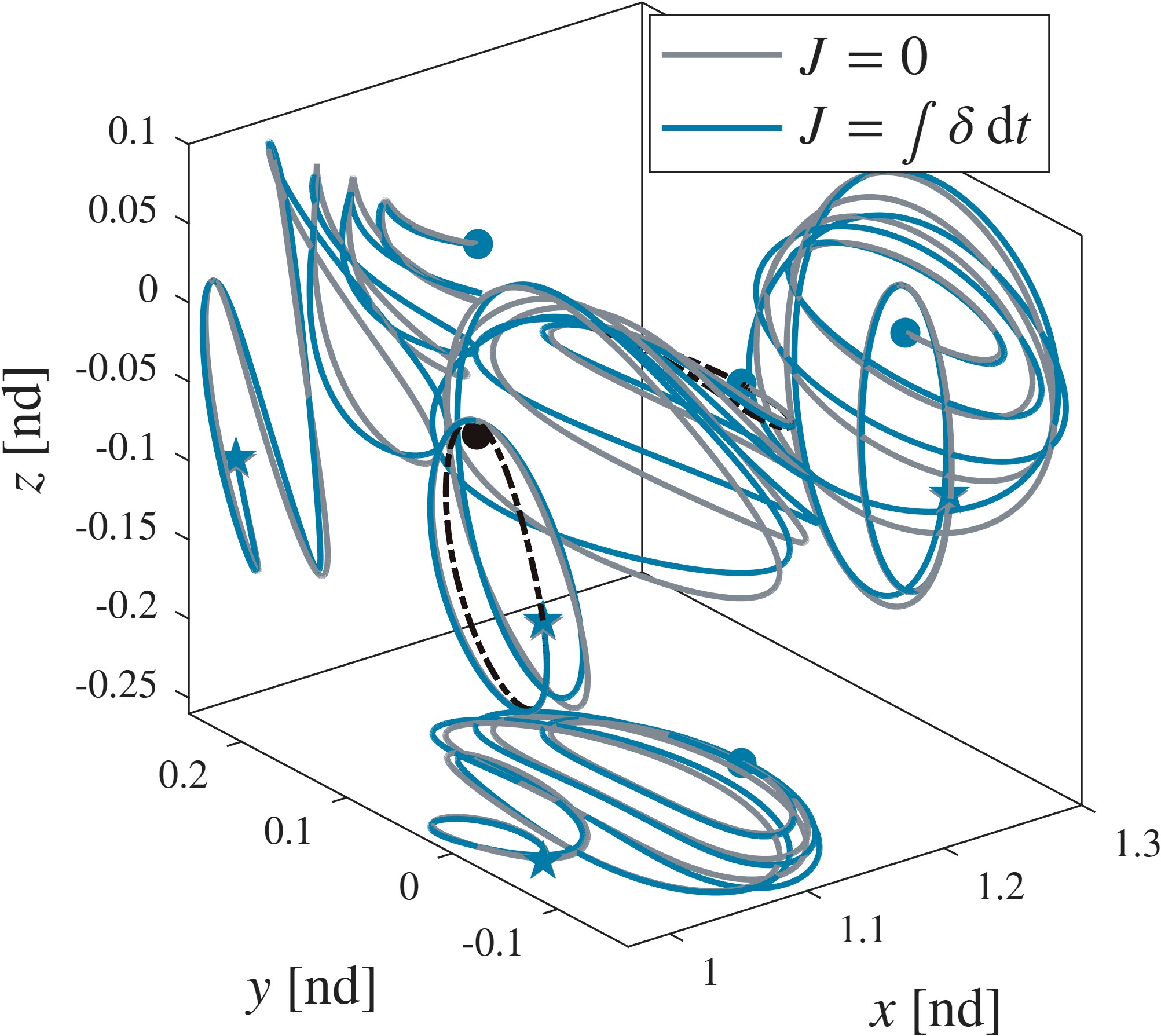}\\
         \includegraphics[width=1\textwidth]{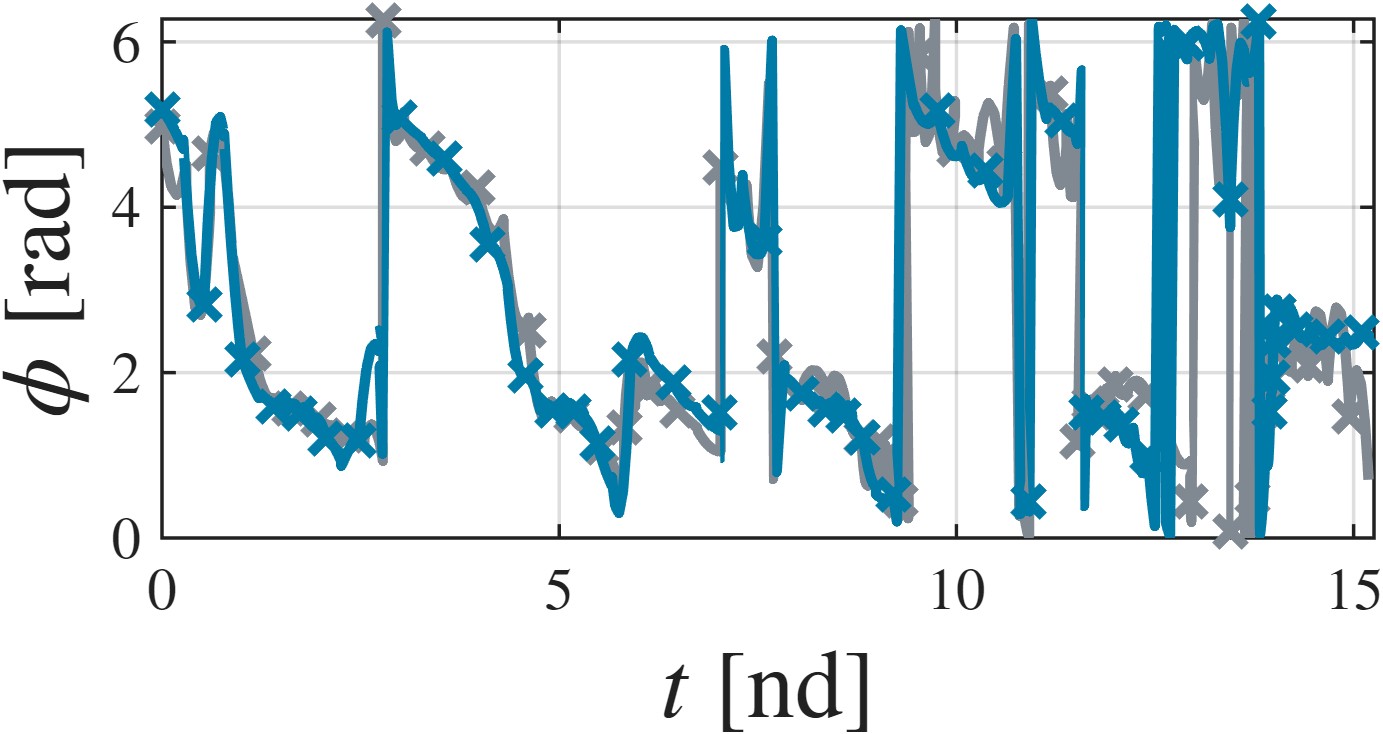}\\
         \includegraphics[width=1\textwidth]{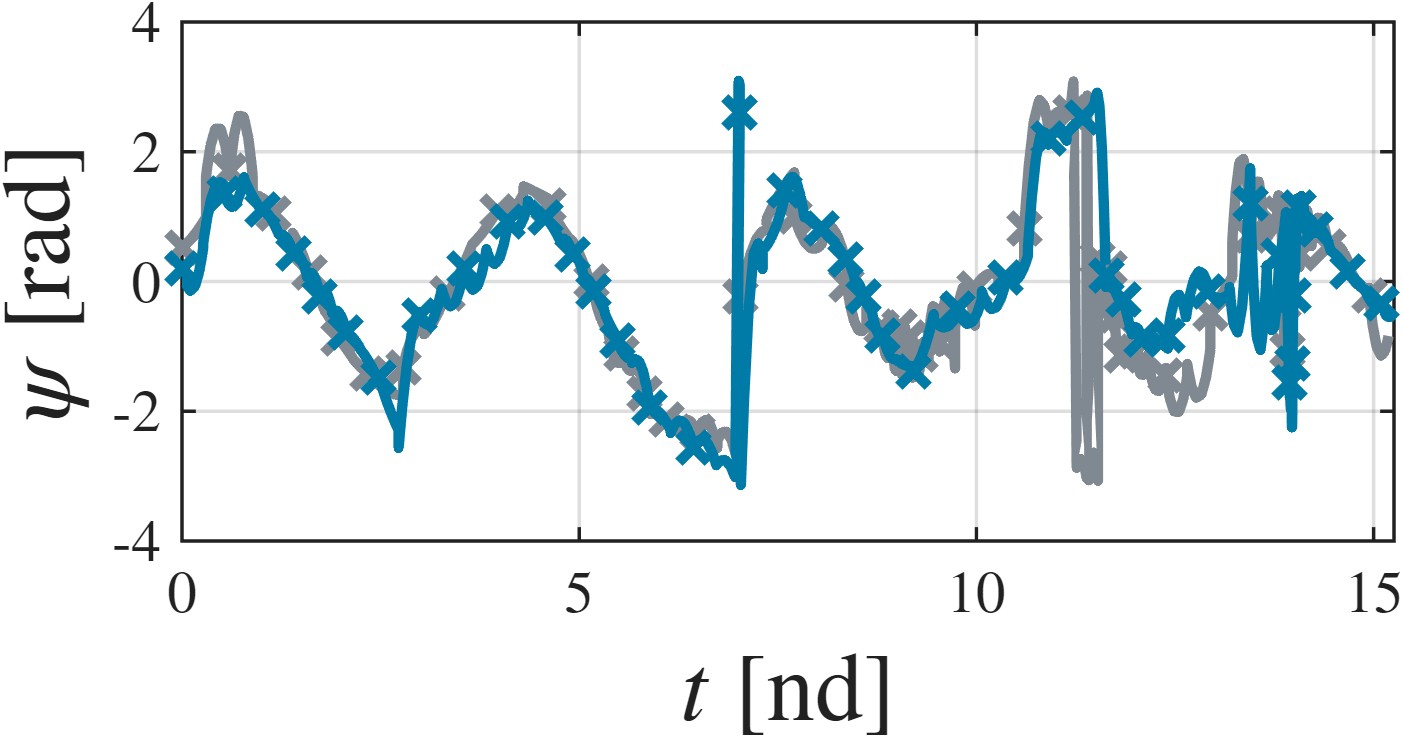}\\
         \includegraphics[width=1\textwidth]{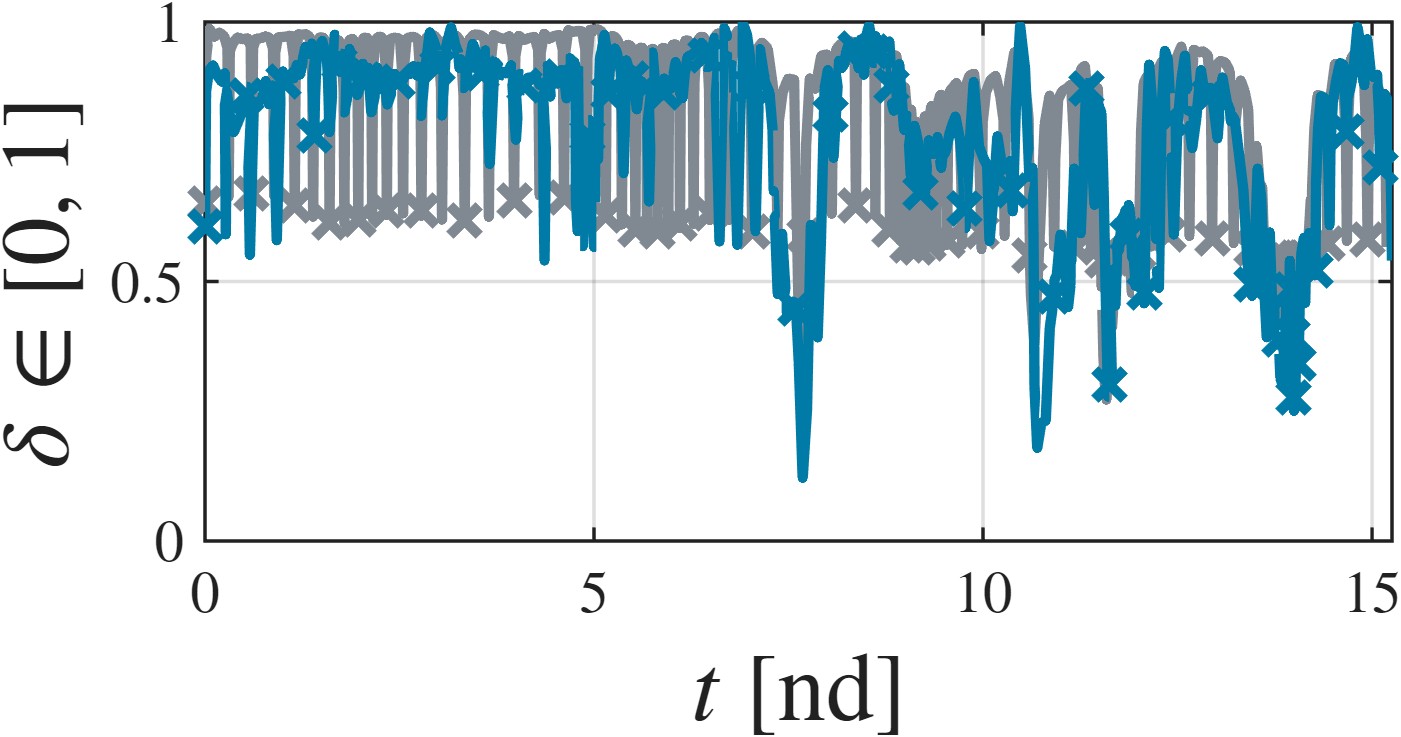}
         \caption{$u_{\max}=0.03$}
     \end{subfigure}
     \caption{Phase-free minimum fuel solutions for $L_2$ halo transfer.}
     \label{fig:halo_minfuel}
\end{figure}

Halo-to-halo transfers are constructed from $\bar p_0 = 0.05$ to $\bar p_f = 0.95$, spanning $90\%$ of the fitted family, at the same four maximum control magnitudes used for the DRO family. The departure orbit is a near-planar halo with a Jacobi constant of $C_J = 3.1504$, an out-of-plane amplitude of $0.0201$ nd, a perilune radius of $50{,}667$ km and a period of $3.411$ TU, or $14.83$ days. The arrival orbit is an NRHO with $C_J = 3.0469$, an out-of-plane amplitude of $0.1818$ nd, a perilune radius of $3{,}190$ km and a period of $1.506$ TU, or $6.55$ days. The transfer therefore traverses a Jacobi constant difference of $\Delta C_J = 0.1035$ while reducing perilune radius by a factor of $16$. Family space initialization trajectories are shown in Fig. \ref{fig:halo_famspace}, minimum time solutions in Fig. \ref{fig:halo_mintime}, minimum fuel solutions in Fig. \ref{fig:halo_minfuel}, and numerical results in Table \ref{tab:halo_mintimefuel}.

The central result of this example is visible when examining predicted vs converged minimum times in Table \ref{tab:halo_mintimefuel}. For the DRO family, the family space prediction of the time of flight was uniformly conservative, with the converged minimum time falling between $32\%$ and $65\%$ of the predicted value and averaging $44\%$. With time regularization and an arclength parameterization, the halo family space prediction lies within $18\%$ of the converged minimum time in the worst case and within $1.3\%$ at $u_{\max} = 0.09$, averaging a ratio of $1.096$ over the four cases. The predicted revolution count is exact in all four cases. The reduced manifold is no longer merely indicating the right solution basin; it is returning a quantitatively accurate estimate of the achievable minimum time before any transcription is performed.

The mechanism appears to be the velocity slip penalty $\alpha$. Across the DRO fit used for Table~1, $\alpha$ ranges over $[0.050,\,0.977]$; the family space traverse is throttled wherever the family function is stiff, and the resulting guess is correspondingly slow. Across the arclength halo fit, $\alpha$ never leaves $[0.757,\,0.903]$, so the guess is close to a genuine
maximum-rate traverse of the family and the predicted time of flight is close to what the phase space can actually achieve. More studies over different orbit families need to be conducted to understand if this success is uniformly applicable.

The family space estimate is no longer conservative: it under-predicts the converged minimum time in three of the four cases. The same minimum time prediction accuracy is directly responsible for the poor minimum fuel performance for this halo example. Because the family space initialization lands so close to the true minimum time, the feasible solution it seeds occupies a region in which nearly all available control is required simply to close the transfer. The minimum fuel problem is time-free and could in principle lengthen $t_f$ and coast, but it instead converges to a local optimum immediately adjacent to the minimum time solution, with duty cycles of $74.7\%$, $93.7\%$, $92.5\%$ and $78.7\%$. The resulting costs are essentially the minimum time costs: $484.4$ against $493.2$ m/s, $440.6$ against $437.4$ m/s, $394.0$ against $409.7$ m/s, and $368.6$ against $446.5$ m/s. At $u_{\max} = 0.09$ the minimum fuel solution is actually more expensive than the minimum time solution despite a longer flight time of $22.21$ against $20.66$ days, which is a strong indicator of a poor local minimum. A straightforward remedy, left to future work, is to seed the minimum fuel problem from a deliberately inflated time of flight rather than from the minimum time estimate, which restores the coast structure the throttle needs. Another avenue of future work could make use of a damping coefficient to tune the velocity slip penalty function $\alpha$ for minimum fuel transfers.

\begin{table}[htbp]
    \centering
    \caption{Phase-free minimum time and minimum fuel solution data for the $L_2$ halo transfer.}
    \begin{tabular}{|l|c|c c c c|}
        \hline
        Solution & $u_{\max}$ &
        \makecell{0.300 [nd] \\ 0.817 [mm/s$^2$]} &
        \makecell{0.090 [nd] \\ 0.245 [mm/s$^2$]} &
        \makecell{0.050 [nd] \\ 0.136 [mm/s$^2$]} &
        \makecell{0.030 [nd] \\ 0.081 [mm/s$^2$]} \\
        \hline\hline

        \multirow{2}{*}{Family Space} & $N_\text{rev}$ & 1 & 2 & 3 & 5\\
                                      & $t_f$ [days]   & 5.914& 20.92&31.74 & 56.57\\
        \hline
        \multirow{5}{*}{Minimum Time} & $N_\text{poly}$ & 7 & 7& 7& 7\\
                                      & $N_\text{seg}$ & 15 &32 &48 & 81\\
                                      & $N_\text{rev}$ &1 &2 & 3& 5\\
                                      & $t_f$ [days]   & 6.987& 20.66& 34.82& 63.25\\
                                      & $\Delta V$ [m/s] & 493.2& 437.4& 409.7 & 446.5\\
        \hline
        \multirow{5}{*}{Minimum Fuel} & $N_\text{poly}$ &7 &7 &7 &7 \\
                                      & $N_\text{seg}$ & 15& 32 & 48& 81\\
                                      & $N_\text{rev}$ &1 & 2& 3& 5\\
                                      & $t_f$ [days]   &9.193 &22.21 &36.20 & 66.34\\
                                      & $\Delta V$ [m/s] &484.4 &440.6 &394.0 & 368.6\\
        \hline
    \end{tabular}
    \label{tab:halo_mintimefuel}
\end{table}
Lastly, the control angle histories in Figs. \ref{fig:halo_mintime} and \ref{fig:halo_minfuel} remain smooth over the majority of each transfer, with the expected concentration of activity at the close lunar passes that dominate the NRHO end.

\section{Conclusion}
This paper presented a robust, unified framework for the design and optimization of many-revolution, low-thrust transfers within periodic orbit families in the Circular Restricted Three-Body Problem (CR3BP). By leveraging a continuous, multi-segment Chebyshev-Fourier analytical representation of the periodic family, the high-dimensional, highly sensitive phase-space dynamics were projected into a smooth, lower-dimensional manifold parameter space. The resulting family-space targeting scheme successfully resolved the complex topological structure and high winding numbers of multi-revolution transfers, providing highly accurate, physically consistent initial guesses for both time- and fuel-optimal trajectories. Time regularization was used to handle orbit families with close passes to primaries.

The subsequent transcription and refinement of these initial guesses using the second-order direct orthogonal polynomial integral collocation method demonstrated strong numerical performance. By parameterizing only the highest-order acceleration derivatives and analytically reconstructing position and velocity states, the collocation framework dramatically reduced the dimension of the nonlinear programming decision space. This decoupling of the internal collocation grid resolution from the state variables yielded highly sparse Jacobian matrices, resulting in superior convergence stability and computational speed. The optimization sweeps for both minimum-time and minimum-fuel cislunar transfers were exemplified using the Moon-centered Distant Retrograde Orbit (DRO) family showcased the method's ability to smoothly transition between widely separated members of a family over multiple revolutions with differing costs. The framework was then extended to families whose dynamics cannot be represented on a uniform time grid. A multi-body time regularization scheme was introduced to handle this case, and its downstream effects on the method were outlined. Transfer examples were drawn from the $L_2$ halo family, where an arclength family parameter coupled with the time regularized family function produced initial guesses for minimum time and revolution count that were highly accurate when compared to the converged phase space solution. Because the halo family space minimum time estimate was no longer conservative (and instead highly accurate), using it to seed minimum fuel problems yielded costs almost indistinguishable from the minimum time solutions, as the solutions fought to simply complete the transfer with the control available. Broadly speaking, the conclusions drawn from both the DRO and halo family examples cannot be used to make general claims about the newly introduced family space initialization method, although the regularized, arclength-parametrized fits show excellent promise. Application to additional periodic orbit families, particularly those with strong hyperbolicity, is necessary to further assess the method's predictive success.
\bibliographystyle{AAS_publication} 
\bibliography{references}
\appendix
\section*{Appendix: Family Function Sensitivities}
{\normalsize
\begin{equation}
    \boldsymbol{\Gamma}({p},\theta) = \boldsymbol{C}^\text{T}\cdot\left[\boldsymbol{\varphi}(\tau({p}))\otimes\text{e}^{\sqrt{-1}\cdot\boldsymbol{k}\cdot\theta}\right]^\text{T}
\end{equation}
\begin{equation}
    \frac{\partial\boldsymbol{\Gamma}}{\partial \theta} = \boldsymbol{C}^\text{T}\cdot\left[\boldsymbol{\varphi}(\tau({p}))\otimes\sqrt{-1}\cdot\text{e}^{\sqrt{-1}\cdot\boldsymbol{k}\cdot\theta}\cdot\text{diag}(\boldsymbol{k})\right]^\text{T}
\end{equation}
\begin{align}
    \frac{\partial\boldsymbol{\Gamma}}{\partial p} &= \boldsymbol{C}^\text{T}\cdot\left[\frac{\partial\boldsymbol{\varphi}}{\partial p}\otimes\text{e}^{\sqrt{-1}\cdot\boldsymbol{k}\cdot\theta}\right]^\text{T},\quad \frac{\partial\boldsymbol{\varphi}}{\partial p}= \frac{\partial\boldsymbol{\varphi}}{\partial \tau}\frac{\partial\tau}{\partial p} \\
    \frac{\partial\varphi_0(\tau)}{\partial\tau} &= 0,\quad \frac{\partial\varphi_1(\tau)}{\partial\tau} = 1\\
    \frac{\partial\varphi_i(\tau)}{\partial\tau} &= 2\varphi_{i-1}(\tau)+2\tau\frac{\partial\varphi_{i-1}(\tau)}{\partial\tau} - \frac{\partial\varphi_{i-2}(\tau)}{\partial\tau},\quad i\geq2 \\
    \frac{\partial\tau}{\partial p} &= \frac{2}{p_{\max} - p_{\min}}
\end{align}
}
\newpage
\section*{Appendix: Sidereal Resonant Distant Retrograde Orbits in the Earth-Moon System}
\begin{figure}[htbp!]
     \centering
     \begin{subfigure}{0.32\textwidth}
         \centering
         \includegraphics[width=0.99\textwidth]{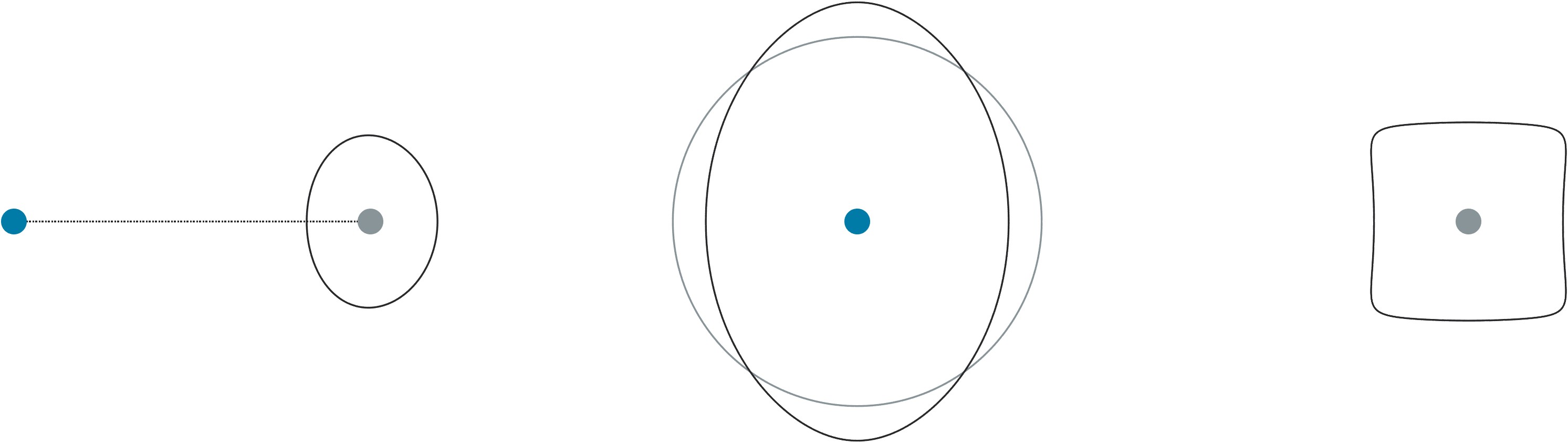}
         \caption{2:1}
     \end{subfigure}\hfill
     \begin{subfigure}{0.32\textwidth}
         \centering
         \includegraphics[width=0.99\textwidth]{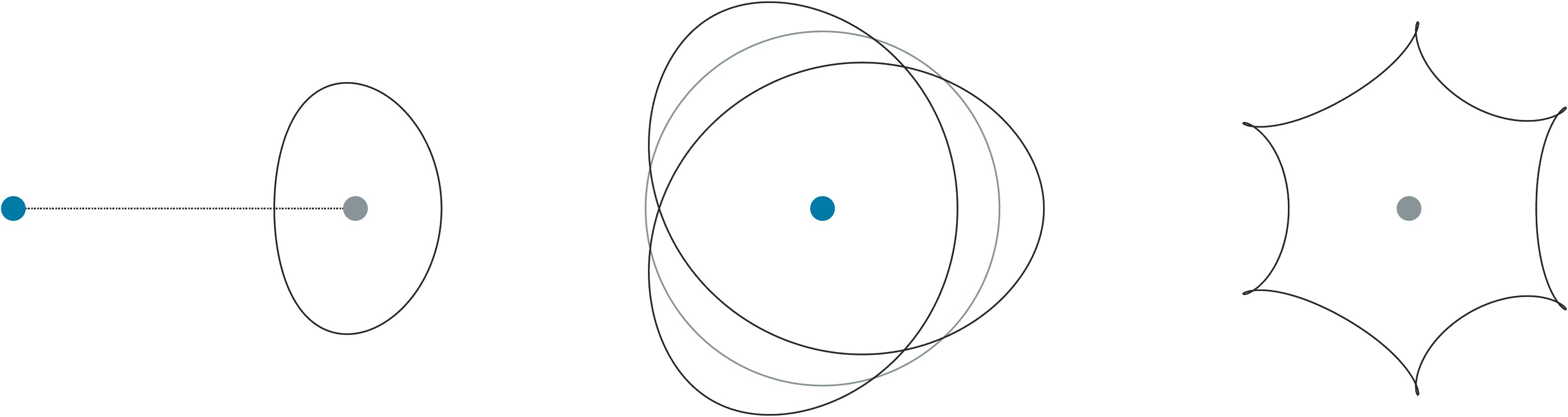}
         \caption{3:2}
     \end{subfigure}\hfill
     \begin{subfigure}{0.32\textwidth}
         \centering
         \includegraphics[width=0.99\textwidth]{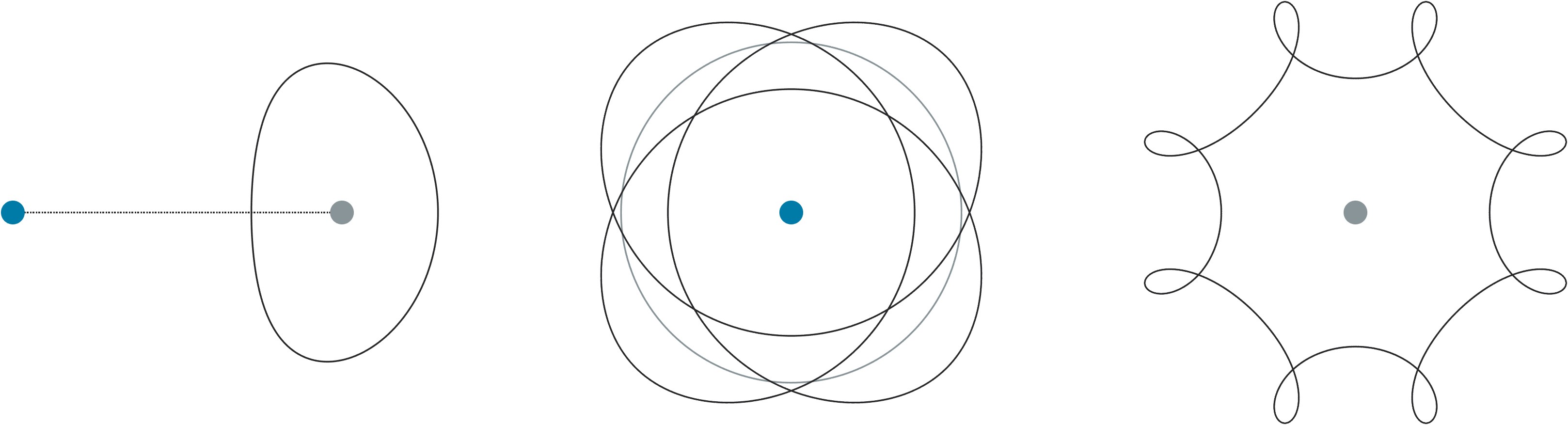}
         \caption{4:3}
     \end{subfigure}\\ % LINE BREAK!
     \begin{subfigure}{0.32\textwidth}
         \centering
         \includegraphics[width=0.99\textwidth]{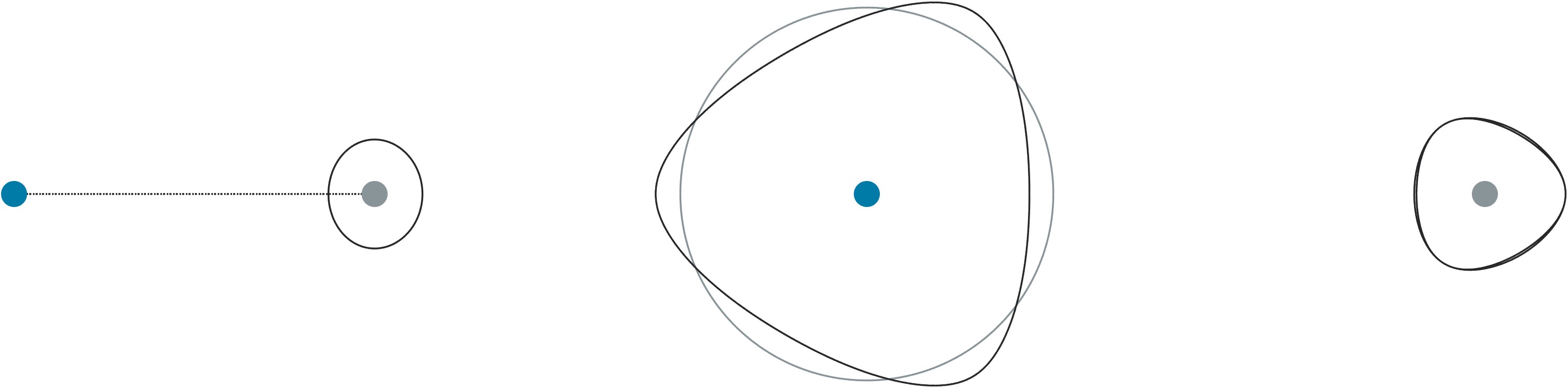}
         \caption{3:1}
     \end{subfigure}\hfill
     \begin{subfigure}{0.32\textwidth}
         \centering
         \includegraphics[width=0.99\textwidth]{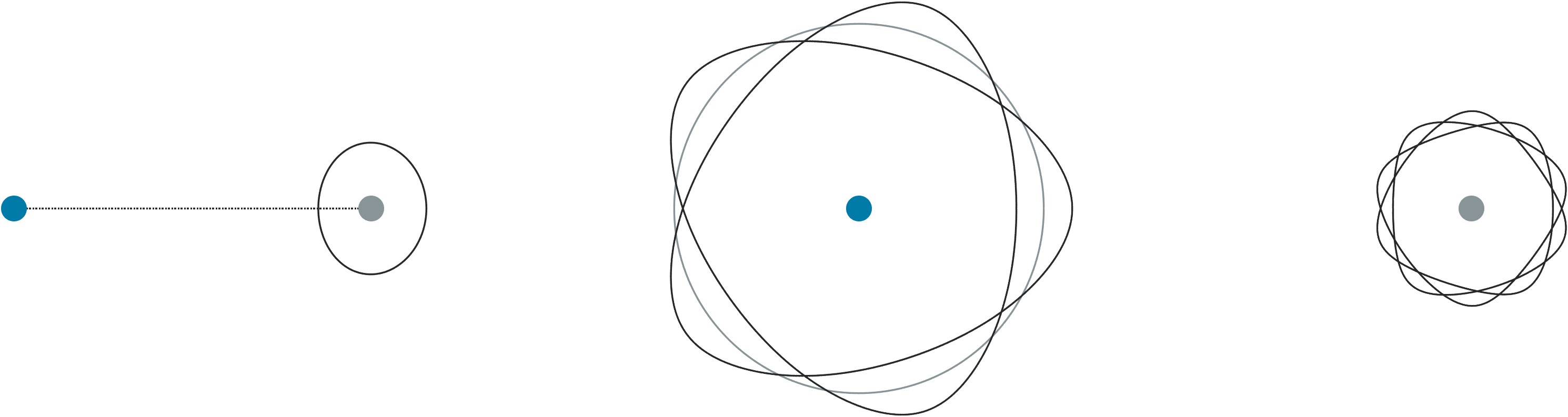}
         \caption{5:2}
     \end{subfigure}\hfill
     \begin{subfigure}{0.32\textwidth}
         \centering
         \includegraphics[width=0.99\textwidth]{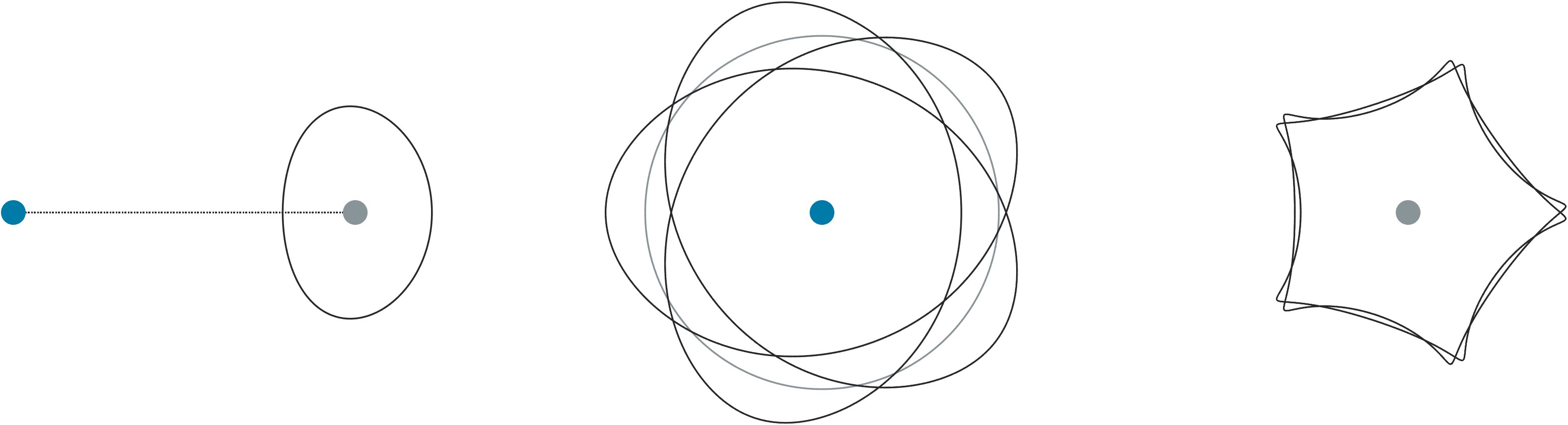}
         \caption{5:3}
     \end{subfigure}\\ % LINE BREAK!
     \begin{subfigure}{0.32\textwidth}
         \centering
         \includegraphics[width=0.99\textwidth]{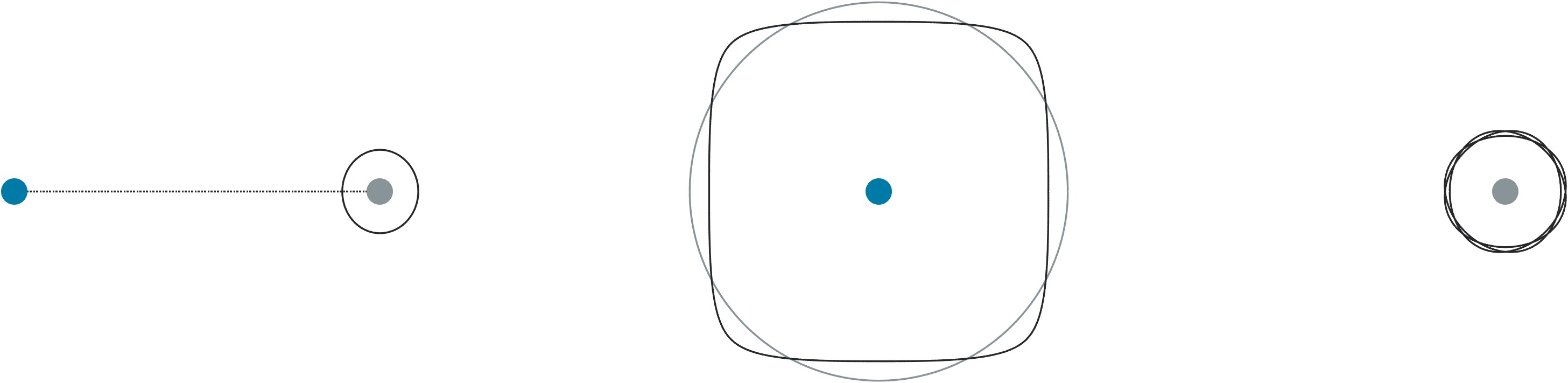}
         \caption{4:1}
     \end{subfigure}\hfill
     \begin{subfigure}{0.32\textwidth}
         \centering
         \includegraphics[width=0.99\textwidth]{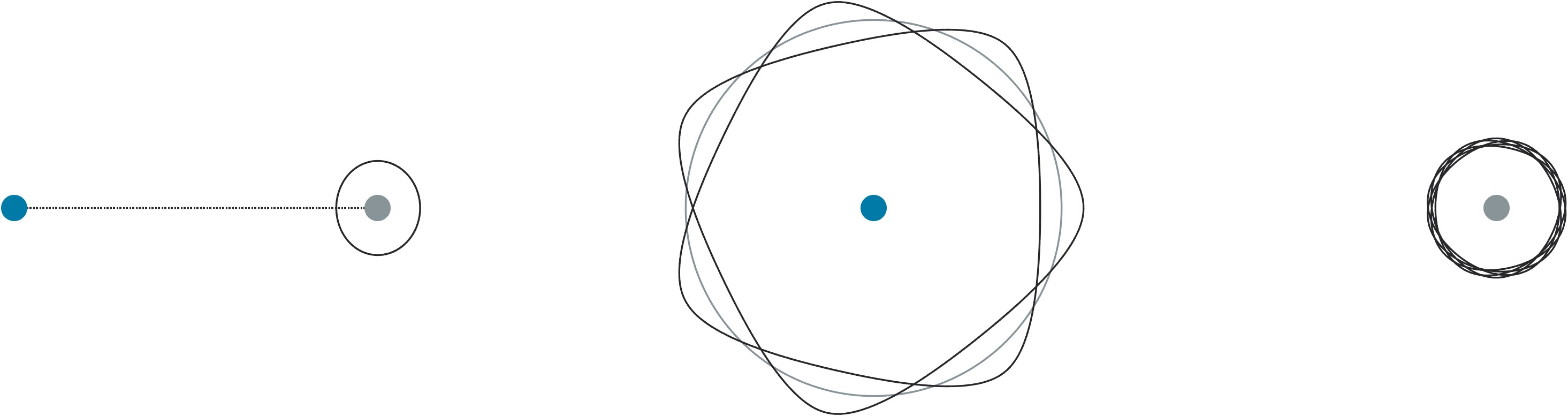}
         \caption{7:2}
     \end{subfigure}\hfill
     \begin{subfigure}{0.32\textwidth}
         \centering
         \includegraphics[width=0.99\textwidth]{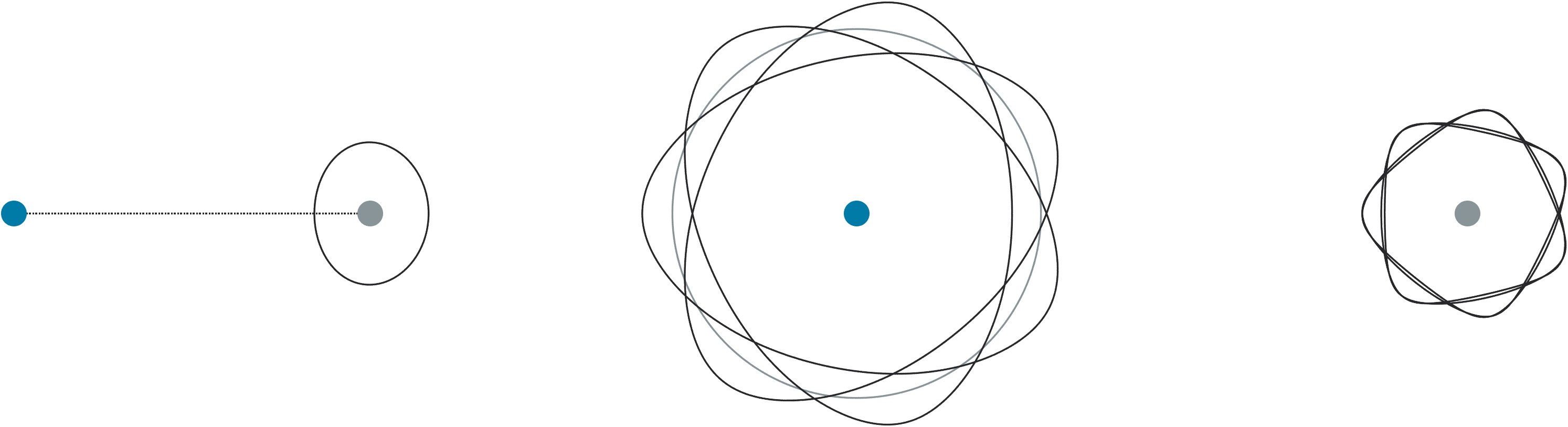}
         \caption{7:3}
     \end{subfigure}\\ % LINE BREAK!
     \begin{subfigure}{0.32\textwidth}
         \centering
         \includegraphics[width=0.99\textwidth]{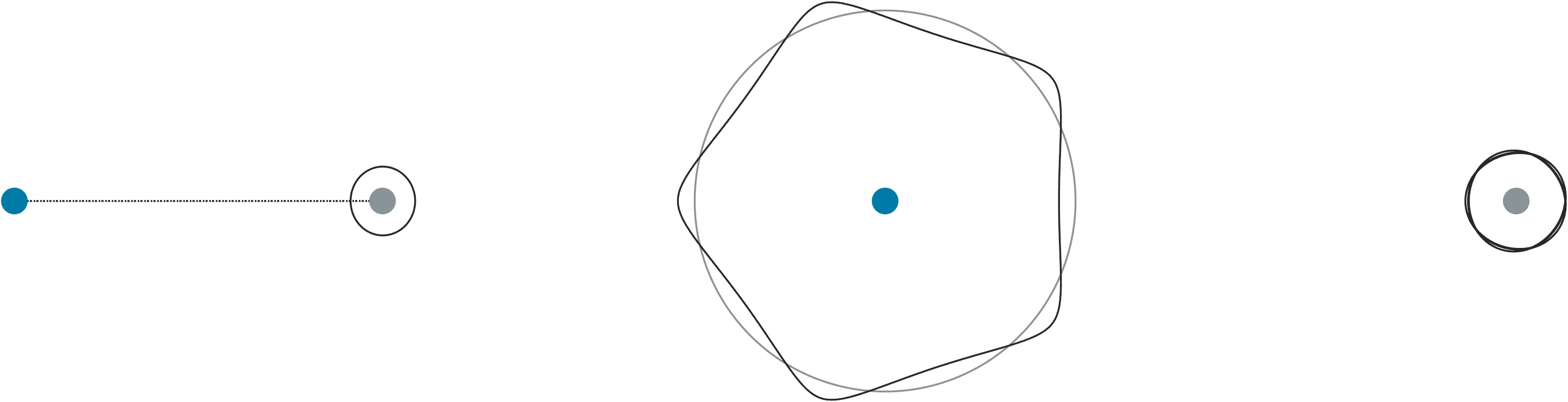}
         \caption{5:1}
     \end{subfigure}\hfill
     \begin{subfigure}{0.32\textwidth}
         \centering
         \includegraphics[width=0.99\textwidth]{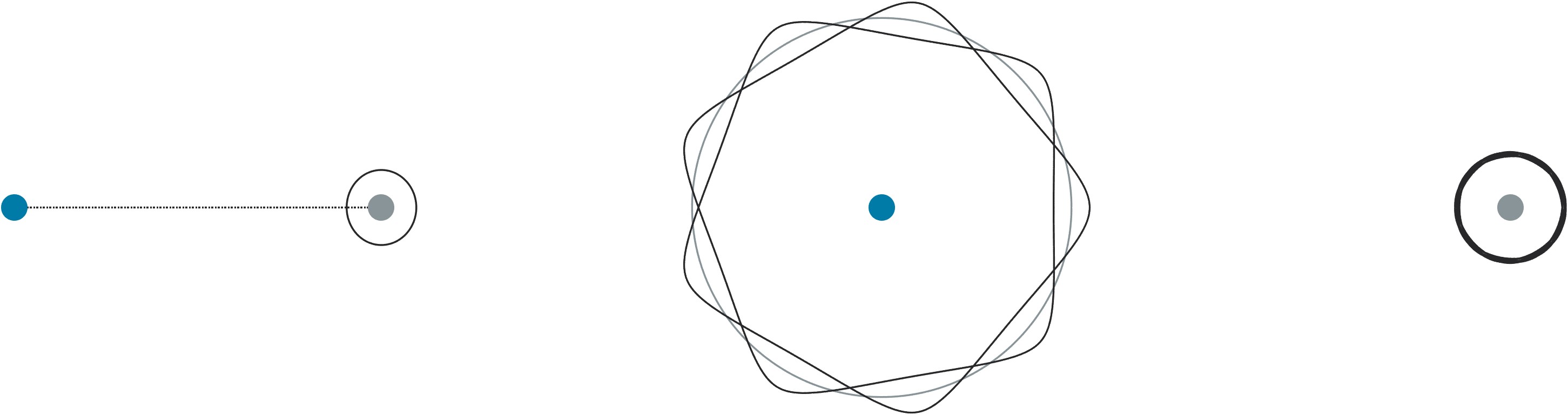}
         \caption{9:2}
     \end{subfigure}\hfill
     \begin{subfigure}{0.32\textwidth}
         \centering
         \includegraphics[width=0.99\textwidth]{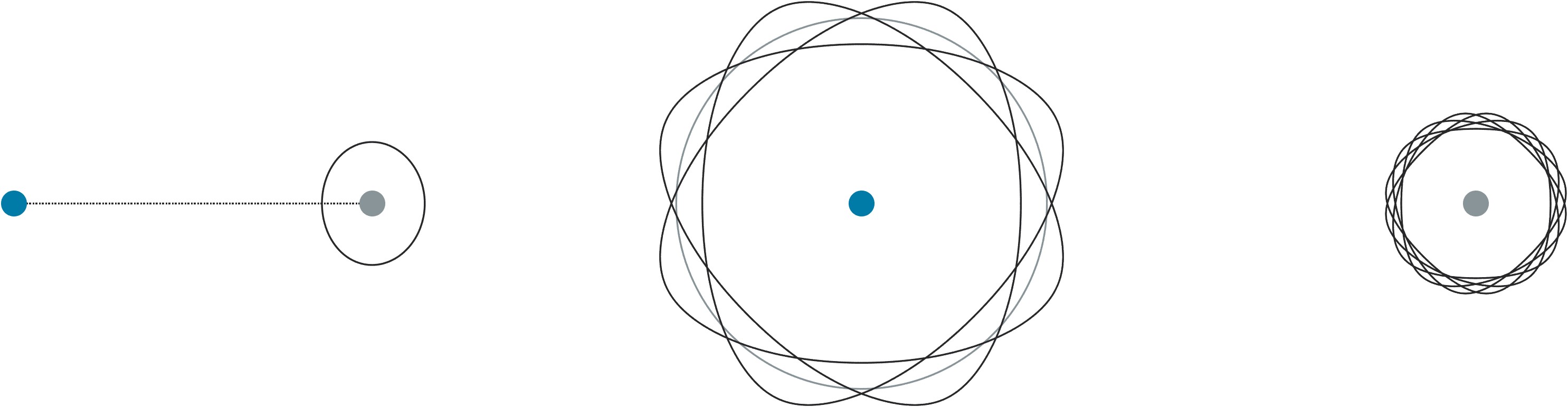}
         \caption{8:3}
     \end{subfigure}
     
     \caption{Sidereal resonant distant retrograde orbits in the Earth-Moon system. Rotating frame (left), Earth-centered inertial frame (center), and Moon-centered inertial frame (right) shown for select lower integer resonances.}
\end{figure}

\end{document}